\documentclass[manuscript,screen]{acmart}
\AtBeginDocument{%
  \providecommand\BibTeX{{%
    \normalfont B\kern-0.5em{\scshape i\kern-0.25em b}\kern-0.8em\TeX}}}

\setcopyright{acmcopyright}
\copyrightyear{2021}
\acmYear{2021}
\acmDOI{10.1145/1122445.1122456}

\acmConference[Performance '21]{Performance '21: International Symposium on Computer Performance, Modeling, Measurements and Evaluation}{Nov 08--12, 2021}{Italy}
\acmBooktitle{Performance '21, Nov 08--12, 2021, Italy}
\acmPrice{15.00}
\acmISBN{978-1-4503-XXXX-X/18/06}

\usepackage{bm}
\usepackage[figuresright]{rotating}

\usepackage{amssymb}
\usepackage{soul}
\usepackage{url}
\usepackage{graphicx}
\usepackage{subfigure}
\usepackage{amsmath}
\usepackage{booktabs}
\usepackage{bm}
\usepackage{bbding}
\usepackage{color}
\usepackage{tablefootnote}
\newtheorem{theorem}{Theorem}
\newtheorem{assumption}{Assumption}
\newtheorem{lemma}{Lemma}

\newtheorem{remark}{Remark}
\usepackage{algorithm}
\usepackage{algorithmic}
\usepackage{adjustbox}
\usepackage{amsthm,amsmath,amssymb}
\usepackage{mathrsfs}
\usepackage{lineno}
\begin{document}

\title{Simultaneously Achieving Sublinear Regret and Constraint Violations for Online Convex Optimization with Time-varying Constraints}


%
\author{Qingsong Liu} \email{liu-qs19@mails.tsinghua.edu.cn} 
\affiliation{ \institution{Institute for Interdisciplinary Information Sciences, Tsinghua University}  \country{China} }
\author{Wenfei Wu} \email{wenfeiwu@outlook.com} 
\affiliation{ \institution{School of Electronics Engineering and Computer Science, Peking University} \country{China} }
\author{Longbo Huang}\email{longbohuang@tsinghua.edu.cn}
\affiliation{ \institution{Institute for Interdisciplinary Information Sciences, Tsinghua University}  \country{China} }
\author{Zhixuan Fang}
\affiliation{ \institution{Institute for Interdisciplinary Information Sciences, Tsinghua University}  \country{China} }
\additionalaffiliation{\institution{Shanghai Qi Zhi Institute} \city{Shanghai} \country{China}}
\email{zfang@mail.tsinghua.edu.cn} 
%\affiliation{\institution{Institute for Interdisciplinary Information Sciences, Tsinghua University} \city{Beijing} \country{China}}{FIRSTAFF}
%\affiliation{\institution{School of Electronics Engineering and Computer Science, Peking University} \city{Beijing} \country{China}}{SECONDAFF}
%\affiliation{\institution{Shanghai Qi Zhi Institute} \city{Shanghai} \country{China}}{THIRDAFF}


\renewcommand{\shortauthors}{Liu and Wu, et al.}

\begin{abstract}
 %Recently, several works \cite{chen2017an,cao2018online,chen2018heterogeneous,chen2019bandit,sharma2020distributed} have studied the online convex optimization (OCO) framework with long-term and time-varying constraints. Compared with classical constrained OCO formulation, where constraint functions are time-invariant and known in advance, it can deal with broader practical applications (has broader applicability).
 %Recently, several works have studied the online convex optimization (OCO) framework with long-term and time-varying constraints. Compared with classical constrained OCO formulation, where constraint functions are time-invariant and known in advance, it can deal with broader practical applications. %(has broader applicability).
 %Online convex optimization problems with long term constraints (constrained OCO) have been investigated extensively in recent years due to their broad applicability. The classical constrained OCO formulation, where constraint functions are time-invariant and known in advance, cannot deal with problems with time-varying constraints, which arise in many practical applications. % (e.g., the schedule of time-varying uncertain power supply in smart grid).
 %This motivates us to study constrained .  
 In this paper, we develop a novel virtual-queue-based online algorithm for online convex optimization (OCO) problems with long-term and time-varying constraints
 %, which has been shown to provide improved performance bounds for time-invariant optimization in the latest works \cite{yu2020a,qiu2020beyond},
  %Both the objective function and the constraint functions are unknown a priori to the agent and revealed only after the corresponding action is chosen. 
  %Recently, virtual queues based algorithms \cite{neely2017online,yu2020a,qiu2020beyond} have been proposed and shown to provide improved performance bounds for OCO with long term constraints. However, these works consider either time-invariant constraints or static optimal benchmark (static regret) while it often behaves poorly in some inherently dynamic environment \cite{chen2017an}.
  %that the constraint functions of the OCO were time invariant which are known in advance or derived performance upper bounds with respect to static optimal benchmark(static regret) while it often behaves poorly in some inherently dynamic environment \cite{chen2017an}. 
  %These limitations made these works too restrictive and not suitable for many applications. Therefore, in this paper we are motivated to develop novel online algorithms based on virtual queues 
%for constrained OCO problems with time-varying constraint functions 
and conduct a performance analysis with respect to the dynamic regret
%benchmark sequence (dynamic regret) 
and constraint violations. %We show that our algorithm is the first parameter-free algorithm to simultaneously achieve sublinear dynamic regret and constraint violations.
We design a new update rule of dual variables and a new way of
incorporating time-varying constraint functions into the dual variables.
To the best of our knowledge, our algorithm is the first parameter-free algorithm to simultaneously achieve sublinear dynamic regret and constraint violations.
 %, as long as the temporal variations of the optimization problems are sublinear, or in the other words, 
%the objective and constraint functions do not vary too drastically across time. 
%the dynamic environment do not vary too drastically across time,
%that the accumulated variations of constraints and per-slot minimizers are sub-linearly increasing with respect to the time horizon, 
%while exiting works cannot always guarantee this.
%The performance bounds of our proposed algorithm outperform those of the state-of-the-art results in many aspects.
Our proposed algorithm also outperforms the state-of-the-art results in many aspects, e.g., our algorithm does not require the Slater condition.
%when constraints variation is not too large. 
%superior to constrained OCO method in most scenarios.
Meanwhile, for a group of practical and widely-studied constrained OCO problems in which the variation of consecutive constraints is smooth enough across time, our algorithm achieves $O(1)$ constraint violations.
Furthermore, we extend our algorithm and analysis to the case when the time horizon $T$ is unknown.
%where the time horizon $T$ is unknown.
%with an unknown time horizon $T$.
%Furthermore, we extend our algorithm and analysis to the cases when more assumptions are satisfied and the time horizon is unknown.
%, and propose a bandit version of our online algorithm which achieves the same performance bounds as the state-of-the art.
Finally, numerical experiments are conducted
to validate the theoretical guarantees of our algorithm, and some applications of our proposed framework will be outlined.
\end{abstract}




\maketitle
\section{Introduction}
%Since many practical optimization problems involved long-term time average constraints, 
%Many models have been proposed to study the optimization problems with time average constraints. Among these formulations, 
Online Convex Optimization (OCO) with long-term constraints has become one of the most popular online learning frameworks in recent years due to its powerful modeling capability for various problems
%: problems from diverse domains 
such as network routing \cite{yu2017online}, online display advertising \cite{hazan2019introduction}, and resources management \cite{chen2017an}.
In the formulation of OCO with long-term constraints, the agent wants to minimize the accumulated loss while satisfying the constraints as much as possible in the long-term. And most existing works considers the scenarios where the constraints are time-invariant \cite{yu2020a,qiu2020beyond}.

However,
%while 
%OCO problems with long-term and 
\emph{time-varying constraints} arise in many practical applications in which the underlying time-varying system is dynamic and uncertain, e.g., smart grid with uncertain renewable energy supply \cite{zhang2016peak} and data centers with dynamic user demands \cite{liu2014pricing}. Thus, this paper considers recently proposed \emph{OCO framework with long-term and time-varying constraints} \cite{cao2018online,chen2019bandit,yi2020distributed}, which is more general and practical than the one with time-variant constraints setting.
% Unlike traditional offline optimization problems, OCO is a sequential decision making procedure of an agent, who needs to commit an action at each round. The time-varying objective functions  are unknown to the agent only after an action is chosen and submitted, the objective function of the current round is revealed to the agent. Because of this lack of real time information, it is difficult for any online algorithm to find the exact optimal action in each round. Instead, related works usually use the regret as a metric to evaluate OCO algorithms, i.e., the performance gap between the actions induced by the algorithm and some benchmark actions (e.g., the static regret) and a sublinear regret is desirable since it implies that the time average performance of the algorithm is no worse than that of the benchmark asymptotically.
%\subsection{Related works}
%
%
%
%
%
%table
%\begin{center}
\begin{small}
\begin{table*}[htbp]
\centering
\resizebox{\textwidth}{30mm}{
\begin{tabular}{cccccc}
\hline
Reference   &  Regret(R) &  Constraint violations(C)
& Parameter-free& Slater condition-free\protect\footnotemark & Simultaneous \\
& & & & & sublinear R\&C \\
\hline
      \cite{chen2017an} & $O(\max\{ { V }_{ x }{ T }^{ a },$
      \\
      &  ${ V }_{ g }{ T }^{ a },{ T }^{ 1-a }\} )$ & $O({ T }^{ 1-a })$ &\Checkmark & \XSolid &\XSolid\\
      \midrule
      \cite{chen2018heterogeneous} & $O({T}^{\frac{7}{8}}{V}_{x})$ & $O(\max\{{T}^{\frac{15}{16}},{T}^{\frac{7}{8}}{V}_{x}\})$ & \Checkmark &\Checkmark & \XSolid \\
     % & &  \small{${T}^{7/8}{V}_{x}\})$} \\
      \midrule
    \cite{chen2019bandit}   & $O({ V }_{ x }{ T }^{ \frac{3}{4} })$  &  $O({ T }^{\frac{3}{4} })$ &\Checkmark & \XSolid &\XSolid
      \\
      \midrule
     \cite{chen2019bandit}   & $O({ V }_{ x }{ T }^{ \frac{1}{2} })$  &  $O({ T }^{ \frac{1}{2} })$ &\Checkmark & \XSolid &\XSolid
      \\
      \midrule
     \cite{cao2018online}  & $O({ V }_{ x }^{ \frac{1}{2} }{ T }^{ \frac{1}{2} })$  & $O({ V }_{ x }^{ \frac{1}{4} }{ T }^{ \frac{3}{4} })$ &\XSolid  & \XSolid &\Checkmark
      \\
      \midrule
    % \small{ \textcolor{red}{\cite{cao2018online}}}  & \small{$O({ V }_{ x }{ T }^{ 1/2 })$}  & \small{ $O(max\{{ T }^{ 3/4 } )$} &\small{\Checkmark}  & \small{\Checkmark}  &\small{\XSolid}
     % \\
     % &  & \small{$,{ V }_{ x }^{ 1/2 }{ T }^{ 1/2 }\}$} \\
     % \midrule
      Thm.1 & $O(\max\{ \sqrt { T{ V }_{ x } } ,{ V }_{ g }\} )$ & $O(\max\{ \sqrt { T } ,{ V }_{ g }\} )$ & \Checkmark & \Checkmark & \Checkmark
      \\
      \midrule
     Thm.1& $O( \sqrt { T{ V }_{ x } } )$ & $O(\max\{ {T}^{\frac{3}{4}} ,{ V }_{ g }\} )$ & \Checkmark & \Checkmark & \Checkmark\\
\hline
\end{tabular}
}
\label{tab:plain}
\caption{ Comparison of performance bounds for OCO with long-term and time-varying constraints w.r.t. dynamic benchmark.}
\vskip -0.3in
\end{table*}
\end{small}
%\end{center}
%\footnotetext{Slater condition-free means without requiring Slater condition holds and [3] assumes a slightly stronger Slater condition.}
\footnotetext{[3] assumes a slightly stronger Slater condition.} 
\subsection{Prior work}
%The broad applicability of OCO with long-term constraints framework motivates several pieces of work on it.
%There is not much study on OCO with long-term constraints until recent years. \cite{mahdavi2012trading} first studied the OCO with long-term and time-invariant constraints and developed an online algorithm with sublinear static regret and accumulated constraint violations. Later \cite{jenatton2016adaptive,yuan2018online} improved the performance bounds in \cite{mahdavi2012trading}. These bounds are further improved in the recent work \cite{yu2020a,wei2020online}, where a state-of-the-art static regret and constraint violations upper bounds are shown under the assumption of the Slater condition. Nevertheless, the aforementioned work either deals with time-invariant constraints or use the static benchmark to define the regret, which is not suitable for many practical scenarios. The time-invariant constraints setting cannot be applied to many practical applications in which the underlying system is dynamic and uncertain. The static benchmark assumes the optimal offline solution of the online learning problem to be constant over the whole time horizon. Hence, the static regret benchmark does not apply to the more general scenario with inherently dynamic environment, e.g., when the offline optimal is time-varying.
%\textbf{The Time-invariant Constraints.} 
OCO with long-term and time-invariant constraints has been extensively studied in the past few years. This branch of literature usually focuses on the minimization of the static regret. \cite{mahdavi2012trading} first studied the OCO with long-term and time-invariant constraints and developed an online algorithm with sublinear static regret and accumulated constraint violations. Later \cite{jenatton2016adaptive,yuan2018online} improved the performance bounds in \cite{mahdavi2012trading}. These bounds are further improved in the recent work \cite{yu2020a,qiu2020beyond}, where a state-of-the-art static regret and constraint violations upper bounds are shown under the assumption of the Slater condition. However, the setting of time-invariant constraints means the constraints will be learned by the agent easily, and hence does not capture the scenarios in which the underlying environment is dynamic and uncertain.
\textbf{The Time-varying constraints.} To overcome the limitations above, recent advances in OCO with long-term constraints considered the time-varying constraints and usually adopt a more practical but challenging metric, the dynamic regret. In this setting, a crucial challenge is to achieve sublinear dynamic regret and constraint violation simultaneously. \cite{cao2018online} studied OCO with long-term and time-varying constraints both in full-information setting and bandit setting with two-point feedback. It is the first work to simultaneously achieve sublinear dynamic regret and constraint violations. But the performance bounds attained in \cite{cao2018online} are only valid when the order of the accumulated variations of the environment is known to the agent in advance, i.e., parameter-dependent. For parameter-free work, \cite{chen2017an} analyzed the performance of a modified online saddle-point (MOSP) method and showed that sublinear dynamic regret and constraints violation may be achieved if the accumulated variations of the environment are sublinear. Later \cite{chen2018heterogeneous} improves upon it in terms of fewer assumptions but incurs a degradation of the performance. \cite{chen2019bandit} proposed a variant of MOSP method for bandit setting with two-point feedback and established the state-of-the-art performance upper bounds. However, all these parameter-free methods do not always guarantee the sublinear regret and constraint violations simultaneously, even given the accumulated variations of the environment is sublinear. Besides, most of them assume the Slater condition holds while it is not true in many scenarios. We list these works in Table 1.

Most related to our work is \cite{yu2020a} and \cite{qiu2020beyond}, which developed virtual-queue-based online algorithms and achieved the best performance bounds 
%with respect to the 
on static regret for the time-invariant constraints setting and the time-varying constraints setting, respectively. 
%These results provide an inspiring insight for OCO with long-term and time-varying constraints that when accumulated variations of the environment are sublinear,could online virtual-queue-based algorithms improve the performance upper bounds with respect to dynamic benchmark compared to the state-of-the-art or simultaneous guarantee the sublinear regret and constraint violations under only common assumptions? 
These results provide an inspiring insights for OCO with long-term constraints.
However, a challenging question remains \emph{if a virtual-queue-based algorithm can improve the state-of-the-art performance on OCO with long-term and time-varying constraints in terms of dynamic regret, and achieve sublinear regret and constraint violations simultaneously under only common assumptions.}
The answer is yes and our main contributions are summarized in the following part.
\subsection{Contributions}
We summarize our main contributions as follows.
\begin{itemize}
    \item We develop and analyze
a novel \emph{parameter-free} virtual-queue-based algorithm
for OCO with long-term and time-varying constraints. Specifically, We prove that our algorithm achieves sublinear dynamic regret and constraint violations \emph{simultaneously} without Slater condition. 
The dynamic regret and constraint violations bounds of our developed algorithm outperform the state-of-the-art %those aforementioned works 
in many aspects. 
%The bounds of the dynamic regret and constraint violations of our developed algorithm outperforms those aforementioned works in many aspects.
See also Table 1 for details.
   \item We show that when the variation of consecutive constraints is smooth enough across time, which holds in many practical applications \cite{chen2017an}, our algorithm can achieve $O(1)$ constraint violations.
   \item To the best of our knowledge, we are the first to consider the unknown time horizon case for OCO with long-term and time-varying constraints. Furthermore, our algorithm with a doubling trick can still preserve the order of performance bounds when the time horizon is unknown. 
   \item We outline some examples of applications, and fit them in the framework of OCO with long-term and time-varying constraints.
\end{itemize}
%, as long as the temporal variations of the optimization problems are sublinear. 
% Our results for both the full information setting and bandit setting and previous works are summarized in Table 2. 
%METHOD  & FEEDBACK & PARAMETER_FREE & REGRET &CONSTRAINT VIOLATIONS & SETUP
%\section{Preliminaries}
\section{Problem setup}
%In this section, we present some preliminaries and assumptions related to OCO with long-term and time-varying constraints.
In this section, we first introduce the OCO problem with long term and time-varying constraints. Then, we present the assumptions in our paper, which are widely-adopted.
\subsection{Formulation}
%Let ${ \{ f }_{ t }(x){ \}  }_{ t=1 }^{ \infty  }$ be time-varying continuous loss functions sequence defined over a closed convex set $\chi \subseteq { R }^{ n }$. For each $i\in\{1,2,...,K\}$, let ${ \{ g }_{ t,i }(x){ \}  }_{ t=1 }^{ \infty  }$ be time-varying continuous constraint functions sequence defined over a closed convex set $\chi$. 
In each round \emph{t}, the agent incurs a loss function ${f}_{t}$ and a constraint requirement ${\bm{g}}_{t}$, i.e., the agent wants to make a decision ${x}_{t}\in \chi$ to minimize the loss ${f}_{t}({x}_{t})$ while satisfying ${\bm{g}}_{t}({x}_{t})\le 0$, where ${\bm{g}}_{t}(x)$ is defined as $[{ g }_{ t,1 }(x),{ g }_{ t,2 }(x),..,{ g }_{ t,K }(x){ ] }^{ T }$. 
In this paper, we assume that ${f}_{t}(x)$ and ${g}_{t,i}(x)$ are defined over a closed convex set $\chi \subseteq { R }^{ n }$. Denote ${ \{ f }_{ t }(x){ \}  }_{ t=1 }^{ \infty  }$ and ${ \{ \bm{g} }_{ t }(x){ \}  }_{ t=1 }^{ \infty  }$ as the sequence of the time-varying loss functions and constraint functions , respectively.
%However, solving this problem is challenging in the online setting since the information about the loss and constraint functions is unknown a priori to the agent.
Thus, the agent's goal is to compute the ${x}_{t}^{*}$ defined as follows:
$$ { x }_{ t }^{ * }=\arg\min _{ x\in \chi }{\{ { f }_{ t }(x)|{ \bm{g} }_{ t }(x)\le \bm{0} \} }.$$
However, solving this problem is challenging in the online setting since the information about the loss and constraint functions is unknown a priori to the agent.
In particular, since ${\bm{g}}_{t}$ is unknown a priori, the constraint ${\bm{g}}_{t}({x}_{t})\le \bm{0}$ is hard to be satisfied in every time slot $t$. Rather, previous work \cite{chen2017an,cao2018online,chen2019bandit} allows instantaneous constraints to be violated at each round, but tries to satisfy the constraints in the
long run. In other words, the agent wants to ensure the long term constraint of 
$\sum _{ t=1 }^{ T }{ { \bm{g} }_{ t }({ x }_{ t }) } \le \bm{0}$
over some given period of
length $T$. This type of long-term constraint is appropriate in many applications (e.g., smart grid with renewable energy supply \cite{cao2018online} ).
Thus, we aim to solve the following online optimization problem.
\begin{equation}
    \tag{P1}
    \begin{split}
      \label{problem}
       \min _{ \{ { x }_{ t }{ \}  }_{ t=1 }^{ T } }{ \sum _{ t=1 }^{ T }{ { f }_{ t }({ x }_{ t }) }  },\; s.t.\;\sum _{ t=1 }^{ T }{ { \bm{g} }_{ t }({ x }_{ t }) } \le \bm{0}.
    \end{split}
\end{equation}
Solving problem \eqref{problem} exactly is still impossible in the online
setting, since the information about the ${f}_{t}$ and ${\bm{g}}_{t}$ is unknown before the action
${x}_{t}$ is chosen. Instead, our goal is to make the total loss $\sum _{ t=1 }^{ T }{ { f }_{ t }({ x }_{ t }) }$
as low as possible compared to the total loss incurred by the benchmark sequence
$\{{x}_{t}^{*}{\}}_{t=1}^{T}$ (${x}_{t}^{*}$ is commonly termed the per-slot minimizer since ${ x }_{ t }^{ * }=\arg\min _{ x\in \chi ,{ \bm{g} }_{ t }(x)\le 0 }{ { f }_{ t }(x) } $)
and meanwhile, to ensure that 
$\sum _{ t=1 }^{ T }{ { \bm{g} }_{ t }({ x }_{ t }) }$ is not too positive, i.e., the long-term constraint is not violated
too much. Therefore, 
for any sequence $\{ { x }_{ t }{ \}  }_{ t=1 }^{ T }$ yielded by online algorithms, 
we define the dynamic regret and the constraint violations, respectively as follow,
\begin{equation}
    \begin{split}
    \label{regret_constraint_violations}
        &Regret=\sum _{ t=1 }^{ T }{ { f }_{ t }({ x }_{ t }) } -\sum _{ t=1 }^{ T }{ { f }_{ t }({ x }_{ t }^{ * }) }. \\
        &{Vio}_{k}=\sum _{ t=1 }^{ T }{ { g }_{ t,k }({ x }_{ t }) },\;k\in \{ 1,2,...,K\}.
    \end{split}
\end{equation}
%, i.e., to contain both the regret (the gap between $\sum _{ t=1 }^{ T }{ { f }_{ t }({ x }_{ t }) }$ and a chosen benchmark) and constraint violation \cite{cao2018online}.
%Since the loss function and constraint function are unknown before the action is chosen, which makes it impossible to solve this problem in the online setting, related works allow instantaneous constraints to be violated at each round, but to maintain sublinear regret and constraint violations in the long-term perspective.
%Specifically, previous work \cite{chen2017an,cao2018online,chen2019bandit} allows instantaneous constraints to be violated at each round, but to maintain sublinear regret and constraint violations in the long-term perspective.
%not too positive.

In this paper, we consider the dynamic regret and constraint violations as the performance metrics. We emphasize that
the definition of dynamic regret and constraint violations in \eqref{regret_constraint_violations}
are prevalent and widely adopted in the literature \cite{chen2017an,chen2018heterogeneous,chen2019bandit,cao2018online}.
%we define the regret and the constraint violation, respectively as follow,
%A good online algorithm with this setting is to choose ${x}_{t}$ in each round $t$ which makes both dynamic regret and constraint violations grow sub-linearly with respect to time horizon $T$. 
%For any sequence $\{ { x }_{ t }{ \}  }_{ t=1 }^{ T }$ yielded by online algorithms, we define the regret and the constraint violation, respectively as follow,
%$$\sum _{ t=1 }^{ T }{ { f }_{ t }({ x }_{ t }) } -\min _{ { x }_{ t }^{ * }\in \chi ,{ \bm{g} }_{ t }({ x }_{ t }^{ * })\le \bm{0} }{ \sum _{ t=1 }^{ T }{ { f }_{ t }({ x }_{ t }^{ * }) }  }$$
%$$\sum _{ t=1 }^{ T }{ { g }_{ t,k }({ x }_{ t }) } ,\;k\in \{ 1,2,...,K\}$$
Our goal is to choose ${x}_{t}$ in each round $t$ such that both the dynamic regret and constraint violations grow sub-linearly with respect to the time horizon $T$. 
Note that the regret defined in \eqref{regret_constraint_violations} may be negative but this also makes sense. 
%It's not that the closer the total cost incurred by the agent and that of a comparator sequence is, the better our online algorithm is, since our goal is to guarantee the total cost incurred by the agent as small as possible. 
This is because we aim to minimize the total cost defined in \eqref{problem} as small as possible, while the comparator sequence can be arbitrarily given.
The significance
of the regret bound guarantee is to make sure the total cost incurred by the agent does not exceed that incurred by a comparator sequence too much, and we would like to see the appearance of the negative regret, i.e, the total cost incurred by the agent is smaller than that incurred by a comparator sequence.
Indeed, negative regret is very common in the standard OCO in terms of universal dynamic regret in which the comparator sequence is arbitrary \cite{zhao2020dynamic,zhang2020online,zhang2018adaptive}. 
%Note that the performance guarantees with respect to the dynamic regret and constraint violations of any online algorithm is strongly correlated with the temporal variations of the dynamic environment \cite{cao2018online}, i.e., time-varying objective and constraint functions, which are unknown until the corresponding actions have been chosen. 

Intuitively, the performance bounds of any online algorithm should depend on how drastically $\{{f}_{t}\}$
and $\{{\bm{g}}_{t}\}$ vary across time, that is, the temporal variations of $\{{f}_{t}\}$ and $\{{\bm{g}}_{t}\}$.
%, which also determines the drift of the per-slot minimizers sequence.
Thus we need to quantify the temporal variations of the dynamic environment. Specifically, we need to quantify the temporal variations of functions sequence. There are mainly two kinds of regularities %used wildly 
in the literature of constrained OCO \cite{chen2017an,cao2018online,chen2019bandit,yi2020distributed,sharma2020distributed,chen2018heterogeneous}.
\begin{itemize}
    \item Path-length: the accumulated variation of per-slot minimizers $\{ { x }_{ t }^{ * }\} $
    \begin{center}
         ${ V }_{ x }=\sum _{ t=2 }^{ T }{ ||{ x }_{ t }^{ * }-{ x }_{ t-1 }^{ * }|| } $
    \end{center}
    \item Function variation: the accumulatd variation of consecutive constraints
    %the variation over consecutive function values
    \begin{center}
        $ { V }_{ \bm{g} }=\sum _{ t=2 }^{ T }{ \underset { x\in \chi  }{ sup } ||{ \bm{g} }_{ t }(x)-{ \bm{g} }_{ t-1 }(x)|| } $
    \end{center}
\end{itemize}
The reason we define the accumulative variation ${V}_{x}$ with respect to ${x}_{t}^{*}$ is that it can quantify the temporal variations of the dynamic environment including loss functions ${f}_{t}$ and constraint functions ${g}_{t}$ since ${ x }_{ t }^{ * }=\arg\min _{ x\in \chi }{\{ { f }_{ t }(x)|{ \bm{g} }_{ t }(x)\le \bm{0} \} }$. While other definitions of it like $\sum _{ t=2 }^{ T }{ \max _{ x\in \chi  }{ |{ f }_{ t }-{ f }_{ t-1 }| }  } $ can only quantify the temporal variations of the loss functions.

We let $||\cdot||$ be the Euclidean norm throughout this paper.
In general, it is challenging to achieve sublinear performance bounds for any online algorithm unless regularity measures are sublinear; that is, the optimization problem is feasible. 
%On the contrary, it is hard to guarantee sublinear dynamic regret and constraint violations for any online algorithm if the accumulated variations are not sublinear. 
%The reason is that 
For example, a non-oblivious adversary may choose a new objective function ${f}_{t}$ and constraint function ${g}_{t}$ such that the current per-slot minimizer ${x}_{t}^{*}$ is at least $O(1)$ distance away from the selected action ${x}_{t}$ at each round $t$ (i.e., the accumulative variations are the of order $T$). In such case, any online algorithm cannot track the per-slot minimizers sequence $\{{x}_{t}^{*}\}$ well and guarantee the sublinear dynamic regret/constraint violations.
\subsection{Assumptions}
%We introduce some common assumptions in the literature of constrained OCO \cite{cao2018online,chen2019bandit,yu2020a}.
After specifying the problem, we introduce some assumptions in this paper, which are also common in the literature of constraint OCO \cite{cao2018online,chen2019bandit,yu2020a}. 
\begin{assumption}
\label{assumption 1}
We make following assumptions with respect to feasible set $\chi$, objective functions  ${ \{ f }_{ t }(x){ \}  }_{ t=1 }^{ T  }$ and constraint functions  ${ \{ \bm{g} }_{ t }(x){ \}  }_{ t=1 }^{ T }$:
\begin{itemize}
    \item The feasible set $\chi$ is closed, convex, and compact with diameter $R$, i.e., $\forall x,y\in \chi$, it holds that $||x-y||\le R$.
    \item The loss functions and constraint functions are convex, and bounded on $\chi$, i.e., there exists a positive constant \emph{F} such that
   $ max\{ |{ f }_{ t }(x)|,{ ||\bm{g} }_{ t }(x)||\} \le F, \forall x \in \chi, t.$
    %\item ${\bm{g}}_{t}$ is Lipschitz continuous with parameter $\beta$, i.e., $||{ \bm{g} }_{ t }(x)-{ \bm{g} }_{ t }(y)||\le \beta ||x-y||, \forall x,y \in \chi, \forall t.$ Here it implies the gradient of ${g}_{t}$ is upper-bounded by $\beta$ over $\chi$.
    \item The gradients of ${g}_{k,t}$ and ${f}_{t}$ are upper-bounded by $G$ over $\chi$, i.e., $ max\{ ||\nabla{ f }_{ t }(x)||,||\nabla{g }_{ k,t }(x)||\} \le G, \forall x \in \chi, k, t$. This is equivalent to ${\bm{g}}_{t}$ is Lipschitz continuous with parameter $\beta$ ($\beta= KG$), i.e., $||{ \bm{g} }_{ t }(x)-{ \bm{g} }_{ t }(y)||\le \beta ||x-y||, \forall x,y \in \chi, t.$
\end{itemize}
\end{assumption}
Under Assumption \ref{assumption 1}, 
%we study the standard setting in constrained OCO, namely, complete feedback or the full-information setting.
%we study two settings, namely, i) complete feedback or the full-information setting, and ii) partial feedback or the bandit setting. 
we study problem (\ref{problem}) %constrained OCO with time-varying constraints 
in the full-information setting; that is,
%, namely, i) complete feedback or the full-information setting, and ii) partial feedback or the bandit setting. 
%In the full-information setting, 
at round $t$, the agent can observe the complete loss and constraint functions after the decision ${x}_{t}$ is submitted. 
%In the bandit setting, the agent only access to the evaluation of ${f}_{t}$ and ${\bm{g}}_{t}$ at ${x}_{t}$ at round $t$. 
In the following sections, we will propose a virtual-queue-based parameter-free algorithm and show that it simultaneously achieves sublinear regret and constraint violations
%sublinear performance upper bounds 
without the Slater condition.
%, which is superior to the known results. 
%Besides, our algorithm is shown to achieve simultaneously sublinear dynamic regret and constraint violations as long as the optimization problems is feasible, i.e., the temporal variations of environment are sub-linearly increasing with respect to the time horizon. While existing works cannot always guarantee this. 
%and establish the same  performance bounds as the state-of-the-art. 
%Next we present the online algorithms for each setting and conduct the corresponding performance analysis.
\section{Algorithm}
\begin{algorithm}[tb]
\caption{VQB}
\label{algorithm1}
\begin{algorithmic}[1] %[1] enables line numbers
\STATE \textbf{Initialize}: ${\alpha}_{1},{\gamma}_{0}>0$, ${\bm{g}}_{0}=\bm{\lambda}(0)=0$, and ${x}_{1} \in \chi$.
\FOR{round $t=1...T-1$}
\STATE Update the dual iterate $\bm{\lambda}(t)$:
\STATE ${ \bm{\lambda}  }(t)=\max\{ { \bm{\lambda}  }(t-1)+{ \gamma}_{t-1}{ \bm{g} }_{ t-1 }({ x }_{ t }),-{ \gamma }_{t-1}{\bm{g} }_{ t-1 }({ x }_{ t })\} $
\STATE Update the primal iterate that satisfies:
\STATE ${ x }_{ t+1 }=\arg{ \min }_{ x\in \chi  } { \nabla { f }_{ t }({ x }_{ t }) }^{ T }(x-{ x }_{ t })+[\bm{\lambda} (t)+{\gamma}_{t-1} { \bm{g} }_{ t-1 }({ x }_{ t }){ ] }^{ T }({\gamma}_{t} { \bm{g} }_{ t }(x))+{ \alpha  }_{t}||x-{ x }_{ t }{ || }^{ 2 }$
\STATE Choose the action ${x}_{t+1}$
\ENDFOR
\end{algorithmic}
\end{algorithm}
In this section, we propose a novel virtual-queue-based algorithm, VQB, which is illustrated in Algorithm \ref{algorithm1}.
%proposed in \cite{yu2020a} to design a time-varying variant, which is illustrated in Algorithm \ref{algorithm1}.
%in the full-information setting. 
%
%doubling trick
%
%
%auxiliary
It introduces a sequence of dual variables $\{ { \bm{\lambda} }_{ t }\} $, which is also called virtual queue. The purpose of
introducing the virtual queue
%this virtual-queue strategy 
is that we can characterize the regret and constraint violations through the drift-plus-penalty expression 
%and then obtain the regret and the constraint violation bounds based on the analysis of it.
and then analysis the regret and the constraint violations based on it.
Similar ideas of updating dual variables based on virtual queues are adopted in several very recent works (e.g., \cite{yu2020a,qiu2020beyond}) for OCO with long-term and time-invariant constraints.
%Compared with previous virtual-queues-based works, e.g., \cite{yu2020a,qiu2020beyond}, 
%but the time-varying constraints bring us a unique and crucial challenge in deciding the update rule of the virtual queues.

But there are some differences between our algorithm and theirs.
First, in order to ensure both regret and constraint violations are simultaneously sublinear for the time-varying constraints setting, we design a new way of involving instantaneous per-slot constraint violation into the virtual queues and decision sequence update.
Moreover, the learning rates of our algorithm, i.e., ${\alpha}_{t}$ and ${\gamma}_{t}$ are time-varying, while the learning rates of algorithm in \cite{yu2017online,yu2020a,qiu2020beyond} are unchanged in the whole time horizon. 
%These differences make 
Therefore, our algorithm needs a new regret and constraint violation analysis due to the new update rule of virtual queues and the time-varying parameters. We will show more details in the theoretical analysis part of section $4$. 
%We find that 

Here we elaborate on the novelty and intuition of the entire algorithmic approach of VQB. 
Note that if there are no constraints ${\bm{g}}_{t}$ (i.e., ${\bm{g}}_{t}=\bm{0}$), then VQB has ${\bm{\lambda}}_{t}=\bm{0},\forall t$ and becomes the OGD algorithm, which is 
wildly-used in standard OCO with learning rate $\eta=\frac{1}{2{\alpha}_{t}}$ since
\begin{small}
\begin{equation}
\begin{split}
\label{1}
    &{ x }_{ t+1 }=\arg\min _{ x\in \chi  }{ \underbrace { { \nabla { f }_{ t }({ x }_{ t }) }^{ T }(x-{ x }_{ t })+{ \alpha  }_{ t }||x-{ x }_{ t }{ || }^{ 2 } }  } ={ \Pi  }_{ \chi  }({ x }_{ t }-\frac { 1 }{ 2{ \alpha  }_{ t } } \nabla { f }_{ t }({ x }_{ t })).\\ 
    &  \quad \quad \quad \quad \quad \quad \quad \quad \quad \quad \quad \quad penalty
    \end{split}
\end{equation}
\end{small}
Call the term marked by an underbrace in \eqref{1} the penalty. Hence, the OGD algorithm is to minimize the penalty term and is a
special case of VQB. In our algorithm VQB, if we define $\bm{Q}(t)={ { \bm{\lambda}  } }(t)+{ \gamma  }_{ t-1 }{ \bm{ g } }_{ t-1 }({ x }_{ t })=\max  \{ { \bm{ \lambda  } }(t-1)+2{ \gamma  }_{ t-1 }{ \bm{ g } }_{ t-1 }({ x }_{ t }),0\} $ to be the vector of virtual queue backlogs and define Lyapunov drift $\Delta (t)=\frac { 1 }{ 2 } ||{ \bm{Q} }(t+1){ || }^{ 2 }-\frac { 1 }{ 2 } ||{ \bm{Q} }(t){ || }^{ 2 }$, the intuition behind VQB is to choose ${x}_{t+1}$ to minimize an upper bound of the following expression (Since ${x}_{t+1}$ has not been determined at round $t$, we replace ${\bm{g}}_{t}({x}_{t+1})$ with ${\bm{g}}_{t}(x)$ in $\Delta (t)$ and omit the constant term.)
\begin{small}
\begin{align*}
    &\underbrace { \Delta (t) } +\underbrace { { \nabla { f }_{ t }({ x }_{ t }) }^{ T }(x-{ x }_{ t })+{ \alpha  }_{ t }||x-{ x }_{ t }{ || }^{ 2 } }. \\ 
    & drift \quad \quad \quad \quad \quad \quad penalty
\end{align*}
\end{small}
 Thus, the intention is to minimize penalty plus the Lyapunov drift, which is a natural method in stochastic network optimization incorporated with the stability condition (e.g., \cite{huang2011utility,huang2012lifo,huang2014power}). The drift term $\Delta(t)$ could be used to evaluate the constraint violations and is closely related to the virtual queues. The penalty term
includes the regularization term $||{x}_{t}-{x}_{t-1}{||}^{2}$ which could smoothen the difference between the coherent actions and make the whole expression strongly-convex. The remaining term describes the optimization problem.

%We find that our algorithm also has similar properties as
Our algorithm also has a close connection with the saddle point methods proposed in the literature of constrained OCO \cite{chen2017an,mahdavi2012trading}, which also incorporates dual variables to the decision-making process. For example, in Algorithm 1,  
${ \bm{\lambda}  }(t)=\max\{ { \bm{\lambda}  }(t-1)+{ \gamma}_{t-1}{ \bm{g} }_{ t-1 }({ x }_{ t }),-{ \gamma }_{t-1}{\bm{g} }_{ t-1 }({ x }_{ t })\} $
is a virtual queue vector for
the constraint violation. The role of $[\bm{\lambda} (t)+{\gamma}_{t-1} { \bm{g} }_{ t-1 }({ x }_{ t }){ ] }^{ T }$ is similar
to a dual variable vector in saddle point-typed OCO algorithms. The main differences between our algorithm and them is the update of dual variables and the way of
incorporating constraint functions into the dual variables (e.g., our algorithm uses a virtual queue to track the constraint violation, and the dual variables in our algorithm are adaptively adjusted by the per-time slot constraint violation).
These differences render our algorithm some advantages over
saddle point methods in terms of performance guarantees.
%The updating rule based on virtual queues makes our proposed algorithms superior to these saddle point methods with respect to (in terms of) performance bounds, and we will show more details in the theoretical analysis part of section $4$. 
%
%
%
%
\section{Results}
In this section, we first present the major theoretical results and analysis of our algorithm. Next, we extend our results to the case when the time horizon is unknown
%them by considering unknown time horizon case 
and the case where the variation of consecutive constraints is smooth enough across time, which captures many practical scenarios and has been frequently considered in \cite{chen2017an,xu2019efficient,amiri2019towards}.
%which is true in many practical constrained OCO problems.
\subsection{Main results}
Within this subsection, we present the upper bounds on the dynamic regret and constraint violations for VQB.
\begin{theorem}
%Consider OCO problems with long-term and time-varying constraints
Consider OCO problem \eqref{problem}
under Assumption \ref{assumption 1}, let $\{ { x }_{ t }^{ * }{ \}  }_{ t=1 }^{ T }$ be the per-slot minimizers sequence which satisfies ${ x }_{ t }^{ * }={ \text{argmin} }_{ x\in \chi ,{ \bm{g} }_{ t }(x)\le 0 }{ f }_{ t }(x)$.
\begin{itemize}
    \item (Case 1) Setting ${\alpha}_{t} =\sqrt { \frac { T }{ R+\sum _{ i\le t }{ ||{ x }_{ i }^{*}-{ x }_{ i-1 }^{*}{ || } }  }  } $ and ${ \gamma  }_{t}^{2}=\frac { 1 }{ 2{ \beta  }^{ 2 } } \frac { 1 }{ \sqrt { 2R }  }$ in VQB, then we have the following performance upper bounds
    %$(including constant factors satisfy ${\alpha}_{t} \ge 2{ \beta  }^{ 2 }{ \gamma  }_{t-1}^{ 2 } $) 
\begin{small}
\begin{equation}
    \begin{split}
        \label{case 1 of theorem 1}
      &regret \le O(\max\{ \sqrt { T{ V }_{ x } } ,{ V }_{ g }\} ),\\
      &{ Vio }_{ k }\le O(\max\{ \sqrt { T } ,{ V }_{ g }\} ),\forall k=1,2,...,K.
    \end{split}
\end{equation}
\end{small}
\item (Case 2) Setting ${\alpha}_{t} =\sqrt { \frac { T }{ R+\sum _{ i\le t }{ ||{ x }_{ i }^{*}-{ x }_{ i-1 }^{*}{ || } }  }  }$ and ${\gamma}_{t}^{2} =\frac { 1 }{ 2{ \beta  }^{ 2 } } \frac { 1 }{ \sqrt { 2R }  }\frac{1}{\sqrt{t+1}}$ in VQB, then we have the following performance upper bounds
\begin{small}
\begin{equation}
    \begin{split}
    \label{case 2 of theorem 1}
      &regret \le O( \sqrt { T{ V }_{ x } } ),\\
      &{ Vio }_{ k }\le O(\max\{ {T}^{\frac{3}{4}},{ V }_{ g }\} ),\forall k=1,2,...,K.
    \end{split}
\end{equation}
\end{small}
\end{itemize}
\end{theorem}

%Furthermore, 
There are several advantages stated as following that makes our results outperform previous studies. First, Theorem 1 implies that 
%under only common assumptions 
VQB can guarantee sublinear regret and constraint violations simultaneously, as long as the accumulated variations of the environment are sublinear, i.e., ${V}_{x}=o(T)$ and ${V}_{\bm{g}}=o(T)$.
%sub-linearly increasing with respect to the time horizon. 
%In addition to requiring stronger conditions(e.g., Slater condition),
Previous studies listed in Table 1 do not always simultaneously guarantee the sublinear performance bounds  since they introduce the $o(T){V}_{x}$ or $o(T){V}_{\bm{g}}$ term in their performance bounds, which may be at least of the order \emph{T} even the optimization problem is feasible, i.e., $\max\{O({V}_{x}), O({V}_{\bm{g}})\}=o(T)$. 

Second, the dynamic regret upper bound guaranteed by both two cases of Theorem 1 could match the state-of-the-art dynamic regret bound $O(\sqrt{T{V}_{x}})$ in general OCO \cite{zhang2018adaptive,zhao2018proximal,zhao2020dynamic}, when the path-length of the benchmark sequence is ${V}_{x}$.

Moreover, 
%Theorem 1 holds no matter whether Slater condition holds or not and 
our algorithm is parameter-free, that is, the parameters in our algorithm do not require prior information of the
%have the dependency on any 
regularities (e.g., ${V}_{x}$ or ${V}_{\bm{g}}$).
Meanwhile, Theorem 1 holds no matter whether the Slater condition holds or not.
The theoretical results of %previous studies
most previous study are valid either under the Slater condition, or the order of the regularities are known prior to the learner.
Only \cite{chen2018heterogeneous} is both parameter-free and independent of the assumption of the Slater condition, however, it 
introduced degraded performance bounds and cannot guarantee the sublinear regret and constraint violations simultaneously. 
Readers could see Table 1 for the detailed comparisons.

We compare the performance bounds of our algorithm with the previous studies listed in Table 1. When ${V}_{x}$ is not too large (e.g., ${V}_{x}=o(\sqrt { T } )$), the regret and constraint violations bounds presented in the first case of Theorem 1 are all no worse than the state-of-the-art results, i.e., $O[\sqrt { T } { V }_{ x },\sqrt { T } ]$ and $O[\sqrt{T{V}_{x}},{V}_{x}^{1/4}{T}^{3/4}]$, established in \cite{chen2019bandit} and \cite{cao2018online}, respectively. Besides, the dynamic regret bound presented in the second case of Theorem 1 is superior to all existing works, and the corresponding constraint violations are also strictly sublinear when the optimization problem is feasible.
%Besides the action sequence $\{ { x }_{ t }\} $ involved in our algorithm, it also

~\\
\textbf{Proof sketch of Theorem 1.}
Within this subsection, we give a proof sketch of Theorem 1. All the proof details of listed lemmas could be found in the Appendix.
Since the drift-plus-penalty expression characterizes the dynamic regret expression, we can translate the bounds of virtual queues $\{ \bm{\lambda} (t)\} $ into bounds of constraint violations. Thus, our proof starts with the analysis of virtual queues properties and drift-plus-penalty expression, that is, the Lyapunov drift term $\Delta(t)=\frac { 1 }{ 2 } [||\bm{\lambda (t+1)}{ || }^{ 2 }-||{ \bm{\lambda}  }(t){ || }^{ 2 }]$ plus the penalty term ${f}_{t}({x}_{t})$, which is associated with the loss value after choosing an action. First, we present the main properties for virtual queues $\{{\lambda}_{t}{\}}^{T}_{t=1}$ introduced in Algorithm 1 and Lyapunov drift term.
\begin{lemma}
\label{lemma 0}
(Properties of virtual queues)
In Algorithm 1, we have the following properties for virtual queues ${\bm{\lambda}}(t)$ and Lyapunov drift term $\Delta(t)$:
\begin{enumerate}%[(i)]
    \item ${\bm{\lambda}}(t)\ge \bm{0}$
    \item ${ \bm{\lambda}  }(t)+{ {\gamma}_{t-1} \bm{g} }_{ t-1 }({ x }_{ t })\ge 0$ 
    \item ${ ||\bm{\lambda}  }(t){ || }\ge { {\gamma}_{t-1} ||\bm{g} }_{ t-1 }({ x }_{ t }){ || }$
    \item ${ {\gamma}_{t-1} \bm{g} }_{ t-1}({ x }_{ t })\le { \bm{\lambda}  }(t)-{ \bm{\lambda}  }(t-1)$, furthermore, ${ ||\bm{\lambda}  }(t){ || }-{ ||\bm{\lambda}  }(t-1){ || }{ \le {\gamma}_{t-1} ||\bm{g} }_{ t-1 }({ x }_{ t }){ || }$ 
    \item $\Delta (t)\le {\gamma}_{t} [\bm{\lambda} (t){ ] }^{ T }{ \bm{g} }_{ t }({ x }_{ t+1 })+{ \gamma  }_{t}^{ 2 }||{ \bm{g} }_{ t }({ x }_{ t+1 }){ || }^{ 2 }$ 
\end{enumerate}
\end{lemma}
The proof of this lemma is motivated by \cite{yu2020a,qiu2020beyond}. However, due to our new algorithm, different constraints setting and fewer assumptions, our proof techniques is slightly different from theirs. Then we present the upper bound of drift-plus-penalty expression in the following lemma.
%constraint violations could be obtained from bounding the virtual queues, 
\begin{lemma}
\label{lemma 1}
(Upper bound of the drift-plus-penalty expression)
Under Assumption \ref{assumption 1}, let $\delta>0$ and $\{ { \alpha  }_{ t }{ \}  }_{ t=1 }^{ T }$, $\{ { \gamma  }_{ t }{ \}  }_{ t=1 }^{ T }$ be any positive non-increasing sequences, if $2{ \gamma  }_{ t }\le { \gamma  }_{ t-1 }+{ \gamma  }_{ t+1 }$ holds for all $t$, then VQB ensures that:
\begin{equation}
    \begin{split}
       &{ f }_{ t }({ x }_{ t })+\Delta (t)\le { \alpha  }_{ t }||{ x }_{ t }^{ * }-{ x }_{ t }{ || }^{ 2 }-{ \alpha  }_{ t+1 }||{ x }_{ t+1 }{ -{ x }_{ t+1 }^{ * }|| }^{ 2 }
       +4R{ \alpha  }_{ t }{ ||{ x }_{ t+1 }^{ * }-{ x }_{ t }^{ * }|| }+({ \beta  }^{ 2 }{ \gamma  }_{t-1}^{ 2 }+\frac { \delta  }{ 2 } -{ \alpha  }_{ t })||{ x }_{ t+1 }-{ x }_{ t }{ || }^{ 2 }\\
       &+\frac { 1 }{ 2\delta  } { G }^{ 2 } +\frac { 1 }{ 2 } { \gamma  }_{t}{\gamma}_{t+1}||{ { { \bm{g} } } }_{ t }({ x }_{ t+1 }){ || }^{ 2 }-\frac{1}{2}{ \gamma  }_{t-1}{\gamma}_{t}||{ { { \bm{g} } } }_{ t-1 }({ x }_{ t }){ || }^{ 2 } 
       +{ \gamma  }_{t-1}{\gamma}_{t}||{ { { \bm{g} } } }_{ t-1 }({ x }_{ t })-{ { { \bm{g} } } }_{ t }({ x }_{ t }){ || }^{ 2 } + { f }_{ t }({ x }_{ t }^{ * }).
    \end{split}
\end{equation}
\end{lemma}
This is the key lemma in our theoretical analysis, which is used to yield the eventual bounds of regret and virtual queues. Next, we bound the dynamic regret as follows based on Lemma \ref{lemma 1}.
\begin{lemma}
\label{lemma 2}
(Regret bound) Under Assumption \ref{assumption 1}, for arbitrary $\delta>0$ which satisfies ${\alpha}_{t}\ge { \beta  }^{ 2 }{ \gamma  }_{t-1}^{ 2 }+\frac { \delta  }{ 2 }$, if ${\gamma}_{t}\le{\gamma}_{t+1}$, ${\alpha}_{t}\le {\alpha}_{t+1}$ and $2{ \gamma  }_{ t }\le { \gamma  }_{ t-1 }+{ \gamma  }_{ t+1 }$ hold for all $t$, then VQB ensures that
%in Algorithm \ref{algorithm1}, then we have:
\begin{equation}
    \begin{split}
        \sum _{ t=1 }^{ T }{ { f }_{ t }({ x }_{ t }) } &\le \sum _{ t=1 }^{ T }{ { f }_{ t }({ x }_{ t }^{ * }) } +{ \alpha  }_{ 1 }{R}^{2}+4R\sum _{ t=1 }^{ T }{ { \alpha  }_{ t }{ ||{ x }_{ t+1 }^{ * }-{ x }_{ t }^{ * }|| } } 
        +\frac { T{ G }^{ 2 } }{ 2\delta  } 
        +\frac { 1 }{ 2 } { \gamma  }_{ T}{ \gamma  }_{ T+1 }||{ { \bm{ g } } }_{ T }({ x }_{ T+1 }){ || }^{ 2 } \\
        &+\frac { 1 }{ 2 } ||{ \bm{\lambda}  }(1){ || }^{ 2 }
        +2F\sum _{ t=1 }^{ T }{ { { \gamma  }_{ t-1 }^{ 2 }||{ \bm{g} }_{ t-1 }({ x }_{ t })-{ \bm{g} }_{ t }({ x }_{ t })|| } }.
        %&\quad \quad \quad \sum _{ t=1 }^{ T }{ { g }_{ k,t }({ x }_{ t }) } \le \frac { { \lambda  }_{k}(T) }{ {\gamma}_{T}  }+{ V }_{ \bm{g} }(T),\forall k=1,2,...,K.
    \end{split}
\end{equation}
Here we define ${x}_{T+1}^{*}={x}_{T}^{*}$.
\end{lemma}
%Please refer to the supplementary material for the detailed proof. 
%In order to obtain the constraint violations upper bounds, 
We further bound the eventual virtual queues length in the following lemma based on Lemma \ref{lemma 1}.
\begin{lemma}
\label{lemma 3}
Under Assumption 1, setting $\delta,{\alpha}_{t} \;and \;{\gamma}_{t}$ to be the same as Lemma \ref{lemma 2}, then  VQB ensures that
\begin{equation}
    \begin{split}
        ||\bm{\lambda} (T){ || }&\le 2\sqrt { F(T-1) } +\sqrt { 2{ \alpha  }_{ 1 }{R}^{2} } +\sqrt { \frac { (T-1){ G }^{ 2 } }{ \delta  }  } 
        +{ \gamma  }_{ T-1 }||{ { { \bm{ g } } } }_{ T-1 }({ x }_{ T }){ || }\\
        &+2\sqrt { 2R\sum _{ t=1 }^{ T-1 }{ { \alpha  }_{ t }||{ x }_{ t+1 }^{ * }-{ x }_{ t }^{ * }{ || } }  }
        +2\sqrt { F\sum _{ t=1 }^{ T-1 }{ { \gamma  }_{ t-1 }^{ 2 }||{ { { \bm{ g } } } }_{ t-1 }({ x }_{ t })-{ { { \bm{ g } } } }_{ t }({ x }_{ t }){ || } }  }.
    \end{split}
\end{equation}
\end{lemma}
%Please refer to the supplementary material for the  detailed proof. 
This is another critical lemma in our theoretical analysis that could be used to yield the constraint violations' upper bounds. We upper bound the constraint violations in the following two lemmas.
\begin{lemma}
\label{lemma 4}
For any non-increasing sequence $\{{\gamma}_{t}\}$, VQB ensures that
\begin{equation}
    \sum _{ t=1 }^{ T }{ { g }_{ k,t }({ x }_{ t }) } \le \frac {|| { \bm{\lambda}  }(T)|| }{ {\gamma}_{T}  }+{ V }_{ \bm{g} },\;\forall k=1,2,...,K.
\end{equation}
\end{lemma}
%Please refer to the supplementary material for the  detailed proof.
Recall that Lemma \ref{lemma 3} bounds the virtual queue length. Thus combining this lemma with the Lemma \ref{lemma 4}, we can bound the constraint violations in the following lemma.
\begin{lemma}
\label{lemma 5}
(Constraint violations' bounds) Setting $\delta,{\alpha}_{t}$ and ${\gamma}_{t}$ to be the same as Lemma \ref{lemma 2}, then VQB ensures that
\begin{equation}
    \begin{split}
        \sum _{ t=1 }^{ T }{ { g }_{ k,t }({ x }_{ t }) } &\le \frac { 2 }{ { \gamma  }_{ T } } \sqrt { F(T-1) } +\frac { 1 }{ { \gamma  }_{ T } } \sqrt { 2{ \alpha  }_{ 1 }{R}^{2} }+{ V }_{ { \bm{g} } } 
        +\frac { { \gamma  }_{ T-1 } }{ { \gamma  }_{ T } } ||{ { { \bm{ g } } } }_{ T-1 }({ x }_{ T }){ || } \\
        &+\frac { 2 }{ { \gamma  }_{ T } } \sqrt { 2R\sum _{ t=1 }^{ T-1 }{ { \alpha  }_{ t }||{ x }_{ t+1 }^{ * }-{ x }_{ t }^{ * }{ || } }  } +\frac { 2 }{ { \gamma  }_{ T } } \sqrt { F\sum _{ t=1 }^{ T-1 }{ { \gamma  }_{ t-1 }^{ 2 }||{ { { \bm{ g } } } }_{ t-1 }({ x }_{ t })-{ { { \bm{ g } } } }_{ t }({ x }_{ t }){ || } }  } +\frac { G }{ { \gamma  }_{ T } } \sqrt { \frac { T-1 }{ \delta  }  }.
    \end{split}
\end{equation}
\end{lemma}
%Please refer to the supplementary material for the detailed proof.
According to Lemmas \ref{lemma 2} and \ref{lemma 5}, with parameters stated in Theorem 1, we could prove the theoretical results of Theorem 1. First we consider the \textbf{Case 1} in Theorem 1, by the setting of ${\alpha}_{t}$ and according to Lemma \ref{pre_lemma 2}, we can obtain
\begin{equation}
    \begin{split}
      \label{1_proof_of_theorem_1}
        &\sum _{ t=1 }^{ T }{ { \alpha  }_{ t }{ ||{ x }_{ t+1 }^{ * }-{ x }_{ t }^{ * }|| } }  
        =\sum _{ t=1 }^{ T }{ \sqrt { \frac { T }{ R+\sum _{ i\le t }{ ||{ x }_{ i }-{ x }_{ i-1 }{ || } }  }  } ||{ x }_{ t+1 }^{ * }-{ x }_{ t }^{ * }|| } \\ 
        &=\sqrt { T } \sum _{ t=1 }^{ T }{ \frac { ||{ x }_{ t+1 }^{ * }-{ x }_{ t }^{ * }|| }{ \sqrt { R+\sum _{ i\le t }{ ||{ x }_{ i }-{ x }_{ i-1 }{ || } }  }  }  }
        \le 2\sqrt { T } \sqrt { \sum _{ t=0 }^{ T }{ ||{ x }_{ t+1 }^{ * }-{ x }_{ t }^{ * }|| }  }  =2\sqrt { T{ V }_{ x } } 
    \end{split}
\end{equation}
By the setting of ${\gamma}_{t}$ and according to Lemma \ref{pre_lemma 1}, we also have
\begin{equation}
    \begin{split}
    \label{2_proof_of_theorem_1}
        &\sum _{ t=1 }^{ T }{ { \gamma  }_{ t-1 }^{ 2 }||{ { \bm{ g } } }_{ t-1 }({ x }_{ t })-{ { \bm{ g } } }_{ t }({ x }_{ t }){ || } } 
        =\frac { 1 }{ 2{ \beta  }^{ 2 } } \frac { 1 }{ \sqrt { 2R }  } \sum _{ t=1 }^{ T }{ ||{ { \bm{ g } } }_{ t-1 }({ x }_{ t })-{ { \bm{ g } } }_{ t }({ x }_{ t }){ || } } \\ 
        &\le \frac { 1 }{ 2{ \beta  }^{ 2 } } \frac { 1 }{ \sqrt { 2R }  } \sum _{ t=1 }^{ T }{ \max _{ x\in \chi  }{ ||{ { \bm{ g } } }_{ t-1 }({ x })-{ { \bm{ g } } }_{ t }({ x }){ || } }  } =\frac { 1 }{ 2{ \beta  }^{ 2 } } \frac { 1 }{ \sqrt { 2R }  } { V }_{ \bm{g} }
    \end{split}
\end{equation}
Setting $\delta=\frac { 1 }{ 2 } \sqrt { \frac { T }{ R+{ V }_{ x } }  }$, it is easy to verity that ${\alpha}_{t}\ge { \beta  }^{ 2 }{ \gamma  }_{t-1}^{ 2 }+\frac { \delta  }{ 2 }$, $2{ \gamma  }_{ t }\le { \gamma  }_{ t-1 }+{ \gamma  }_{ t+1 }$, and both $\{{\alpha}_{t}\}$ and $\{{\gamma}_{t}\}$ are non-increasing sequences.
Thus combing Lemma 2 with (\ref{1_proof_of_theorem_1}), (\ref{2_proof_of_theorem_1}) and rearranging terms yields
\begin{equation}
    \begin{split}
        &\sum _{ t=1 }^{ T }{ { f }_{ t }({ x }_{ t }) } -\sum _{ t=1 }^{ T }{ { f }_{ t }({ x }_{ t }^{ * }) } \\
        &\le 4R\sum _{ t=1 }^{ T }{ { \alpha  }_{ t }{ ||{ x }_{ t+1 }^{ * }-{ x }_{ t }^{ * }|| } } +{ \alpha  }_{ 1 }{ R }^{ 2 }+\frac { T{ G }^{ 2 } }{ 2\delta  } +\frac { 1 }{ 2 } { \gamma  }_{ T }{ \gamma  }_{ T+1 }||{ { \bm{ g } } }_{ T }({ x }_{ T+1 }){ || }^{ 2 }
        +2F\sum _{ t=1 }^{ T }{ { \gamma  }_{ t-1 }^{ 2 }||{ { \bm{ g } } }_{ t-1 }({ x }_{ t })-{ { \bm{ g } } }_{ t }({ x }_{ t }){ || } } +\frac { 1 }{ 2 } ||\bm{\lambda} (1){ || }^{ 2 }\\ 
        &\le 8R\sqrt { T{ V }_{ x } } +\sqrt { \frac { T }{ R }  } { R }^{ 2 }+\frac { T{ G }^{ 2 } }{ 2\delta  } +\frac { 1 }{ 4{ \beta  }^{ 2 } } \frac { 1 }{ \sqrt { 2R }  } { F }^{ 2 }+\frac { F }{ { \beta  }^{ 2 } } \frac { { V }_{ \bm{g} } }{ \sqrt { 2R }  } +\frac { 1 }{ 2 } ||\bm{\lambda} (1){ || }^{ 2 }\\ &\overset { (a) }{ \le  } 8R\sqrt { T{ V }_{ x } } +\sqrt { \frac { T }{ R }  } { R }^{ 2 }+T{ G }^{ 2 }\sqrt { \frac { R+{ V }_{ x } }{ T }  } +\frac { 1 }{ 4{ \beta  }^{ 2 } } \frac { 1 }{ \sqrt { 2R }  } { F }^{ 2 }+\frac { F }{ { \beta  }^{ 2 } } \frac { { V }_{ \bm{g} } }{ \sqrt { 2R }  } +\frac { 1 }{ 2 } ||\bm{\lambda} (1){ || }^{ 2 }\\ 
        &=8R\sqrt { T{ V }_{ x } } +{ R }^{ 2/3 }\sqrt { T } +{ G }^{ 2 }\sqrt { T(R+{ V }_{ x }) } +\frac { 1 }{ 4{ \beta  }^{ 2 } } \frac { 1 }{ \sqrt { 2R }  } { F }^{ 2 }+\frac { F }{ { \beta  }^{ 2 } } \frac { { V }_{ \bm{g} } }{ \sqrt { 2R }  } +\frac { 1 }{ 2 } ||\bm{\lambda} (1){ || }^{ 2 }\\ 
        &=O(max\{ \sqrt { T{ V }_{ x } } ,{ V }_{ \bm{g} }\} )
    \end{split}
\end{equation}
Where (a) holds since we set $\delta=\frac { 1 }{ 2 } \sqrt { \frac { T }{ R+{ V }_{ x } }  }$. According to Lemma $5$ and Assumption 1, we have
\begin{align*}
    &\sum _{ t=1 }^{ T }{ { g }_{ t,k }({ x }_{ t }) } \le \frac { 2 }{ { \gamma  }_{ T } } \sqrt { F(T-1) } +\frac { 1 }{ { \gamma  }_{ T } } \sqrt { 2{ \alpha  }_{ 1 }{ R }^{ 2 } } +\frac { { \gamma  }_{ T-1 } }{ { \gamma  }_{ T } } ||{ { {\bm { g } } } }_{ T-1 }({ x }_{ T }){ || }+\frac { G }{ { \gamma  }_{ T } } \sqrt { \frac { T-1 }{ \delta  }  } \\
        &+\frac { 2 }{ { \gamma  }_{ T } } \sqrt { 2R\sum _{ t=1 }^{ T-1 }{ { \alpha  }_{ t }||{ x }_{ t+1 }^{ * }-{ x }_{ t }^{ * }{ || } }  } +\frac { 2 }{ { \gamma  }_{ T } } \sqrt { F\sum _{ t=1 }^{ T-1 }{ { \gamma  }_{ t-1 }^{ 2 }||{ { { \bm{ g } } } }_{ t-1 }({ x }_{ t })-{ { {\bm { g } } } }_{ t }({ x }_{ t }){ || } }  } +{ V }_{ \bm{ g } }\\ 
        &\overset { (a) }{ \le  } \sqrt { \frac { 1 }{ 2\sqrt { 2R } { \beta  }^{ 2 } }  } [2\sqrt { F(T-1) } +R\sqrt { 2\sqrt { \frac { T }{ R }  }  } +G\sqrt { 2\frac { (T-1)\sqrt { R+{ V }_{ x } }  }{ \sqrt { T }  }  } ]
\end{align*}
\begin{equation}
    \begin{split}
        &+F+\sqrt { \frac { 2 }{ \sqrt { 2R } { \beta  }^{ 2 } }  } [\sqrt { 4R\sqrt { T{ V }_{ x } }  } +\sqrt { \frac { F }{ 2{ \beta  }^{ 2 } } \frac { { V }_{ \bm{g} } }{ \sqrt { 2R }  }  } ]+{ V }_{ \bm{ g } }\\ 
        &\overset { (b) }{ \le  } \sqrt { \frac { 1 }{ 2\sqrt { 2R } { \beta  }^{ 2 } }  } [2\sqrt { F(T-1) } +R\sqrt { 2\sqrt { \frac { T }{ R }  }  } +G\sqrt { 2R\sqrt { T(1+T } ) } ]\\ &+F+\sqrt { \frac { 2 }{ \sqrt { 2R } { \beta  }^{ 2 } }  } [\sqrt { 4R\sqrt { R } T } +\sqrt { \frac { F }{ 2{ \beta  }^{ 2 } } \frac { { V }_{ \bm{g} } }{ \sqrt { 2R }  }  } ]+{ V }_{ \bm{ g } }
        =O(max\{ \sqrt { T } ,{ V }_{ \bm{g} }\} )
    \end{split}
\end{equation}
Where (a) is due to (\ref{1_proof_of_theorem_1}) and (\ref{2_proof_of_theorem_1}); (b) is due to ${V}_{x}\le RT$. For the
\textbf{Case 2} in Theorem 1, since the settings of ${\alpha}_{t}$ in both two cases are the same, we can also derive that
\begin{equation}
\label{3_proof_of_theorem_1}
    \sum _{ t=1 }^{ T }{ { \alpha  }_{ t }{ ||{ x }_{ t+1 }^{ * }-{ x }_{ t }^{ * }|| } }\le 2\sqrt { T{ V }_{ x } } 
\end{equation}
For term $\sum _{ t=1 }^{ T }{ { \gamma  }_{ t-1 }^{ 2 }||{ {\bm { g } } }_{ t-1 }({ x }_{ t })-{ {\bm { g } } }_{ t }({ x }_{ t }){ || } } $, according to Lemma \ref{pre_lemma 1} and by the setting of ${\gamma}_{t}$, we have
\begin{equation}
    \begin{split}
    \label{4_proof_theorem_1}
        &\sum _{ t=1 }^{ T }{ { \gamma  }_{ t-1 }^{ 2 }||{ { \bm{ g } } }_{ t-1 }({ x }_{ t })-{ { \bm{ g } } }_{ t }({ x }_{ t }){ || } } =\frac { 1 }{ 2{ \beta  }^{ 2 } } \frac { 1 }{ \sqrt { 2R }  } \sum _{ t=1 }^{ T }{ \frac { 1 }{ \sqrt { t }  } ||{ { \bm{ g } } }_{ t-1 }({ x }_{ t })-{ { \bm{ g } } }_{ t }({ x }_{ t }){ || } } \\
        &\le \frac { F }{ { \beta  }^{ 2 } } \frac { 1 }{ \sqrt { 2R }  } \sum _{ t=1 }^{ T }{ \frac { 1 }{ \sqrt { t }  }  } \le \frac { F }{ { \beta  }^{ 2 } } \frac { \sqrt { 2T }  }{ \sqrt { R }  } 
    \end{split}
\end{equation}
Setting $\delta=\frac { 1 }{ 2 } \sqrt { \frac { T }{ R+{ V }_{ x } }  }$, it is easy to verity that ${\alpha}_{t}\ge { \beta  }^{ 2 }{ \gamma  }_{t-1}^{ 2 }+\frac { \delta  }{ 2 }$, $2{ \gamma  }_{ t }\le { \gamma  }_{ t-1 }+{ \gamma  }_{ t+1 }$, and both $\{{\alpha}_{t}\}$ and $\{{\gamma}_{t}\}$ are non-increasing sequences. Based on Lemma 2, combining (\ref{3_proof_of_theorem_1}), (\ref{4_proof_theorem_1}) and Assumption 1 gives
\begin{equation}
    \begin{split}
        &\sum _{ t=1 }^{ T }{ { f }_{ t }({ x }_{ t }) } -\sum _{ t=1 }^{ T }{ { f }_{ t }({ x }_{ t }^{ * }) } \\
        &\le 4R\sum _{ t=1 }^{ T }{ { \alpha  }_{ t }{ ||{ x }_{ t+1 }^{ * }-{ x }_{ t }^{ * }|| } } +{ \alpha  }_{ 1 }{ R }^{ 2 }+\frac { T{ G }^{ 2 } }{ 2\delta  } +\frac { 1 }{ 2 } { \gamma  }_{ T }{ \gamma  }_{ T+1 }||{ { \bm{ g } } }_{ T }({ x }_{ T+1 }){ || }^{ 2 }
        +2F\sum _{ t=1 }^{ T }{ { \gamma  }_{ t-1 }^{ 2 }||{ { \bm{ g } } }_{ t-1 }({ x }_{ t })-{ { \bm{ g } } }_{ t }({ x }_{ t }){ || } } +\frac { 1 }{ 2 } ||\bm{\lambda} (1){ || }^{ 2 }\\ 
        &\le 8R\sqrt { T{ V }_{ x } } +{ R }^{ 3/2 }\sqrt { T } +T{ G }^{ 2 }\sqrt { \frac { R+{ V }_{ x } }{ T }  } +\frac { { F }^{ 2 } }{ 4{ \beta  }^{ 2 } } \frac { 1 }{ \sqrt { 2R }  } \frac { 1 }{ \sqrt { T }  } +\frac { 2{ F }^{ 2 } }{ { \beta  }^{ 2 } } \frac { \sqrt { 2T }  }{ \sqrt { R }  } +\frac { 1 }{ 2 } ||\bm{\lambda} (1){ || }^{ 2 }
        =O(\sqrt { T{ V }_{ x } } )
    \end{split}
\end{equation}
Furthermore, based on Lemma 5, we obtain the bounds of constraint violations as follows
\begin{align*}
    &\sum _{ t=1 }^{ T }{ { g }_{ t,k }({ x }_{ t }) } \le \frac { 2 }{ { \gamma  }_{ T } } \sqrt { F(T-1) } +\frac { 1 }{ { \gamma  }_{ T } } \sqrt { 2{ \alpha  }_{ 1 }{ R }^{ 2 } } +\frac { { \gamma  }_{ T-1 } }{ { \gamma  }_{ T } } ||{ { { \bm{ g } } } }_{ T-1 }({ x }_{ T }){ || }+\frac { G }{ { \gamma  }_{ T } } \sqrt { \frac { T-1 }{ \delta  }  } \\
        &+\frac { 2 }{ { \gamma  }_{ T } } \sqrt { 2R\sum _{ t=1 }^{ T-1 }{ { \alpha  }_{ t }||{ x }_{ t+1 }^{ * }-{ x }_{ t }^{ * }{ || } }  } +\frac { 2 }{ { \gamma  }_{ T } } \sqrt { F\sum _{ t=1 }^{ T-1 }{ { \gamma  }_{ t-1 }^{ 2 }||{ { { \bm{ g } } } }_{ t-1 }({ x }_{ t })-{ { { \bm{ g } } } }_{ t }({ x }_{ t }){ || } }  } +{ V }_{ \bm{ g } }\\ 
        &\overset { (a) }{ \le  } { T }^{ 1/4 }\sqrt { \frac { 1 }{ 2\sqrt { 2R } { \beta  }^{ 2 } }  } [2\sqrt { F(T-1) } +R\sqrt { 2\sqrt { \frac { T }{ R }  }  } +G\sqrt { 2\frac { (T-1)\sqrt { R+{ V }_{ x } }  }{ \sqrt { T }  }  } ]
\end{align*}
\begin{equation}
    \begin{split}
        &+{ (\frac { T+1 }{ T } ) }^{ 1/4 }F+{ T }^{ 1/4 }\sqrt { \frac { 2 }{ \sqrt { 2R } { \beta  }^{ 2 } }  } [\sqrt { 4R\sqrt { T{ V }_{ x } }  } +\frac { F }{ { \beta  } } \sqrt { \frac { \sqrt { 2T }  }{ \sqrt { R }  }  } ]+{ V }_{ \bm{ g } }\\ 
        &\overset { (b) }{ \le  } { T }^{ 1/4 }\sqrt { \frac { 1 }{ 2\sqrt { 2R } { \beta  }^{ 2 } }  } [2\sqrt { F(T-1) } +R\sqrt { 2\sqrt { \frac { T }{ R }  }  } +G\sqrt { 2R\sqrt { T(1+T } ) } ]\\ 
        &+{ 2 }^{ 1/4 }F+{ T }^{ 1/4 }\sqrt { \frac { 2 }{ \sqrt { 2R } { \beta  }^{ 2 } }  } [\sqrt { 4R\sqrt { T{ V }_{ x } }  } +\frac { F }{ { \beta  } } \sqrt { \frac { \sqrt { 2T }  }{ \sqrt { R }  }  } ]+{ V }_{ \bm{ g } }
        =O(max\{ { T }^{ 3/4 },{ V }_{ \bm{g} }\} )
    \end{split}
\end{equation}
Where (a) follows from (\ref{3_proof_of_theorem_1}) and (\ref{4_proof_theorem_1}); (b) holds by the fact that  ${V}_{x}\le RT$ and $\frac{T+1}{T}\le 2$. It completes the proof.
\begin{remark}
When feasible set $\chi$ is time-varying, i.e., $\chi(t)$, our algorithm VQB is valid and we can verify that its corresponding theoretical results also hold.
\end{remark}
\subsection{Slater condition}
In the previous section, we have shown that as long as the optimization problem is feasible, our algorithm could simultaneously achieve sublinear dynamic regret and constraint violation 
%without introducing any more assumptions than the common assumptions in the literature of constrained OCO
with only limited common assumptions in the literature of constrained OCO. 
Meanwhile, \cite{chen2017an} pointed out that in many practical constrained OCO problems, the variation of consecutive constraints is smooth across time.
 Thus, we examine whether the smoothness of the dynamic environment's temporal variations can lead to better bounds of constraint violations for VQB. %virtual-queues-based-algorithms.
 %Theorem $1$ shows that 
%it is known that the best bounds of constraint violations for OCO with long term and time-invariant constraints are $O(1)$ \cite{yu2020a,qiu2020beyond}.
%, achieved by virtual-queues-based algorithm. 
%Thus 
 Within this subsection, we consider a slightly stronger Slater condition that has been considered in \cite{chen2017an}. We will show that our variant of VQB, illustrated in Algorithm 2, could guarantee the $O(1)$ constraint violations under this assumption. The difference between VQB and Algorithm 2 is the way of incorporating constraints into the virtual queues updates and decision iterations, i.e., Algorithm 2 uses ${\bm{g}}_{t}({x}_{t})$ instead of ${\bm{g}}_{t-1}({x}_{t})$ to update the dual iterate ${\bm{\lambda}}_{t}$ and primal iterate ${x}_{t}$ compared with VQB.
 Technically, the update step in Algorithm 2 can yield much lower constraint violations when the variation of consecutive constraints is smooth across time, as shown in Theorem 2. First, we give the definition of the Slater condition.
 \begin{algorithm}[tb]
\caption{}
\label{algorithm2}
%\textbf{Output}: decision sequence $\{ { x }_{ t }\}$
\begin{algorithmic}[1] %[1] enables line numbers
\STATE \textbf{Initialize}: ${\alpha},\gamma>0$, ${\bm{g}}_{0}=\bm{\lambda}(0)=0$, and ${x}_{1} \in \chi$.
\FOR{round $t=1...T-1$}
\STATE Update the dual iterate $\bm{\lambda}(t)$:
\STATE ${ \bm{\lambda}  }(t)=max\{ { \bm{\lambda}  }(t-1)+{ \gamma \bm{g} }_{ t }({ x }_{ t }),-{ \gamma \bm{g} }_{ t }({ x }_{ t })\} $
\STATE Update the primal iterate that satisfies:
\STATE ${ x }_{ t+1 }={ argmin }_{ x\in \chi  } { \nabla { f }_{ t }({ x }_{ t }) }^{ T }(x-{ x }_{ t })+[\bm{\lambda} (t)+\gamma { \bm{g} }_{ t }({ x }_{ t }){ ] }^{ T }(\gamma { \bm{g} }_{ t }(x))+{ \alpha  }||x-{ x }_{ t }{ || }^{ 2 }$
\STATE Choose the action ${x}_{t+1}$
\ENDFOR
\end{algorithmic}
\end{algorithm}
\begin{assumption}
\label{assumption 2}
(Slater condition). There exists %a constant
$\epsilon>0$ and $\hat{x}\in \chi$ such that ${ \bm{g} }_{ t }(\hat { x } )\le -\epsilon I,\forall t$.
\end{assumption}
Assumption \ref{assumption 2} is known as the interior point condition or Slater condition, which is also used widely in the literature of OCO with time-varying constraints \cite{chen2017an,chen2019bandit,sharma2020distributed,yi2020distributed}. Based on Assumption 2, we next introduce a slightly stronger Slater condition assumption, which is valid in many practical scenarios \cite{chen2017an,xu2019efficient,amiri2019towards}.
\begin{assumption}
\label{assumption 3}
The slater constant $\epsilon$ is larger than the maximum variation of consecutive  constraints, i.e., $\epsilon >{ \bar { V }  }_{ \bm{g} }=\underset { t }{ \max }  \underset { x\in \chi  }{ \max } ||{ \bm{g} }_{ t+1 }(x)-{ \bm{g} }_{ t }(x)||$. %which is equivalent to that $\epsilon =\underset { k,t }{ \min } \underset { x\in \chi  }{\max } [-{ g }_{ t,k }(x){ ] }^{ + }>\underset { t }{ \max }  \underset { x\in \chi  }{ \max } ||{ \bm{g} }_{ t+1 }(x)-{ \bm{g} }_{ t }(x)||$.
\end{assumption}
Note that this assumption was adopted in \cite{chen2017an}, which is valid when the region defined by $\{ x|x\in \chi ,{ g }_{ t }(x)\le 0\} $ is large enough, or the variation of consecutive constraints is smooth enough across time.

With a similar intuition of VQB, if we define $\bm{Q}(t)={ { \bm{\lambda}  } }(t)+{ \gamma  }_{ t}{ \bm{ g } }_{ t }({ x }_{ t })=\max  \{ { \bm{ \lambda  } }(t-1)+2{ \gamma  }_{ t }{ \bm{ g } }_{t}({ x }_{ t }),0\} $ as the vector of virtual queue backlogs and let parameters ${\alpha}_{t},{\gamma}_{t}$ be time-invariant, Algorithm 2 also chooses ${x}_{t+1}$ to minimize an upper bound of the following expression:
\begin{small}
\begin{align*}
    &\underbrace { \Delta (t) } +\underbrace { { \nabla { f }_{ t }({ x }_{ t }) }^{ T }(x-{ x }_{ t })+{ \alpha  }_{ t }||x-{ x }_{ t }{ || }^{ 2 } }. \\ 
    & drift \quad \quad \quad \quad \quad \quad penalty
\end{align*}
\end{small}
The reason to replace ${\bm{g}}_{t-1}({x}_{t})$ with ${\bm{g}}_{t}({x}_{t})$ in Algorithm 2 is motivated by the observation that ${\bm{g}}_{t}({x}_{t})$ could be directly accumulated
into queue ${\bm{\lambda}}(t)$ (recall that ${ \bm{\lambda}  }(t)=\max\{ { \bm{\lambda}  }(t-1)+{ \gamma \bm{g} }_{ t }({ x }_{ t }),-{ \gamma \bm{g} }_{ t }({ x }_{ t })\} $) and we intend to have small queue backlogs when the variation of consecutive constraint functions is smooth across time. This is important
for a much tighter analysis of the constraint violations under the strongly Slater condition. If we use ${\bm{g}}_{t}({x}_{t})$ instead of ${\bm{g}}_{t-1}({x}_{t})$ in 
VQB, we could get 
${ \gamma g }_{ t,k }({ x }_{ t })\le { \lambda  }_{ k }(t)-{ \lambda  }_{ k }(t-1)$. Then, we can characterize the constraint violations only by the bounds of virtual queues $\{ \bm{\lambda} (t)\} $ without the term ${V}_{g}$ (comparing with the Lemma 5), i.e., 
$\sum _{ t=1 }^{ T }{ { g }_{ t,k }({ x }_{ t }) } \le \frac { { \lambda  }_{ k }(T) }{ \gamma  } \le \frac{||\bm{\lambda}(T)||}{\gamma},\forall k.$
In such case, the length of virtual queues $\{||{\bm{\lambda}}(t)||{\}}_{t=1}^{T}$
 is upper bounded by a constant under the strongly Slater condition (Lemma 11). Then we could obtain an $O(1)$ bound of constraint violations, shown in the following theorem.
%
%
%
% Put Algorithm 3 here
\begin{theorem}
Under Assumption \ref{assumption 1},\ref{assumption 2} and \ref{assumption 3}, setting ${\alpha} ={T}^{a}$ and ${\gamma}^{2}=\frac{1}{2{\beta}^{2}}{T}^{a}$ in Algorithm \ref{algorithm2}, the dynamic regret and constraint violations are upper bounded by 
\begin{equation}
    \begin{split}
        &regret\le O(\max\{{T}^{a}{V}_{x},{T}^{a}{V}_{g},{T}^{1-a}\} )\\
        &{ Vio }_{ k }\le O(1),\; \forall k=1,2,...,K
    \end{split}
\end{equation}
In particular, the performance upper bounds become $O(\max$ $\{\sqrt { T }{V}_{x},\sqrt{T}{V}_{g}\} )$ and $O(1)$ if we set ${\alpha} =2{\beta}^{2}{\gamma}^{2} =\sqrt { T  }$.
\end{theorem}
Note that our performance bounds established by Algorithm 2
are strictly better than \cite{chen2017an} under the same assumptions. Besides, the constraint violations for Algorithm \ref{algorithm2} can decrease into $O(1)$ when the variations of consecutive constraints are smooth enough across time.
%algorithm here
%algorithm3

~\\
\textbf{Proof sketch of Theorem 2.} Here we give a proof sketch of Theorem 2. All the proof details of listed lemmas could be found in the Appendix. Note that both $\{{\alpha}_{t}\}$ and $\{{\gamma}_{t}\}$ are constant sequences in Algorithm 2, thus here we omit the subscript $t$. Similar as Lemma \ref{lemma 0}, in Algorithm 2, we have the following lemma for the properties of virtual queues and Lyapunov drift term.
\begin{lemma}
\label{lemma 6}
In algorithm 2, at each round $t$, we have
\begin{enumerate}
    \item ${\bm{\lambda}}(t)\ge \bm{0}$ 
    \item ${ \bm{\lambda}  }(t)+{ \gamma \bm{g} }_{ t }({ x }_{ t })\ge \bm{0}$ 
    \item ${ ||\bm{\lambda}  }(t){ || }\ge { \gamma ||\bm{g} }_{ t-1 }({ x }_{ t }){ || }$
    \item ${ \gamma \bm{g} }_{ t }({ x }_{ t })\le { \bm{\lambda}  }(t)-{ \bm{\lambda}  }(t-1)$, furthermore, ${ ||\bm{\lambda}  }(t){ || }-{ ||\bm{\lambda}  }(t-1){ || }{ \le \gamma ||\bm{g} }_{ t }({ x }_{ t }){ || }$
    \item $\Delta (t)\le \gamma [\bm{\lambda} (t){ ] }^{ T }{ \bm{g} }_{ t+1 }({ x }_{ t+1 })+{ \gamma  }^{ 2 }||{ \bm{g} }_{ t+1 }({ x }_{ t+1 }){ || }_{ 2 }^{ 2 }$ \label{5_ in lemma 12}
\end{enumerate}
\end{lemma}
The proof of this Lemma is similar as the proof of Lemma \ref{lemma 0} and hence we omit the details.
\begin{lemma}
\label{lemma 7}
Under Assumption 1,2 and 3, setting $\alpha,\gamma,\delta$ such that $\alpha\ge \frac{1}{2}({\beta}^{2}{\gamma}^{2}+\delta)$, then  Algorithm 2 ensures that
\begin{equation}
    \begin{split}
     &{ f }_{ t }({ x }_{ t })+\Delta (t)\\ 
     &\le { f }_{ t }({ x }_{ t }^{ * })+{ \alpha  }||{ x }_{ t }^{ * }-{ x }_{ t }{ || }^{ 2 }-{ \alpha  }||{ x }_{ t+1 }{ -{ x }_{ t+1 }^{ * }|| }^{ 2 }+4R{ \alpha  }{ ||{ x }_{ t+1 }^{ * }-{ x }_{ t }^{ * }|| }^{ 2 }+{ \gamma  }^{ 2 } F||{ { { \bm{ g } } } }_{ t+1 }({ x }_{ t+1 })-{ { { { \bm{ g } } } }_{ t }({ x }_{ t+1 })|| } \\ 
     &+\frac { 1 }{ 2\delta  } { G }^{ 2 }+\frac { 1 }{ 2 } { \gamma  }^{ 2 }||{ { { \bm{ g } } } }_{ t+1 }({ x }_{ t+1 }){ || }^{ 2 }-\frac { 1 }{ 2 } { \gamma  }^{ 2 }||{ { { \bm{ g } } } }_{ t }({ x }_{ t }){ || }^{ 2 }+\gamma ||{ { \bm{\lambda}  } }(t){ || }||{ { { \bm{ g } } } }_{ t+1 }({ x }_{ t+1 })-{ { { \bm{ g } } } }_{ t }({ x }_{ t+1 }){ || }
    \end{split}
\end{equation}
\end{lemma}
 Take a similar derivation process as the proof of Theorem 1, we also characterize the regret and constraint violations through the bound of drift-plus-penalty expression stated above. Therefore, we bound the dynamic regret and constraint violations in the following lemmas, respectively.
\begin{lemma}
\label{lemma 8}
Under the Assumption 1, 2 and 3, setting $\alpha,\gamma,\delta$ such that $\alpha\ge \frac{1}{2}{\beta}^{2}{\gamma}^{2}+\frac{1}{2}\delta$, then Algorithm 2 ensures that
\begin{equation}
    \begin{split}
        &\sum _{ t=1 }^{ T }{ { f }_{ t }({ x }_{ t }) } \le \sum _{ t=1 }^{ T }{ { f }_{ t }({ x }_{ t }^{ * }) } +{ \alpha  }||{ x }_{ 1 }^{ * }-{ x }_{ 1 }{ || }^{ 2 }+4R\alpha {V}_{x} 
        +\frac { T{ G }^{ 2 } }{ 2\delta  } +{ \gamma  }^{ 2 } F{ V }_{ g }+L(1)+\gamma \max _{ t }{ ||{ { \bm{\lambda}  } }(t){ || } } { V }_{ \bm{ g } }
    \end{split}
\end{equation}
\end{lemma}
\begin{lemma}
\label{lemma 9}
In Algorithm 2, we have
\begin{equation}
    \sum _{ t=1 }^{ T }{ { g }_{ t,k }({ x }_{ t }) } \le \frac { { \lambda  }_{ k }(T) }{ \gamma  }\le \frac{||\bm{\lambda}(T)||}{\gamma},\; \forall k\in \{ 1,2,...,K\}. 
\end{equation}
\end{lemma}
Note that the above two lemmas show that the final bounds of dynamic regret and constraint violations can be obtained by bounding the ${\bm{\lambda}}(t)$.
Since we are allowed to introduce the assumption of strongly Slater condition, we will show that the length of virtual queues is upper bounded by a constant in this case. Hence we adopt different techniques for the virtual queues analysis compared with the proof of Lemma \ref{lemma 3}, and bound them by accomplishing the following lemma.
%Hence next we do this by accomplishing the following lemma.
\begin{lemma}
\label{lemma 10}
In Algorithm 2, we have
\begin{equation}
    ||{ \bm{\lambda}  }(t){ || }\le \gamma  F+\frac { GR+{ \gamma  }^{ 2 }\epsilon  F+2{ \gamma  }^{ 2 }{ F }^{ 2 }+{ \alpha  }{ R }^{ 2 }  }{ \gamma (\epsilon -{ \bar { V }  }_{ g }) },\;\forall t.
\end{equation}
\end{lemma}
Finally, based on the above lemmas, we prove the Theorem 2 as follows. We setting $\delta=\frac { 1 }{ 2 }{T}^{a}$, and it is easy to verity that ${\alpha}_{t}\ge \frac{1}{2}{ \beta  }^{ 2 }{ \gamma  }^{ 2 }+\frac { \delta  }{ 2 }$. Combining the Lemma \ref{lemma 8} with Lemma \ref{lemma 10}, we have
\begin{equation}
    \begin{split}
        &\sum _{ t=1 }^{ T }{ { f }_{ t }({ x }_{ t }) } -\sum _{ t=1 }^{ T }{ { f }_{ t }({ x }_{ t }^{ * }) } \\ 
        &\le { \alpha  }||{ x }_{ 1 }^{ * }-{ x }_{ 1 }{ || }^{ 2 }+4R\alpha { V }_{ x }+\frac { T{ G }^{ 2 } }{ 2\delta  } +{ \gamma  }^{ 2 }F{ V }_{ \bm{g} }+\frac { 1 }{ 2 } ||{ \bm{\lambda}  }(1){ || }^{ 2 }+\gamma \max _{ t }{ ||{ { \bm{\lambda}  } }(t){ || } } { V }_{ \bm{ g } }\\ 
        &\le \sum _{ t=1 }^{ T }{ { f }_{ t }({ x }_{ t }^{ * }) } +{ \alpha  }{ R }^{ 2 }+4R\alpha { V }_{ x }+\frac { T{ G }^{ 2 } }{ 2\delta  } +{ \gamma  }^{ 2 }F{ V }_{ \bm{g} }+\frac { 1 }{ 2 } ||{ \bm{\lambda}  }(1){ || }^{ 2 }+({ \gamma  }^{ 2 }F+\frac { GR+{ \gamma  }^{ 2 }\epsilon F+2{ \gamma  }^{ 2 }{ F }^{ 2 }+{ \alpha  }{ R }^{ 2 } }{ (\epsilon -{ \bar { V }  }_{ \bm{g} }) } ){ V }_{ \bm{g} }\\ 
        &=O(max\{ {T}^{a} { V }_{ x },{ \gamma  }^{ 2 }{ V }_{ \bm{g} },{ T }^{ 1-a }\} )=O(max\{ {T}^{a} { V }_{ x },{ V }_{ \bm{g} }{ T }^{ a },{ T }^{ 1-a }\} )
    \end{split}
\end{equation}
Furthermore, combining the Lemma \ref{lemma 9} with Lemma \ref{lemma 10}, we can obtain
\begin{equation}
    \begin{split}
        &\sum _{ t=1 }^{ T }{ { g }_{ t,k }({ x }_{ t }) } \le \frac { ||{ \bm{\lambda}  }(T)|| }{ \gamma  } \le F+\frac { GR+{ \gamma  }^{ 2 }\epsilon F+2{ \gamma  }^{ 2 }{ F }^{ 2 }+{ \alpha  }{ R }^{ 2 } }{ { \gamma  }^{ 2 }(\epsilon -{ \bar { V }  }_{ \bm{g} }) } \\ 
        &=O(\frac { \alpha  }{ { \gamma  }^{ 2 } } )=O({ T }^{ 1-a }),\; \forall k=1,2,...,K.
    \end{split}
\end{equation}
 In particular, when setting $\alpha=\frac{1}{2}\sqrt{T}$ and ${\gamma}^{2} =\frac{1}{2{\beta}^{2}}\sqrt { T  } $ , the performance upper bounds become $O(max\{\sqrt { T }{V}_{x},\sqrt{T}{V}_{\bm{g}}\} )$ and $O(1)$. This completes the proof.
\subsection{Unknown time horizon \emph{T}}
%Put Algorithm 2 here
\begin{algorithm}[tb]
\caption{The Doubling Trick for Online Algorithm A}
\label{algorithm_doubling_trick}
\begin{algorithmic}[1] %[1] enables line numbers
\STATE Let $i=1$.
\WHILE{not reach the end of the time horizon}
\STATE Reset \emph{A} with parameters chosen for $T={2}^{i}$.
\STATE Run \emph{A} for ${2}^{i}$ rounds.
\STATE Let $i=i+1$.
\ENDWHILE
\end{algorithmic}
\end{algorithm}
In this subsection, we extend our algorithm and analysis 
%in the previous parts by considering unknown time horizon $T$ case.
%to the unknown time horizon case.
the case when time horizon is unknown.
Recall that the parameters of VQB, i.e., ${\gamma}_{t}$ and ${\alpha}_{t}$, and previous methods for OCO with long-term and time-varying constraints all depend on the time horizon \emph{T}, while the total rounds $T$ is not known prior to the learner in many practical scenarios. In such cases, we use the doubling trick strategy to tune the parameters for our algorithm. To the best of our knowledge, we are the first to consider the unknown time horizon case for OCO with long-term and time-varying constraints.
For any online Algorithm \emph{A} whose parameters depend on the time horizon $T$, the doubling trick is described in Algorithm 3.
%  preserve the order of our algorithms' performance upper bounds
%
\begin{theorem}
Under Assumption \ref{assumption 1}, for any unknown time horizon T, run VQB until reaching the end of the time horizon. Let ${t}_{i}={2}^{i}$ be the index of the first round of i-th epoch.
\begin{itemize}
    \item (Case 1) Setting $ {\alpha}_{t}=\sqrt { \frac { {2}^{i} }{ R+\sum _{ {t}_{i}\le i\le t }{ ||{ x }_{ i }^{*}-{ x }_{ i-1 }^{*}{ || } }  }  }$ and ${ \gamma  }_{t}^{2}=\frac { 1 }{ 2{ \beta  }^{ 2 } } \frac { 1 }{ \sqrt { 2R }  }$   for $t\in [{ t }_{ i },{ t }_{ i+1 }-1]$, then VQB with doubling trick ensures 
\begin{equation}
    \begin{split}
      &regret \le O(\max\{ \sqrt { T{ V }_{ x } } ,{ V }_{ g }\} ),\\
      &{ Vio }_{ k }\le O(\max\{ \sqrt { T } ,{ V }_{ g }\} ),\forall k=1,2,...,K.
    \end{split}
\end{equation}
   \item (Case 2) Setting ${\alpha}_{t} =\sqrt { \frac { {2}^{i} }{ R+\sum _{ {t}_{i}\le i\le t }{ ||{ x }_{ i }^{*}-{ x }_{ i-1 }^{*}{ || }}  }  } $ and ${\gamma}_{t}^{2} =\frac { 1 }{ 2{ \beta  }^{ 2 } } \frac { 1 }{ \sqrt { 2R }  }\frac{1}{\sqrt{t-{t}_{i}+2}}$ for $t\in [{ t }_{ i },{ t }_{ i+1 }-1]$, then VQB with doubling trick ensures 
\begin{equation}
    \begin{split}
      &regret \le O( \sqrt { T{ V }_{ x } } ),\\
      &{ Vio }_{ k }\le O(\max\{ {T}^{\frac{3}{4}},{ V }_{ g }\} ),\forall k=1,2,...,K.
    \end{split}
\end{equation}
\end{itemize}
\end{theorem}
%Theorem 2 shows that using the doubling trick our algorithm can still preserve the order of
Theorem 3 shows that our algorithm with the doubling trick can still preserve the order of
dynamic regret and constraint violations bounds even though the time horizon $T$ is unknown. 
%While for parameter-dependent methods like \cite{cao2018online}, it is hard to use the doubling trick for them since the order of the environment's temporal variations is unknown in advance, which shows the superiority of parameter-free algorithms.
Note that our algorithm adapts to the doubling trick because of the property of parameter-free, while parameter-dependent methods (e.g., \cite{cao2018online}) cannot do this.

~\\
\textbf{Proof sketch of Theorem 3.}
Here we give a proof of Theorem 3. For the \textbf{Case 1} in Theorem 3, since the \emph{i}-th epoch consists of at most ${2}^{i}$ rounds, the time horizon is divided into $N=\lceil {log}_{2}T\rceil$ epochs. Let ${ \Delta  }_{ x }^{ i }=\sum _{ t={ t }_{ i } }^{ t={ t }_{ i }-1 }{ ||{ x }_{ t }^{ * }-{ x }_{ t-1 }^{ * }|| } $ and ${ \Delta  }_{ g }^{ i }=\sum _{ t={ t }_{ i } }^{ t={ t }_{ i }-1 }{ \underset { x\in \chi  }{ sup } ||{ \bm{g} }_{ t }(x)-{ \bm{g} }_{ t-1 }(x)|| } $.
By Theorem 1, in the \emph{i}-th epoch there exists a constant \emph{C} such that the dynamic regret and constraint violations are at most $max\{ C\sqrt { T{ \Delta  }_{ x }^{ i } } ,C{ \Delta  }_{ g }^{ i }\}$ and $max\{ C\sqrt { {2}^{i} } ,C{ \Delta  }_{ g }^{ i }\} $ respectively. The final bound could be obtained by summing the individual bounds over all the epochs. Therefore, we could upper bound the total dynamic regret and constraint violations as follows
\begin{equation}
    \begin{split}
        \label{regret bound 1 with unknown time horizon}
        Regret &\le max\{ \sum _{ i=0 }^{ N }{ C\sqrt { { 2 }^{ i }{ \Delta  }_{ x }^{ i } }  } ,\sum _{ i=0 }^{ N }{ C{ \Delta  }_{ g }^{ i } } \} 
        \overset { (a) }{ \le  } max\{ C\sqrt { \sum _{ i=0 }^{ N }{ { \Delta  }_{ x }^{ i } }  } \sqrt { \sum _{ i=0 }^{ N }{ { 2 }^{ i } }  } ,C\sum _{ i=0 }^{ N }{ { \Delta  }_{ g }^{ i } } \} \\ 
        &\le max\{ C\sqrt { { V }_{ x } } \sqrt { { 2 }^{ N+1 }-1 } ,C\sum _{ i=0 }^{ N }{ { \Delta  }_{ g }^{ i } } \} 
        =max\{ C\sqrt { { V }_{ x } } \sqrt { { 2 }^{ \lceil { log }_{ 2 }T\rceil +1 }-1 } ,C{ V }_{ \bm{g} }\} \\ 
        &\le max\{ C\sqrt { { V }_{ x } } \sqrt { { 2 }^{ { log }_{ 2 }T+2 } } ,C{ V }_{ \bm{g} }\} =max\{ 2C\sqrt { T{ V }_{ x } } ,C{ V }_{ \bm{g} }\} 
    \end{split}
\end{equation}
%Conducting similar analysis as before, 
Where (a) is due to the Cauchy-schwarz inequality. And the total constraint violations are at most
\begin{equation}
    \begin{split}
        { Vio }_{ k }&\le max\{ \sum _{ i=0 }^{ N }{ C\sqrt { { 2 }^{ i } }  } ,\sum _{ i=0 }^{ N }{ C{ \Delta  }_{ g }^{ i } } \} 
        =max\{ \frac { C }{ \sqrt { 2 } -1 } [{ \sqrt { 2 }  }^{ N+1 }-1],C{ V }_{ \bm{g} }\} \\
        &\le max\{ \frac { 2C }{ \sqrt { 2 } -1 } { \sqrt { 2 }  }^{ { log }_{ 2 }T },C{ V }_{ \bm{g} }\} 
        =max\{ \frac { 2C }{ \sqrt { 2 } -1 } \sqrt { T } ,D{ V }_{ \bm{g} }\},\forall k=1,2,...,K.
    \end{split}
\end{equation}
For the \textbf{Case 2}, we conducting similar analysis as Case 1. By Theorem 1, in the \emph{i}-th epoch there exists a constant \emph{D} such that the dynamic regret and constraint violations are at most $D\sqrt { T{ \Delta  }_{ x }^{ i } } $ and $max\{ D{ 2 }^{ \frac { 3 }{ 4 } i },D{ \Delta  }_{ g }^{ i }\}  $ respectively. According to (\ref{regret bound 1 with unknown time horizon}), the total regret is still at most the order of $\sqrt{T{V}_{x}}$ without changing. For the total constraint violations, we also have
\begin{equation}
    \begin{split}
        { Vio }_{ k }&\le max\{ D\sum _{ i=0 }^{ N }{ { 2 }^{ \frac { 3 }{ 4 } i } } ,\sum _{ i=0 }^{ N }{ D{ \Delta  }_{ g }^{ i } } \}
        \le max\{ \frac { D }{ { 2 }^{ 3/4 }-1 } ({ 2 }^{ \frac { 3 }{ 4 } \lceil { log }_{ 2 }T\rceil  }-1),D{ V }_{ \bm{g} }\} \\ 
        &\le max\{ \frac { 2D }{ { 2 }^{ 3/4 }-1 } { 2 }^{ \frac { 3 }{ 4 } { log }_{ 2 }T },D{ V }_{ \bm{g} }\} 
        =max\{ \frac { 2D }{ { 2 }^{ 3/4 }-1 } { T }^{ 3/4 },D{ V }_{ \bm{g} }\},\forall k=1,2,...,K.
    \end{split}
\end{equation}
This completes the proof.
\section{Numerical experiments}
In this section, we conduct numerical experiments to validate the theoretical performance of our algorithm. Specifically, we consider the online ridge regression (ORR) problem \cite{arce2012online} as the numerical example. We compare the time-averaged regrets and constraint violations of our algorithm 
%(under two different parameter settings) 
with previous work in two different dynamic environments.
%We observe that the time average regrets and constraint violations of Algorithm $1$ and previous works all converge to zero as time goes to infinity, as long as the temporal variations of the environment is not too large, and our algorithm's convergence rate is superior to those previous works. 
The problem formulation of ORR at round $t$ is as follows.
\begin{small}
\begin{equation}
    \begin{split}
        &\underset { { x }_{ t }\in \chi  }{ \text{Minimize} } \sum _{ i=1 }^{ n }{ ({ x }_{ t }^{ T }{ p }_{ i,t }+{ b }-{ q }_{ i,t }{ ) }^{ 2 } } \\ 
        &\quad \quad \quad \; s.t.\quad ||{ x }_{ t }||\le { a }_{ t }
    \end{split}
\end{equation}
\end{small}
Here $\{ { p }_{ i,t },{ q }_{ i,t }{ \}  }_{ i=1 }^{ n }$ are the training data at round $t$ and ${a}_{t}$ characterizes the $t$-th round constraint on the ${l}_{2}$ norm of the decision variable, i.e., weight vector. We define $\{ x|||x{ || }_{ \infty  }\le C,x\in {R}^{k}\} $ as the feasible set.
The above ORR formulation %arises from
could be applied in accurate and reliable forecasting of traffic in intelligent transportation systems \cite{haworth2014local}. The training data $\{ { p }_{ i,t },{ q }_{ i,t }{ \}  }_{ i=1 }^{ n }$ and %real-time
constraint ${a}_{t}$ may not be known prior to the
agent at round $t$ due to the delayed arrival training data.

\textbf{Experimental setting.} At round $t$, 
%we first update the true weight vector as follows.
we generate the parameters $\{{p}_{i,t},{q}_{i,t},\forall i\}$, ${a}_{t}$ and the per-slot minimizer ${x}_{t}^{*}$ in the following way. Let ${ x }_{ t }^{*}={ \Pi  }_{ \chi  }({ x }_{ t-1 }^{*}+{ \tau  }_{ t })$, where each entry of ${\tau}_{t}$ is a uniform random variable, sampled from a time-varying set ${B}_{t}$ (We will specialize it later). %ii) ${b}_{t}^{*}=b$, where $b$ is known prior to the agent. 
Then we generate ${a}_{t}$ and $\{ { p }_{ i,t },{ q }_{ i,t }{ \}  }_{ i=1 }^{ n }$ as follows. i) ${p}_{i,t}={p}_{i,t-1}+{u}_{i,t}$, where each entry of ${u}_{i,t}$ is i.i.d, uniformly sampled from  set ${B}_{t}$. ii) ${ q }_{ i,t }={ p }_{ i,t }^{T}{ x }_{ t }^{ * }+b$. iii) ${ a }_{ t }=||{ x }_{ t }^{ *}||$.

Next, we introduce the baslines \cite{chen2017an,chen2018heterogeneous,cao2018online,chen2019bandit} for comparison. The algorithms in \cite{chen2017an,chen2018heterogeneous} are based on MOSP method. Although \cite{chen2019bandit} only considered the bandit setting, their algorithm and theoretical guarantees are also valid in the full-information setting. Meanwhile, note that the theoretical guarantees in \cite{cao2018online} are valid only when the agent has prior knowledge of ${V}_{x}$ (or the order of it).
For fair comparison, we set the learning rates in their algorithm to be parameter-free, and obtain the $O({V}_{x}{T}^{1/2})$ regret and $O(\max\{{V}_{x}^{1/2}{T}^{1/2},{T}^{3/4}\})$ constraint violations. %The parameters settings of these methods could be found in the supplementary  material.
Finally, we introduce our experimental details. In our experiment, we let $n=k=5,\;C=7$. The parameters of our algorithm and other baselines are presented in Table 2.
\begin{table}[htp]
\centering
\caption{Parameters of our algorithm and baselines}
 \begin{tabular}{p{4.5cm}p{8.5cm}}
\toprule   
  %\multicolumn{1}{c}{N=1} & \quad \quad N=3 &\quad \quad N=5 \\ 
 \multicolumn{1}{l}{\centering Methods} & \multicolumn{1}{l}{\centering Parameters}  \\
  \midrule
  Baseline \cite{cao2018online} & $\delta=8n{C}^{2}+1,\;\eta=\frac{2}{\sqrt{T}}$    \\  
  Baseline  \cite{chen2017an}& $\alpha=\mu={T}^{\frac{1}{3}}$    \\    
  Baseline \cite{chen2018heterogeneous} &$\delta=1$,\;${\lambda}_{1}=4\sqrt{2}{T}^{\frac{1}{8}}$   \\
 Baseline \cite{chen2019bandit}& $\mu={T}^{-\frac{1}{2}},\;\alpha=2 {T}^{-\frac{1}{2}}$ \\
 VQB(Case 1) & ${\alpha}_{t} =\sqrt { \frac { T }{ R+\sum _{ i\le t }{ ||{ x }_{ i }^{*}-{ x }_{ i-1 }^{*}{ || } }  }  } $, ${ \gamma  }_{t}^{2}=\frac { 1 }{ 2{ \beta  }^{ 2 } } \frac { 1 }{ \sqrt { 2R }  }$\\
 VQB(Case 2) & ${\alpha}_{t} =\sqrt { \frac { T }{ R+\sum _{ i\le t }{ ||{ x }_{ i }^{*}-{ x }_{ i-1 }^{*}{ || } }  }  }$, ${\gamma}_{t}^{2} =\frac { 1 }{ 2{ \beta  }^{ 2 } } \frac { 1 }{ \sqrt { 2R }  }\frac{1}{\sqrt{t+1}}$ \\
  \bottomrule  
\end{tabular}
\end{table}

\textbf{Results and analysis.} We first consider the case when
${V}_{x}={V}_{\bm{g}}=O(\log(T))$. To do this, we set ${B}_{t}$ to be $[-\frac{1}{2t},\frac{1}{2t}]$. From figure 1(a) and (b), we can see that our algorithm VQB achieves lowest time-averaged regret $\frac{regret(t)}{t}$ and constraint violation $\frac{Vio(t)}{t}$, which validates our theoretical results. %since ${V}_{\bm{g}}$ is small. 
Moreover, we can also see that the regrets achieved by VQB under two parameter settings are very close, which is consistent with the theoretical results in Theorem 1 that the regret upper bounds between them are identical by noting that ${V}_{\bm{g}}=O(\log(T))$.
\begin{figure*}
\centering
\begin{minipage}[c]{0.5\textwidth}
\centering
\includegraphics[height=5cm,width=6.5cm]{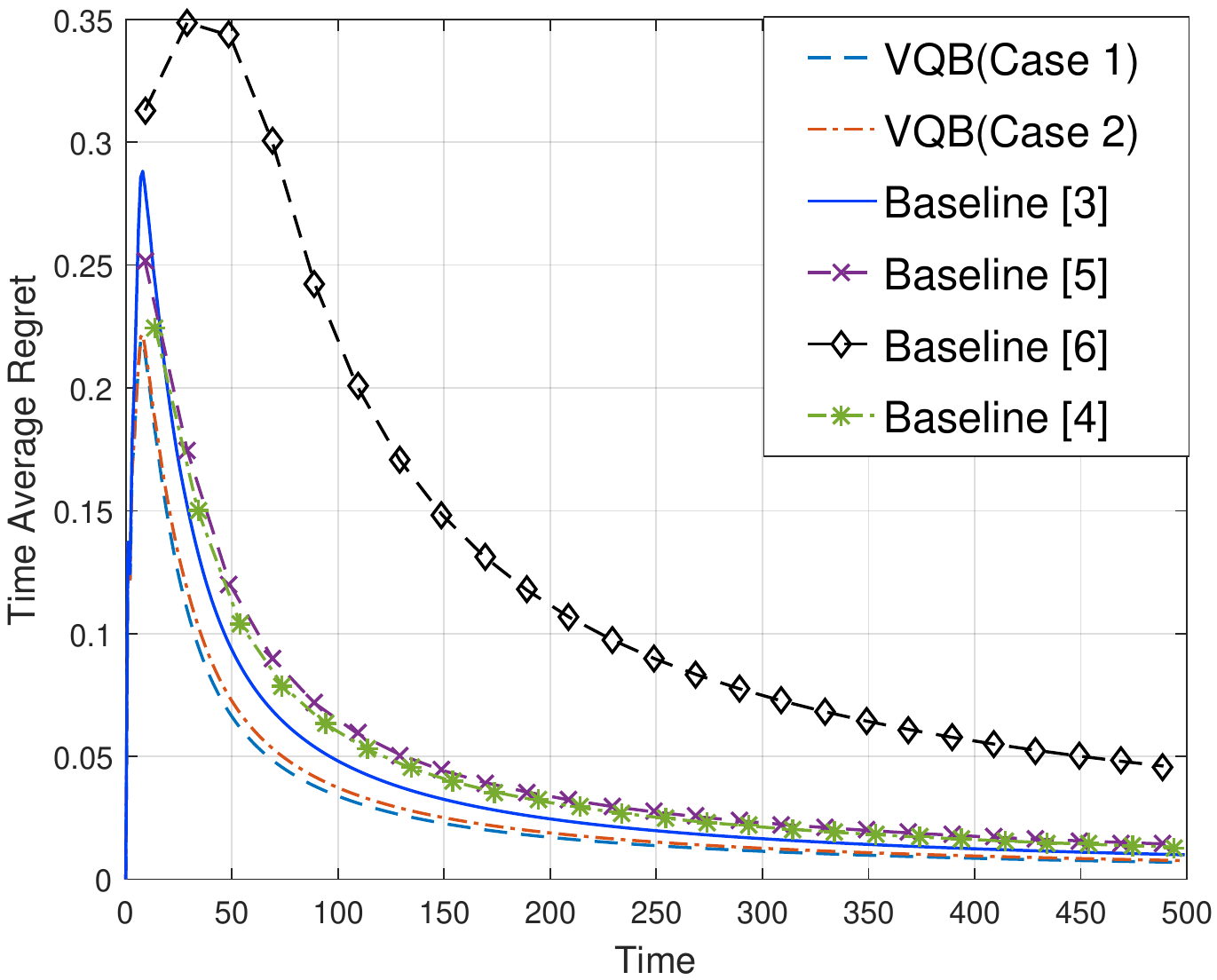}
\end{minipage}%
\begin{minipage}[c]{0.5\textwidth}
\centering
\includegraphics[height=5cm,width=6.5cm]{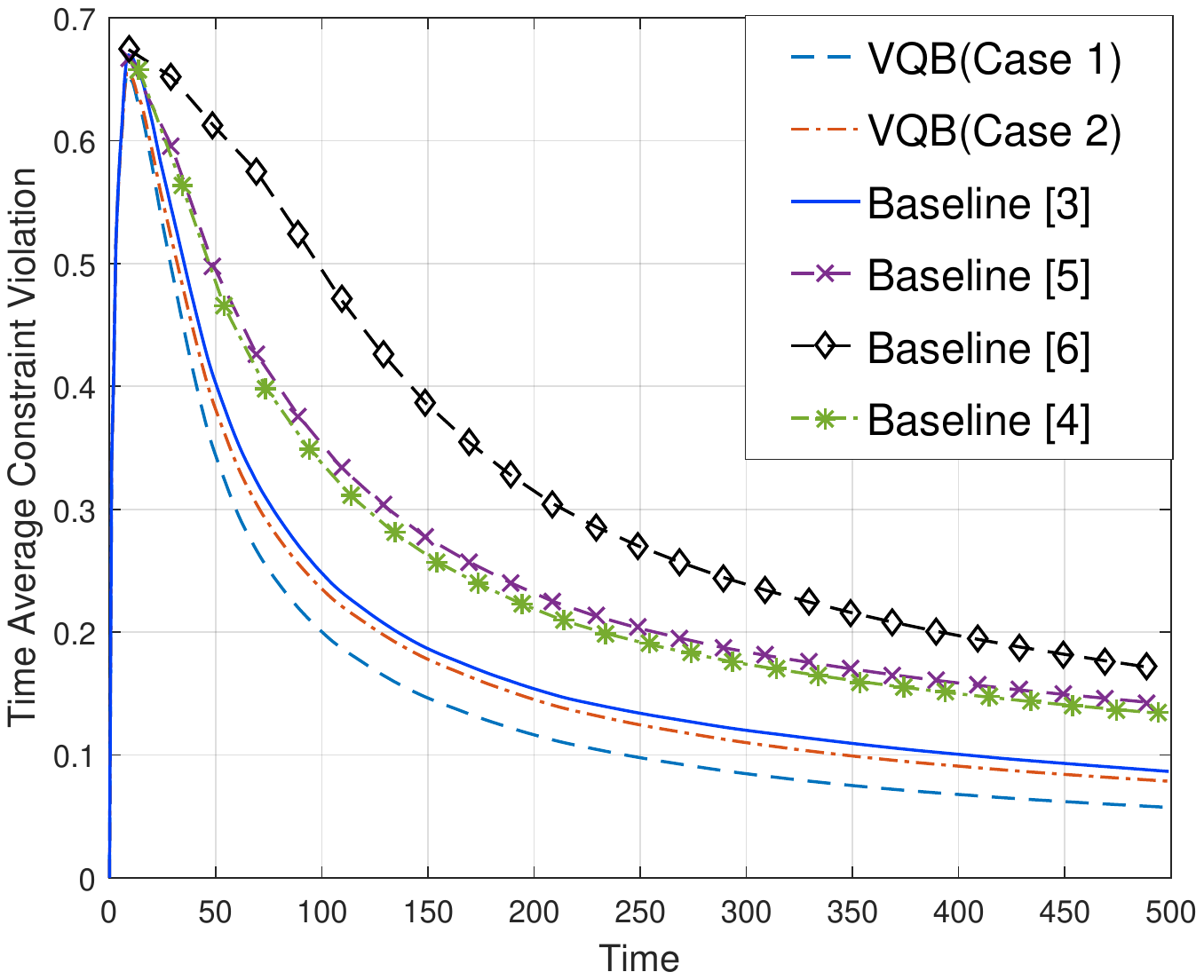}
\end{minipage}
\caption{(a) Evolutions of $Regret(t)/t$; (b) Evolations of $Vio(t)/t$.}
\end{figure*}
\begin{figure*}
\centering
\begin{minipage}[c]{0.5\textwidth}
\centering
\includegraphics[height=5cm,width=6.5cm]{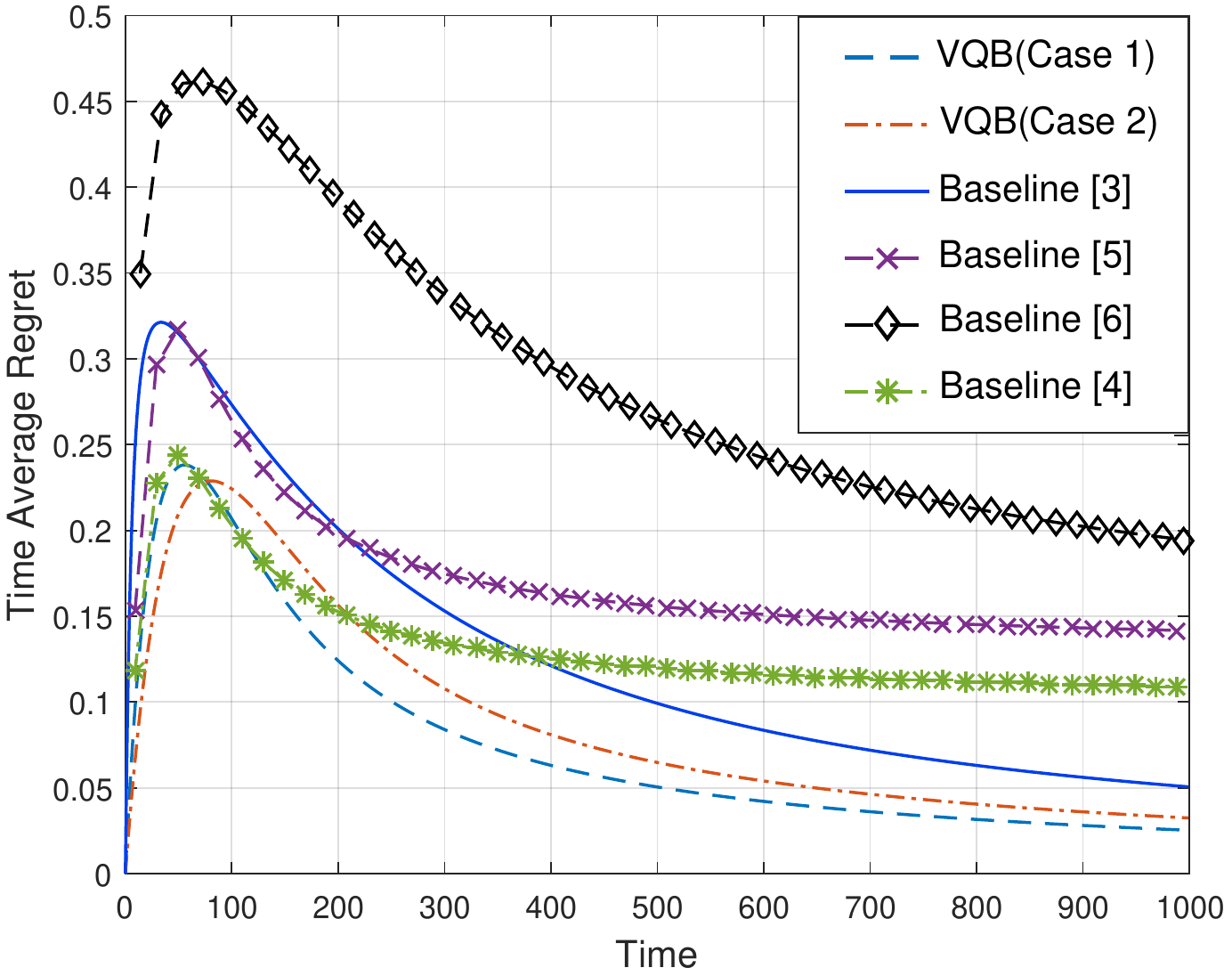}
\end{minipage}%
\begin{minipage}[c]{0.5\textwidth}
\centering
\includegraphics[height=5cm,width=6.5cm]{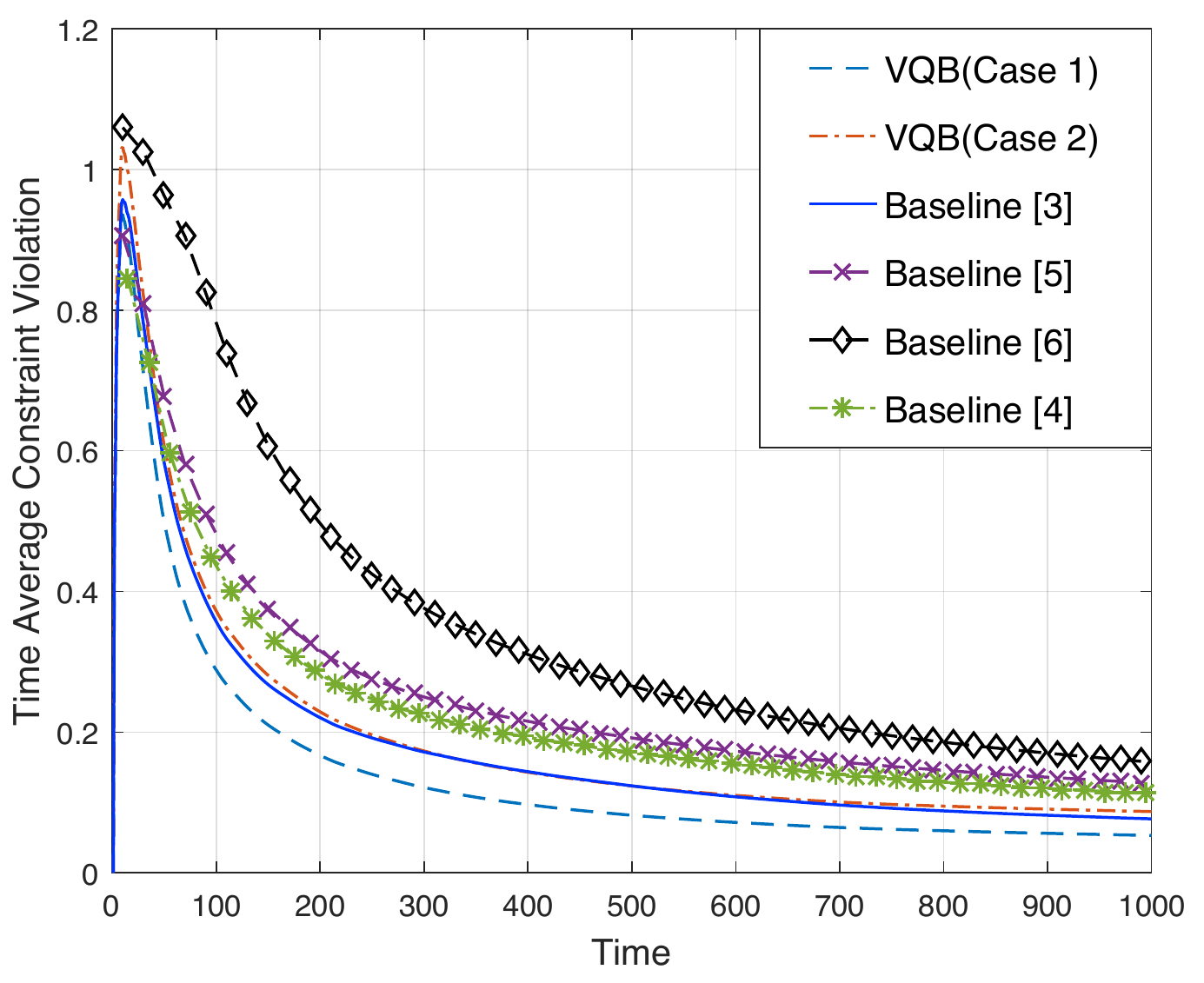}
\end{minipage}
\caption{(a) Evolutions of $Regret(t)/t$; (b) Evolations of $Vio(t)/t$.}
\end{figure*}

We also consider the case when ${V}_{x}={V}_{\bm{g}}=O(\sqrt{T})$. To do this, we set  ${B}_{t}$ to be $[-\frac{1}{2\sqrt{t}},\frac{1}{2\sqrt{t}}]$. In this case, the regret bounds of all baselines are at least the order of $T$. 
From Figure 2, we notice that all methods can guarantee sublinear constraint violation in this case, which matches the theoretical results listed in Table 1. Figure 2 also shows that VQB can achieve simultaneous sublinear regret and constraint violation, while other baselines (\cite{chen2017an,chen2018heterogeneous,chen2019bandit}) cannot, which matches their theoretical results. We observe that baseline \cite{cao2018online} achieves a near sublinear regret in this setting, yet this may not always be the case due to its $O(T)$ regret bound, or the performance bounds established by \cite{cao2018online} may not be tight. Besides, the regrets of the VQB are better, which also coincides our theoretical bounds in this case.
%because of the length of the space. 

\section{Applications}
%In this section, we will outline two examples of applications of the proposed framework.
In this section, we show several applications of our formulation to diverse problems across resource allocation and job scheduling. We emphasize that none of these applications would be possible without a simultaneously achieving sublinear regret and constraint violations algorithm, which has not been attainable with previous approaches.
\begin{figure}[ht]
\centerline{\includegraphics[width=0.6\columnwidth]{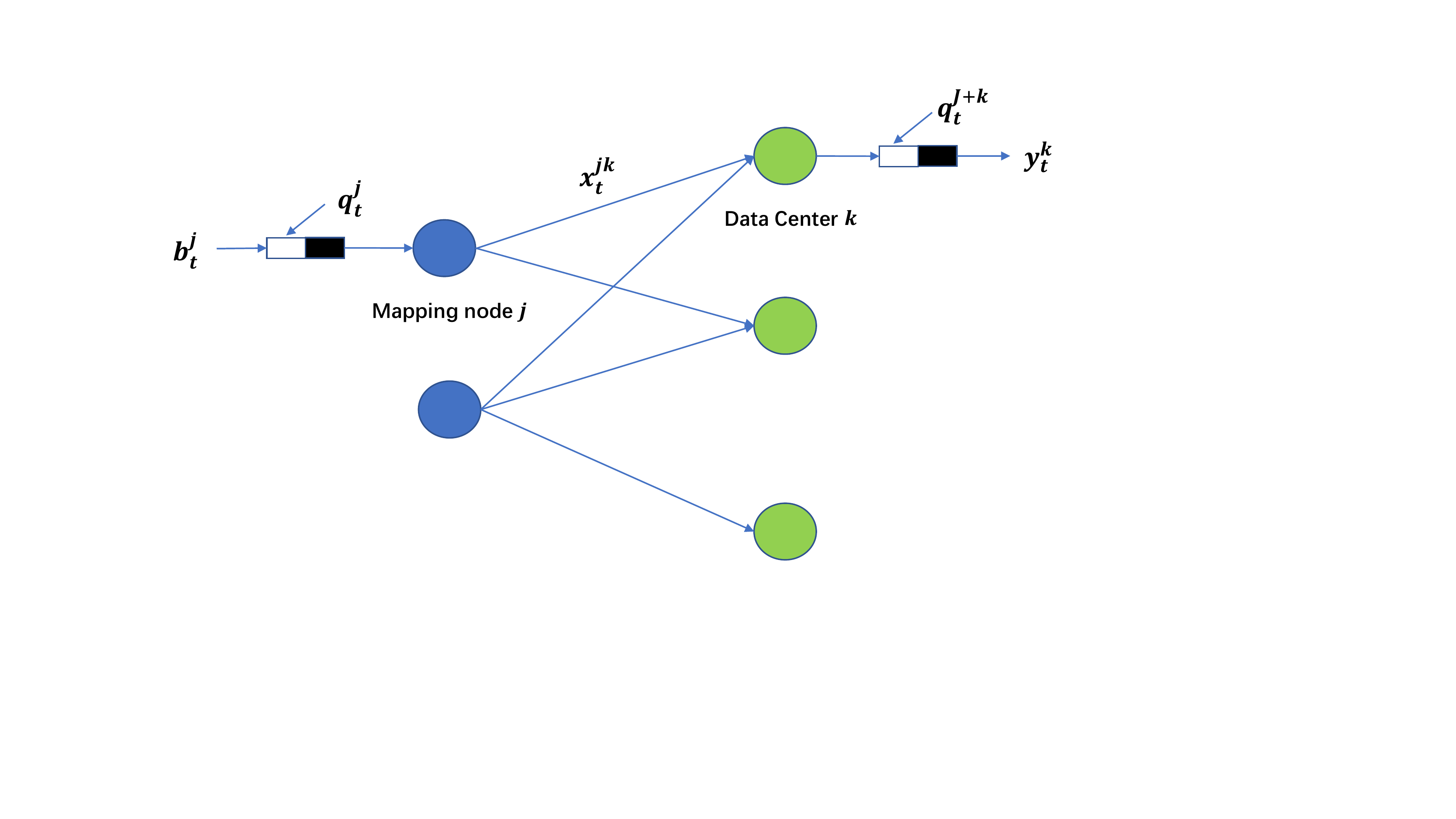}}
\caption{The overview of the system.}
\end{figure}
\subsection{Online Network Resource Allocation}
Within this subsection, we consider an online resource allocation problem over a cloud network \cite{chen2016dglb,chen2017an}. The network consists of mapping nodes $\mathcal{J}=\{1,...,J\}$ and data centers $\mathcal{K}=\{1,...,K\}$. We use a directed graph $\mathcal{G}=(\mathcal{I},\mathcal{\epsilon} )$ to represent it, where $\mathcal{I}=\mathcal{J}\cup \mathcal{K} $, $|\mathcal{I}|=J+K$, and $|\mathcal{\epsilon }|=E$. $\mathcal{\epsilon}$ includes all the links which connect mapping nodes with data centers, and the "virtual" exogenous edges coming out of the data centers.
At each time slot $t$, each mapping node $j$ receives a data request 
${b}_{t}^{j}$ from exogenous user, and schedules ${x}_{t}^{jk}$ workload to data center $k$. Each data center $k$ serves workload ${y}_{t}^{k}$ based on its source availability. We assume each node (including data center and mapping node) could buffer the unserved workloads into its local queue. Next we describe the workflow of the overall system, which is illustrated in Figure 3. Specifically,
at each time slot $t$, mapping node $j$ has an exogenous workload ${b}_{t}^{j}$ plus that stored in its local queue ${q}_{t}^{j}$, then it schedules workload ${x}_{t}^{jk}$ to data center $k$. Data center $k$ has a received workload of the amount of $\sum _{ j=1 }^{ J }{ { x }_{ t }^{ jk } }$ plus that stored in its local queue ${q}_{t}^{J+k}$, and 
serves an amount of workload ${y}_{t}^{k}$.

We define a resource allocation vector ${ \bm{x} }_{ t }=[{ x }_{ t }^{ 11 },...,{ x }_{ t }^{ JK },{ y }_{ t }^{ 1 },...,{ y }_{ t }^{ K }{ ] }^{ T } \in { R }_{ + }^{ E }$, and load arrival vector ${ \bm{b} }_{ t }=[{ b }_{ t }^{ 1 },...,{ b }_{ t }^{ J },0,...,0{ ] }^{ T }$ $\in { R }_{ + }^{ I }$
to represent the exogenous load arrival rates of all nodes at time slot $t$.
We also define $I\times E$ node-incidence matrix $\bm{A}$, where $(i,j)$-th entry ${\bm{A}}_{(i,j)}=1$ if link $j$ enters node $i$, or ${\bm{A}}_{(i,j)}=-1$ if link $j$ leaves node $i$, otherwise ${\bm{A}}_{(i,j)}=0$.
Hence the vector $\bm{A}{\bm{x}}_{t}+{ \bm{b} }_{ t }$ represents the aggregate workloads of all nodes. There is service residual at node $i$ if $(\bm{A}{\bm{x}}_{t}+{ \bm{b} }_{ t }{)}_{i}>0$, otherwise the current load of node $i$ exceeds its service capacity. 
%We assume each node (including data center and mapping node) could buffer the unserved workloads into its local queue. 
At each time slot $t$, the queue length vector of all nodes is given by ${ \bm{q} }_{ t }=[{ q }_{ t }^{ 1 },...,{ q }_{ t }^{ J+K }{ ] }^{ T }$, and its update rule of ${\bm{q}}_{t}$ is ${ \bm{q} }_{ t+1 }=({ \bm{q} }_{ t }+\bm{A}{ \bm{x} }_{ t }+{ \bm{b} }_{ t }{ ) }^{ + }$. We denote ${B}_{jk}$ be the maximum bandwidth of link $(j,k)$, and ${C}_{k}$ be the resource capacity of data center $k$. Thus the feasible set is $\chi =\{ x|0\le x\le c\} $, where $ c=[{ B }_{ 11 },...,{ B }_{ JK },{ C }_{ 1 },...,{ C }_{ K }]$.

Here we formulate the accumulated cost of the overall system. We divide it into two parts, the one is power cost, the other is bandwidth cost. The power cost characterizes the energy price and renewable generation, and the bandwidth cost characterizes the transmission delay.
The power cost of each data center $k$ at time slot $t$ is ${ f }_{ t,k }({ y }_{ t }^{ k })$. The bandwidth cost of link $(j,k)$ is ${f}_{t,jk}({x}_{t}^{jk})$. Both of them are unknown before the resource allocation at time slot $t$.
Hence at each time slot $t$, the instantaneous cost of the overall system is
\begin{equation}
     { f }_{ t }({ \bm{x} }_{ t })=\sum _{ k\in \mathcal{K} }{ { f }_{ t,k }({ y }_{ t }^{ k }) } +\sum _{ j\in \mathcal{J} }{ \sum _{ k\in \mathcal{K} }{ { f }_{ t,jk }({ x }_{ t }^{ jk }) }  } 
\end{equation}
Our goal is to minimize the accumulated cost of the overall system while ensuring all workloads are served, shown in the following optimization problem ${\bm{P}}_{1}$:
\begin{equation}
    \begin{split}
   {\bm{P}}_{1}\bm{:} \quad &\min _{ \{ {\bm{x}}_{t}\in \chi\}  }{ \sum _{ t=1 }^{ T }{ { f }_{ t }({ x }_{ t }) }  } \\ 
        &s.t.\quad { \bm{q} }_{ t+1 }=[{ \bm{q} }_{ t }+\bm{A}{ \bm{x} }_{ t }+{ \bm{b} }_{ t }{ ] }^{ + },\forall t,\\ 
        &\quad \quad \quad { \bm{q} }_{ 1 }={ \bm{q} }_{ T+1 }=0,
    \end{split}
\end{equation}
where the initial queue length is given by ${\bm{q}}_{1}$, and ${ \bm{q} }_{ T+1 }=0$ implies that all workloads should be served before the end of scheduling horizon $T$. However, solving ${P}_{1}$ is generally challenging using traditional methods since the future workload arrivals are not known a prior. Therefore, we could relax the first constraint in ${\bm{P}}_{1}$ as follows:
\begin{equation}
    { \bm{q} }_{ T+1 }\ge { \bm{q} }_{ T }+\bm{A}{ \bm{x} }_{ T }+{ \bm{b} }_{ T }\ge...\ge { \bm{q} }_{ 1 }+\sum _{ t=1 }^{ T }{ (\bm{A}{ \bm{x} }_{ t }+{ \bm{b} }_{ t }) } \Rightarrow \sum _{ t=1 }^{ T }{ (\bm{A}{ \bm{x} }_{ t }+{ \bm{b} }_{ t }) } \le { \bm{q} }_{ T+1 }-{ \bm{q} }_{ 1 }=0
\end{equation}
Then we transform ${\bm{P}}_{1}$ into the following optimization problem ${\bm{P}}_{2}$:
%\begin{equation} \begin{split} &\min _{ \{ { x }_{ t }\in \chi \}  }{ \sum _{ t=1 }^{ T }{ { f }_{ t }({ x }_{ t }) }  } \\  &s.t.\quad \sum _{ t=1 }^{ T }{ (A{ x }_{ t }+{ b }_{ t }) } \le 0, \end{split} \end{equation}
\begin{equation}
     {\bm{P}}_{2}\bm{:} \quad \min _{ \{ { x }_{ t }\in \chi \}  }{ \sum _{ t=1 }^{ T }{ { f }_{ t }({ x }_{ t }) }  } \quad s.t.\quad \sum _{ t=1 }^{ T }{ (\bm{A}{ x }_{ t }+{ \bm{b} }_{ t }) } \le 0,
\end{equation}
which could be solved by our framework of OCO with long-term and time-varying constraints.
\subsection{Online Job Scheduling}
Within this subsection, we consider an online job scheduling problem \cite{im2016scheduling,im2016competitive,liu2019online}, in which the computing cluster consists of multiple servers with heterogeneous computation resources. Specifically, consider a computing cluster, and it consists of $M$ servers, which indexed from $1$ to $M$. We assume server $i$ has ${C}_{i}$ CPU cores and can process multiple jobs simultaneously unless the total demand of its execution jobs exceeds ${C}_{i}$. Time is slotted and job $j$ arrives the cluster at time slot ${a}_{j}$. The total number of jobs is $N$. Each incoming job joins a global queue managed by a scheduler, waiting to be assigned to an available server(s) for execution in the subsequent time slots. At the beginning of each time slot, the scheduler has to decide which job(s) to schedule and which sever(s) assigned to it(them).
%together with their server assignment. 
We assume job $j$ requires ${d}_{j}$ CPU cores to run and ${p}_{j}$
units of time to finish when its demand of ${d}_{j}$ CPU cores are fully satisfied, where both ${d}_{j}$ and ${p}_{j}$ are integers and will be reported to the scheduler once job $j$ arrives at the cluster. Thus
the quantity ${v}_{j}={p}_{j}{d}_{j}$ is the volume of job $j$. %We also assume preemption and migration are allowed, i.e., a running job can be check-pointed, preempted, and resumed later on the same or a different server.
We also assume that preemption and migration are allowed, i.e., a running job can be check-pointed, preempted, and then recovered on the same server or on a different server.

We denote ${u}_{j}^{i}(t)$ 
be the number of CPU cores allocated to job $j$ on server $i$ at time-slot $t$. By time-division-multiplexing of CPU cores, any job could also be processed even if it is allocated fewer than ${d}_{j}$ CPU cores, but it needs to take more than ${p}_{j}$ time-slots to finish. We say job $j$ finishes when its completion time  ${c}_{j}$ satisfies $\sum _{ t={ a }_{ j }+1 }^{ { c }_{ j } }{ \sum _{ i=1 }^{ M }{ { u }_{ j }^{ i }(t) }  } \ge { d }_{ j }{ p }_{ j }$, that is, its volume is fully served, and the flowtime of job $j$ is ${c}_{j}-{a}_{j}$. We assume that any job $j$ cannot benefit from the extra number of cores (i.e., it is allocated more than ${d}_{j}$ cores). To avoid the waste of resources, we have $\sum _{ i=1 }^{ M }{ { u }_{ j }^{ i }(t) } \le { d }_{ j }$. The online scheduler strikes a balance between fairness and job latency. Hence, like \cite{liu2019online}, we adopt the ${l}_{k}$ norm of job flowtime \cite{liu2019online} to represent job's "cost". Then our goal is to minimize the sum of ${l}_{k}$ norm of all jobs' flowtime while satisfying some constraints, shown in the
following optimization problem ${\bm{P}}_{3}$:
\begin{equation}
    \begin{split}
        {\bm{P}}_{3}\bm{:} \quad &\min _{ \{ { u }_{ j }^{ i }(t)\}  }{ \sum _{ j=1 }^{ N }{ ({ c }_{ j }-{ a }_{ j }{ ) }^{ k } }  } \\ & s.t.\quad \sum _{ t={ a }_{ j }+1 }^{ { c }_{ j } }{ \sum _{ i=1 }^{ M }{ { u }_{ j }^{ i }(t) }  } \ge { d }_{ j }{ p }_{ j },\forall j,\\ & \quad \quad \sum _{ i }{ { u }_{ j }^{ i }(t) } \le { d }_{ j },\forall j,t,\\ \quad \quad &\quad \quad \sum _{ j:t>{ a }_{ j } }{ { u }_{ j }^{ i }(t) } \le { C }_{ i },\forall i,t,\\
        &\quad \quad \; {u}_{j}^{i}(t)\in N, 
    \end{split}
\end{equation}
where the third constraint means the total number of allocated CPU
cores on server $i$ cannot exceed its capacity ${C}_{i}$ at any time slot. However, it can be verified that ${\bm{P}}_{3}$ is NP-hard since it is an integer programming problem. Hence here we could adopt the approximation algorithm to solve ${\bm{P}}_{3}$. Specifically, we approximate ${l}_{k}$ norm of flowtime $({c}_{j}-{a}_{j}{)}^{k}$ with its fractional job flowtime counterpart \cite{im2015temporal}, that is,
$$({ c }_{ j }-{ a }_{ j }{ ) }^{ k }\approx \sum _{ t={ a }_{ j }+1 }{ \sum _{ i=1 }^{ M }{ { ({ (t-{ a }_{ j }{ ) }^{ k } }/{ { p }_{ j } }+{ p }_{ j }^{ k-1 }){ u }_{ j }^{ i }(t) }/{ { d }_{ j } } }  }. $$
We define ${ y }_{ j }(t)={ \sum _{ i }{ { u }_{ j }^{ i }(t) }  }/{ { d }_{ j }(t) }$ be the total CPU cores rates allocated to
job $j$ at time slot $t$. We could transform ${\bm{P}}_{3}$ into the following optimization problem ${\bm{P}}_{4}$:
\begin{align*}
    \begin{split}
        {\bm{P}}_{4}\bm{:} \quad \min _{ \{ { u }_{ j }^{ i }(t)\}  }{ \sum _{ t=1 } \sum _{ j:t\ge { a }_{ j }+1 }{ ({ (t-{ a }_{ j }{ ) }^{ k } }/{ { p }_{ j } }+{ p }_{ j }^{ k-1 }){ y }_{ j }(t) }  }
    \end{split}
\end{align*}
\begin{align}
    \begin{split}
         &s.t.\quad \sum _{ t={ a }_{ j }+1 }{ { y }_{ j }(t) } \ge { p }_{ j },\forall j,\\ 
        & \quad \quad \quad 0\le { y }_{ j }(t)\le 1,\forall j,t,\\  &\quad \quad \quad \sum _{ j:t>{ a }_{ j } }{ { d }_{ j }{ y }_{ j }(t) } \le \sum _{ i=1 }^{ M }{ { C }_{ i } } ,\forall t,\\
        & \quad \quad \quad {d}_{j}{y}_{j}(t)\in N,\forall j,t.
    \end{split}
\end{align}
 We denote by ${OPT}_{{\bm{P}}_{3}}^{*}$ and ${OPT}_{{\bm{P}}_{4}}^{*}$ the optimal objective values of optimization problem ${\bm{P}}_{3}$ and ${\bm{P}}_{4}$ respectively, then by using the same argument as \cite{im2015temporal}, we have ${OPT}_{{\bm{P}}_{4}}^{*}\le 2{OPT}_{{\bm{P}}_{3}}^{*}$. Next we show that we could use the framework of OCO with long-term and time-varying constraints to solve problem ${\bm{P}}_{4}$. 
 
 \textbf{Solve ${\bm{P}}_{4}$ using the framework of OCO with long-term and time-varying constraints.} We formulate the third constraint in ${\bm{P}}_{4}$ as the short-term constraint which needs to be satisfied strictly at each time slot. We also notice that the first constraint could be formulated as the long-term constraint. Thus, we separate the first constraint in ${\bm{P}}_{4}$ into each time slot constraint:
 \begin{equation}
 \label{per slot time constraint}
     { g }_{ t,j }({ y }_{ j }(t))=\frac { { p }_{ j } }{ T-{ a }_{ j } } -{ y }_{ j }(t)\le 0,
 \end{equation}
 where $T$ is the predicted completion time for all jobs, which could be known or estimated ahead of time in many scenarios. The per time slot constraint \eqref{per slot time constraint} could be violated in some time slots but the accumulated constraint violations should be controlled. We relax the integer constraint of ${y}_{j}(t)$ and define the feasible set of it as:
 $$\chi (t)=\{ \bm{y}|\sum _{ j:t>{ a }_{ j } }{ { d }_{ j }{ y }_{ j } } \le \sum _{ i=1 }^{ M }{ { C }_{ i } } ,0\le { y }_{ j }\le 1.\}.$$
Indeed, we could still make the resultant online algorithm satisfies the integer constraint in the sequel. Therefore, the optimization problem ${\bm{P}}_{4}$ could be transformed into the following optimization problem ${\bm{P}}_{5}$:
 \begin{equation}  
 \begin{split} 
  {\bm{P}}_{5}\bm{:} \quad&\min _{ \bm{y}(t)\in \chi (t) }{ {f}_{t}(\bm{y}(t))=\sum _{ t=1 }^{ T } \sum _{ j:t\ge { a }_{ j }+1 }{ ({ (t-{ a }_{ j }{ ) }^{ k } }/{ { p }_{ j } }+{ p }_{ j }^{ k-1 }){ y }_{ j }(t) }  } \\ &\quad s.t.\quad \sum _{ t={ a }_{ j }+1 }^{ T }{ { g }_{ t,j }({ y }_{ j }(t)) } \le 0,\forall j,  
 \end{split} 
 \end{equation}
 which could be solved by our framework of OCO with long-term and time-varying constraints.
 As stated before, although feasible set $\chi$ is a time-varying
set, our algorithm is valid and the corresponding theoretical results also hold.
\section{Conclusion and future work}
% the setting of parameters are omitted due to space limitation. see XX for details.
In this paper, we develop and analyze a novel algorithm for OCO with long term and time-varying constraints. To the best of our knowledge, our algorithm is the first parameter-free algorithm to simultaneously achieve sublinear dynamic regret and violation under common assumptions. We then extend our algorithm and analysis to some practical cases.
%(copy from other works) In light of this current work and the recent work(cite), it is interesting if the combination of smoothness and strong convexity can lead to a better performance bounds. 
%And it would be in interesting to derive lower bounds on the dynamic regret or constraint violations of such algorithms, perhaps with online variants of standard constructions used in the offline setting(see for instance [2])
%
For future work, It is a good direction to investigate sharper performance bounds for OCO with long-term and time-varying constraints. Moreover, whether incorporating other properties, like strong convexity and smoothness, can lead to better performance bounds is also an open question. 
%In the future, we will investigate sharper performance bounds for OCO with long-term and time-varying constraints. Moreover, we will study whether incorporating other properties, like strong convexity and smoothness, can lead to better performance bounds is also an open question. 
%consider incorporating other properties, like strong convexity and smoothness, to further enhance the dynamic regret and the constraint violations

\clearpage
\bibliographystyle{ACM-Reference-Format}
\bibliography{sample-base}


\begin{thebibliography}{32}

%%% the Web address.

\ifx \showCODEN    \undefined \def \showCODEN     #1{\unskip}     \fi
\ifx \showDOI      \undefined \def \showDOI       #1{#1}\fi
\ifx \showISBNx    \undefined \def \showISBNx     #1{\unskip}     \fi
\ifx \showISBNxiii \undefined \def \showISBNxiii  #1{\unskip}     \fi
\ifx \showISSN     \undefined \def \showISSN      #1{\unskip}     \fi
\ifx \showLCCN     \undefined \def \showLCCN      #1{\unskip}     \fi
\ifx \shownote     \undefined \def \shownote      #1{#1}          \fi
\ifx \showarticletitle \undefined \def \showarticletitle #1{#1}   \fi
\ifx \showURL      \undefined \def \showURL       {\relax}        \fi
% invisible to TeX
\providecommand\bibfield[2]{#2}
\providecommand\bibinfo[2]{#2}
\providecommand\natexlab[1]{#1}
\providecommand\showeprint[2][]{arXiv:#2}

\bibitem[\protect\citeauthoryear{Amiri}{Amiri}{2019}]%
        {amiri2019towards}
\bibfield{author}{\bibinfo{person}{Maryam Amiri}.}
  \bibinfo{year}{2019}\natexlab{}.
\newblock \emph{\bibinfo{title}{Towards enhancing QoE for software defined
  networks based cloud gaming services}}.
\newblock \bibinfo{thesistype}{Ph.D. Dissertation}.
  \bibinfo{school}{Universit{\'e} d'Ottawa/University of Ottawa}.
\newblock


\bibitem[\protect\citeauthoryear{Arce and Salinas}{Arce and Salinas}{2012}]%
        {arce2012online}
\bibfield{author}{\bibinfo{person}{Paola Arce} {and} \bibinfo{person}{Luis
  Salinas}.} \bibinfo{year}{2012}\natexlab{}.
\newblock \showarticletitle{Online ridge regression method using sliding
  windows}. In \bibinfo{booktitle}{\emph{2012 31st International Conference of
  the Chilean Computer Science Society}}. IEEE, \bibinfo{pages}{87--90}.
\newblock


\bibitem[\protect\citeauthoryear{Cao and Liu}{Cao and Liu}{2018}]%
        {cao2018online}
\bibfield{author}{\bibinfo{person}{Xuanyu Cao} {and} \bibinfo{person}{KJ~Ray
  Liu}.} \bibinfo{year}{2018}\natexlab{}.
\newblock \showarticletitle{Online convex optimization with time-varying
  constraints and bandit feedback}.
\newblock \bibinfo{journal}{\emph{IEEE Trans. Automat. Control}}
  \bibinfo{volume}{64}, \bibinfo{number}{7} (\bibinfo{year}{2018}),
  \bibinfo{pages}{2665--2680}.
\newblock


\bibitem[\protect\citeauthoryear{Chen and Giannakis}{Chen and
  Giannakis}{2019}]%
        {chen2019bandit}
\bibfield{author}{\bibinfo{person}{Tianyi Chen} {and}
  \bibinfo{person}{Georgios~B Giannakis}.} \bibinfo{year}{2019}\natexlab{}.
\newblock \showarticletitle{Bandit Convex Optimization for Scalable and Dynamic
  IoT Management}.
\newblock \bibinfo{journal}{\emph{IEEE INTERNET OF THINGS JOURNAL}}
  \bibinfo{volume}{6}, \bibinfo{number}{1} (\bibinfo{year}{2019}).
\newblock


\bibitem[\protect\citeauthoryear{Chen, Ling, and Giannakis}{Chen
  et~al\mbox{.}}{2017}]%
        {chen2017an}
\bibfield{author}{\bibinfo{person}{Tianyi Chen}, \bibinfo{person}{Qing Ling},
  {and} \bibinfo{person}{Georgios~B Giannakis}.}
  \bibinfo{year}{2017}\natexlab{}.
\newblock \showarticletitle{An Online Convex Optimization Approach to Proactive
  Network Resource Allocation}.
\newblock \bibinfo{journal}{\emph{IEEE Transactions on Signal Processing}}
  \bibinfo{volume}{65}, \bibinfo{number}{24} (\bibinfo{year}{2017}),
  \bibinfo{pages}{6350--6364}.
\newblock


\bibitem[\protect\citeauthoryear{Chen, Ling, Shen, and Giannakis}{Chen
  et~al\mbox{.}}{2018}]%
        {chen2018heterogeneous}
\bibfield{author}{\bibinfo{person}{Tianyi Chen}, \bibinfo{person}{Qing Ling},
  \bibinfo{person}{Yanning Shen}, {and} \bibinfo{person}{Georgios~B
  Giannakis}.} \bibinfo{year}{2018}\natexlab{}.
\newblock \showarticletitle{Heterogeneous Online Learning for
  “Thing-Adaptive” Fog Computing in IoT}.
\newblock \bibinfo{journal}{\emph{IEEE Internet of Things Journal}}
  \bibinfo{volume}{5}, \bibinfo{number}{6} (\bibinfo{year}{2018}),
  \bibinfo{pages}{4328--4341}.
\newblock


\bibitem[\protect\citeauthoryear{Chen, Marques, and Giannakis}{Chen
  et~al\mbox{.}}{2016}]%
        {chen2016dglb}
\bibfield{author}{\bibinfo{person}{Tianyi Chen}, \bibinfo{person}{Antonio~G
  Marques}, {and} \bibinfo{person}{Georgios~B Giannakis}.}
  \bibinfo{year}{2016}\natexlab{}.
\newblock \showarticletitle{DGLB: Distributed stochastic geographical load
  balancing over cloud networks}.
\newblock \bibinfo{journal}{\emph{IEEE Transactions on Parallel and Distributed
  Systems}} \bibinfo{volume}{28}, \bibinfo{number}{7} (\bibinfo{year}{2016}),
  \bibinfo{pages}{1866--1880}.
\newblock


\bibitem[\protect\citeauthoryear{Fang, Harvey, Portella, and Friedlander}{Fang
  et~al\mbox{.}}{2020}]%
        {fang2020online}
\bibfield{author}{\bibinfo{person}{Huang Fang}, \bibinfo{person}{Nicholas~JA
  Harvey}, \bibinfo{person}{Victor~S Portella}, {and}
  \bibinfo{person}{Michael~P Friedlander}.} \bibinfo{year}{2020}\natexlab{}.
\newblock \showarticletitle{Online mirror descent and dual averaging: keeping
  pace in the dynamic case}.
\newblock \bibinfo{journal}{\emph{arXiv preprint arXiv:2006.02585}}
  (\bibinfo{year}{2020}).
\newblock


\bibitem[\protect\citeauthoryear{Haworth, Shawe-Taylor, Cheng, and
  Wang}{Haworth et~al\mbox{.}}{2014}]%
        {haworth2014local}
\bibfield{author}{\bibinfo{person}{James Haworth}, \bibinfo{person}{John
  Shawe-Taylor}, \bibinfo{person}{Tao Cheng}, {and} \bibinfo{person}{Jiaqiu
  Wang}.} \bibinfo{year}{2014}\natexlab{}.
\newblock \showarticletitle{Local online kernel ridge regression for
  forecasting of urban travel times}.
\newblock \bibinfo{journal}{\emph{Transportation research part C: emerging
  technologies}}  \bibinfo{volume}{46} (\bibinfo{year}{2014}),
  \bibinfo{pages}{151--178}.
\newblock


\bibitem[\protect\citeauthoryear{Hazan}{Hazan}{2019}]%
        {hazan2019introduction}
\bibfield{author}{\bibinfo{person}{Elad Hazan}.}
  \bibinfo{year}{2019}\natexlab{}.
\newblock \showarticletitle{Introduction to online convex optimization}.
\newblock \bibinfo{journal}{\emph{arXiv preprint arXiv:1909.05207}}
  (\bibinfo{year}{2019}).
\newblock


\bibitem[\protect\citeauthoryear{Huang, Liu, and Hao}{Huang
  et~al\mbox{.}}{2014}]%
        {huang2014power}
\bibfield{author}{\bibinfo{person}{Longbo Huang}, \bibinfo{person}{Xin Liu},
  {and} \bibinfo{person}{Xiaohong Hao}.} \bibinfo{year}{2014}\natexlab{}.
\newblock \showarticletitle{The power of online learning in stochastic network
  optimization}. In \bibinfo{booktitle}{\emph{The 2014 ACM international
  conference on Measurement and modeling of computer systems}}.
  \bibinfo{pages}{153--165}.
\newblock


\bibitem[\protect\citeauthoryear{Huang, Moeller, Neely, and
  Krishnamachari}{Huang et~al\mbox{.}}{2012}]%
        {huang2012lifo}
\bibfield{author}{\bibinfo{person}{Longbo Huang}, \bibinfo{person}{Scott
  Moeller}, \bibinfo{person}{Michael~J Neely}, {and} \bibinfo{person}{Bhaskar
  Krishnamachari}.} \bibinfo{year}{2012}\natexlab{}.
\newblock \showarticletitle{LIFO-backpressure achieves near-optimal
  utility-delay tradeoff}.
\newblock \bibinfo{journal}{\emph{IEEE/ACM Transactions On Networking}}
  \bibinfo{volume}{21}, \bibinfo{number}{3} (\bibinfo{year}{2012}),
  \bibinfo{pages}{831--844}.
\newblock


\bibitem[\protect\citeauthoryear{Huang and Neely}{Huang and Neely}{2011}]%
        {huang2011utility}
\bibfield{author}{\bibinfo{person}{Longbo Huang} {and}
  \bibinfo{person}{Michael~J Neely}.} \bibinfo{year}{2011}\natexlab{}.
\newblock \showarticletitle{Utility optimal scheduling in processing networks}.
\newblock \bibinfo{journal}{\emph{Performance Evaluation}}
  \bibinfo{volume}{68}, \bibinfo{number}{11} (\bibinfo{year}{2011}),
  \bibinfo{pages}{1002--1021}.
\newblock


\bibitem[\protect\citeauthoryear{Im, Kulkarni, and Moseley}{Im
  et~al\mbox{.}}{2015}]%
        {im2015temporal}
\bibfield{author}{\bibinfo{person}{Sungjin Im}, \bibinfo{person}{Janardhan
  Kulkarni}, {and} \bibinfo{person}{Benjamin Moseley}.}
  \bibinfo{year}{2015}\natexlab{}.
\newblock \showarticletitle{Temporal fairness of round robin: Competitive
  analysis for lk-norms of flow time}. In \bibinfo{booktitle}{\emph{Proceedings
  of the 27th ACM symposium on Parallelism in Algorithms and Architectures}}.
  \bibinfo{pages}{155--160}.
\newblock


\bibitem[\protect\citeauthoryear{Im, Kulkarni, Moseley, and Munagala}{Im
  et~al\mbox{.}}{2016a}]%
        {im2016competitive}
\bibfield{author}{\bibinfo{person}{Sungjin Im}, \bibinfo{person}{Janardhan
  Kulkarni}, \bibinfo{person}{Benjamin Moseley}, {and} \bibinfo{person}{Kamesh
  Munagala}.} \bibinfo{year}{2016}\natexlab{a}.
\newblock \showarticletitle{A competitive flow time algorithm for heterogeneous
  clusters under polytope constraints}. In
  \bibinfo{booktitle}{\emph{Approximation, Randomization, and Combinatorial
  Optimization. Algorithms and Techniques (APPROX/RANDOM 2016)}}. Schloss
  Dagstuhl-Leibniz-Zentrum fuer Informatik.
\newblock


\bibitem[\protect\citeauthoryear{Im, Naghshnejad, and Singhal}{Im
  et~al\mbox{.}}{2016b}]%
        {im2016scheduling}
\bibfield{author}{\bibinfo{person}{Sungjin Im}, \bibinfo{person}{Mina
  Naghshnejad}, {and} \bibinfo{person}{Mukesh Singhal}.}
  \bibinfo{year}{2016}\natexlab{b}.
\newblock \showarticletitle{Scheduling jobs with non-uniform demands on
  multiple servers without interruption}. In \bibinfo{booktitle}{\emph{IEEE
  INFOCOM 2016-The 35th Annual IEEE International Conference on Computer
  Communications}}. IEEE, \bibinfo{pages}{1--9}.
\newblock


\bibitem[\protect\citeauthoryear{Jenatton, Huang, and Archambeau}{Jenatton
  et~al\mbox{.}}{2016}]%
        {jenatton2016adaptive}
\bibfield{author}{\bibinfo{person}{Rodolphe Jenatton}, \bibinfo{person}{Jim
  Huang}, {and} \bibinfo{person}{C{\'e}dric Archambeau}.}
  \bibinfo{year}{2016}\natexlab{}.
\newblock \showarticletitle{Adaptive algorithms for online convex optimization
  with long-term constraints}. In \bibinfo{booktitle}{\emph{International
  Conference on Machine Learning}}. \bibinfo{pages}{402--411}.
\newblock


\bibitem[\protect\citeauthoryear{Liu, Xu, and Lau}{Liu et~al\mbox{.}}{2019}]%
        {liu2019online}
\bibfield{author}{\bibinfo{person}{Yang Liu}, \bibinfo{person}{Huanle Xu},
  {and} \bibinfo{person}{Wing~Cheong Lau}.} \bibinfo{year}{2019}\natexlab{}.
\newblock \showarticletitle{Online job scheduling with resource packing on a
  cluster of heterogeneous servers}. In \bibinfo{booktitle}{\emph{IEEE INFOCOM
  2019-IEEE Conference on Computer Communications}}. IEEE,
  \bibinfo{pages}{1441--1449}.
\newblock


\bibitem[\protect\citeauthoryear{Liu, Liu, Low, and Wierman}{Liu
  et~al\mbox{.}}{2014}]%
        {liu2014pricing}
\bibfield{author}{\bibinfo{person}{Zhenhua Liu}, \bibinfo{person}{Iris Liu},
  \bibinfo{person}{Steven Low}, {and} \bibinfo{person}{Adam Wierman}.}
  \bibinfo{year}{2014}\natexlab{}.
\newblock \showarticletitle{Pricing data center demand response}.
\newblock \bibinfo{journal}{\emph{ACM SIGMETRICS Performance Evaluation
  Review}} \bibinfo{volume}{42}, \bibinfo{number}{1} (\bibinfo{year}{2014}),
  \bibinfo{pages}{111--123}.
\newblock


\bibitem[\protect\citeauthoryear{Mahdavi, Jin, and Yang}{Mahdavi
  et~al\mbox{.}}{2012}]%
        {mahdavi2012trading}
\bibfield{author}{\bibinfo{person}{Mehrdad Mahdavi}, \bibinfo{person}{Rong
  Jin}, {and} \bibinfo{person}{Tianbao Yang}.} \bibinfo{year}{2012}\natexlab{}.
\newblock \showarticletitle{Trading regret for efficiency: online convex
  optimization with long term constraints}.
\newblock \bibinfo{journal}{\emph{The Journal of Machine Learning Research}}
  \bibinfo{volume}{13}, \bibinfo{number}{1} (\bibinfo{year}{2012}),
  \bibinfo{pages}{2503--2528}.
\newblock


\bibitem[\protect\citeauthoryear{Qiu and Wei}{Qiu and Wei}{2020}]%
        {qiu2020beyond}
\bibfield{author}{\bibinfo{person}{Shuang Qiu} {and} \bibinfo{person}{Xiaohan
  Wei}.} \bibinfo{year}{2020}\natexlab{}.
\newblock \showarticletitle{Beyond O($\sqrt{T}$) Regret for Constrained Online
  Optimization: Gradual Variations and Mirror Prox}.
\newblock \bibinfo{journal}{\emph{arXiv preprint arXiv:2006.12455}}
  (\bibinfo{year}{2020}).
\newblock


\bibitem[\protect\citeauthoryear{Sharma, Khanduri, Shen, Bucci~Jr, and
  Varshney}{Sharma et~al\mbox{.}}{2020}]%
        {sharma2020distributed}
\bibfield{author}{\bibinfo{person}{Pranay Sharma}, \bibinfo{person}{Prashant
  Khanduri}, \bibinfo{person}{Lixin Shen}, \bibinfo{person}{Donald~J Bucci~Jr},
  {and} \bibinfo{person}{Pramod~K Varshney}.} \bibinfo{year}{2020}\natexlab{}.
\newblock \showarticletitle{On distributed online convex optimization with
  sublinear dynamic regret and fit}.
\newblock \bibinfo{journal}{\emph{arXiv preprint arXiv:2001.03166}}
  (\bibinfo{year}{2020}).
\newblock


\bibitem[\protect\citeauthoryear{Xu, Liu, Lau, Guo, and Liu}{Xu
  et~al\mbox{.}}{2019}]%
        {xu2019efficient}
\bibfield{author}{\bibinfo{person}{Huanle Xu}, \bibinfo{person}{Yang Liu},
  \bibinfo{person}{Wing~Cheong Lau}, \bibinfo{person}{Jun Guo}, {and}
  \bibinfo{person}{Alex Liu}.} \bibinfo{year}{2019}\natexlab{}.
\newblock \showarticletitle{Efficient online resource allocation in
  heterogeneous clusters with machine variability}. In
  \bibinfo{booktitle}{\emph{IEEE INFOCOM 2019-IEEE Conference on Computer
  Communications}}. IEEE, \bibinfo{pages}{478--486}.
\newblock


\bibitem[\protect\citeauthoryear{Yi, Li, Yang, Xie, Chai, and Johansson}{Yi
  et~al\mbox{.}}{2020}]%
        {yi2020distributed}
\bibfield{author}{\bibinfo{person}{Xinlei Yi}, \bibinfo{person}{Xiuxian Li},
  \bibinfo{person}{Tao Yang}, \bibinfo{person}{Lihua Xie},
  \bibinfo{person}{Tianyou Chai}, {and} \bibinfo{person}{Karl~H Johansson}.}
  \bibinfo{year}{2020}\natexlab{}.
\newblock \showarticletitle{Distributed bandit online convex optimization with
  time-varying coupled inequality constraints}.
\newblock \bibinfo{journal}{\emph{IEEE Trans. Automat. Control}}
  (\bibinfo{year}{2020}).
\newblock


\bibitem[\protect\citeauthoryear{Yu, Neely, and Wei}{Yu et~al\mbox{.}}{2017}]%
        {yu2017online}
\bibfield{author}{\bibinfo{person}{Hao Yu}, \bibinfo{person}{Michael Neely},
  {and} \bibinfo{person}{Xiaohan Wei}.} \bibinfo{year}{2017}\natexlab{}.
\newblock \showarticletitle{Online convex optimization with stochastic
  constraints}. In \bibinfo{booktitle}{\emph{Advances in Neural Information
  Processing Systems}}. \bibinfo{pages}{1428--1438}.
\newblock


\bibitem[\protect\citeauthoryear{Yu and Neely}{Yu and Neely}{2020}]%
        {yu2020a}
\bibfield{author}{\bibinfo{person}{Hao Yu} {and} \bibinfo{person}{Michael~J
  Neely}.} \bibinfo{year}{2020}\natexlab{}.
\newblock \showarticletitle{A Low Complexity Algorithm with O($\sqrt{T}$)
  Regret and O(1) Constraint Violations for Online Convex Optimization with
  Long Term Constraints}.
\newblock \bibinfo{journal}{\emph{Journal of Machine Learning Research}}
  \bibinfo{volume}{21}, \bibinfo{number}{1} (\bibinfo{year}{2020}),
  \bibinfo{pages}{1--24}.
\newblock


\bibitem[\protect\citeauthoryear{Yuan and Lamperski}{Yuan and
  Lamperski}{2018}]%
        {yuan2018online}
\bibfield{author}{\bibinfo{person}{Jianjun Yuan} {and} \bibinfo{person}{Andrew
  Lamperski}.} \bibinfo{year}{2018}\natexlab{}.
\newblock \showarticletitle{Online convex optimization for cumulative
  constraints}. In \bibinfo{booktitle}{\emph{Advances in Neural Information
  Processing Systems}}. \bibinfo{pages}{6137--6146}.
\newblock


\bibitem[\protect\citeauthoryear{Zhang}{Zhang}{[n.d.]}]%
        {zhang2020online}
\bibfield{author}{\bibinfo{person}{Lijun Zhang}.}
  \bibinfo{year}{[n.d.]}\natexlab{}.
\newblock \showarticletitle{Online Learning in Changing Environments.}
\newblock


\bibitem[\protect\citeauthoryear{Zhang, Lu, and Zhou}{Zhang
  et~al\mbox{.}}{2018}]%
        {zhang2018adaptive}
\bibfield{author}{\bibinfo{person}{Lijun Zhang}, \bibinfo{person}{Shiyin Lu},
  {and} \bibinfo{person}{Zhi-Hua Zhou}.} \bibinfo{year}{2018}\natexlab{}.
\newblock \showarticletitle{Adaptive online learning in dynamic environments}.
  In \bibinfo{booktitle}{\emph{Advances in neural information processing
  systems}}. \bibinfo{pages}{1323--1333}.
\newblock


\bibitem[\protect\citeauthoryear{Zhang, Hajiesmaili, Cai, Chen, and Zhu}{Zhang
  et~al\mbox{.}}{2016}]%
        {zhang2016peak}
\bibfield{author}{\bibinfo{person}{Ying Zhang}, \bibinfo{person}{Mohammad~H
  Hajiesmaili}, \bibinfo{person}{Sinan Cai}, \bibinfo{person}{Minghua Chen},
  {and} \bibinfo{person}{Qi Zhu}.} \bibinfo{year}{2016}\natexlab{}.
\newblock \showarticletitle{Peak-aware online economic dispatching for
  microgrids}.
\newblock \bibinfo{journal}{\emph{IEEE transactions on smart grid}}
  \bibinfo{volume}{9}, \bibinfo{number}{1} (\bibinfo{year}{2016}),
  \bibinfo{pages}{323--335}.
\newblock


\bibitem[\protect\citeauthoryear{Zhao, Zhang, Zhang, and Zhou}{Zhao
  et~al\mbox{.}}{2020}]%
        {zhao2020dynamic}
\bibfield{author}{\bibinfo{person}{Peng Zhao}, \bibinfo{person}{Yu-Jie Zhang},
  \bibinfo{person}{Lijun Zhang}, {and} \bibinfo{person}{Zhi-Hua Zhou}.}
  \bibinfo{year}{2020}\natexlab{}.
\newblock \showarticletitle{Dynamic regret of convex and smooth functions}.
\newblock \bibinfo{journal}{\emph{Advances in Neural Information Processing
  Systems}}  \bibinfo{volume}{33} (\bibinfo{year}{2020}).
\newblock


\bibitem[\protect\citeauthoryear{Zhao, Qiu, and Liu}{Zhao
  et~al\mbox{.}}{2018}]%
        {zhao2018proximal}
\bibfield{author}{\bibinfo{person}{Yawei Zhao}, \bibinfo{person}{Shuang Qiu},
  {and} \bibinfo{person}{Ji Liu}.} \bibinfo{year}{2018}\natexlab{}.
\newblock \showarticletitle{Proximal Online Gradient is Optimum for Dynamic
  Regret}.
\newblock \bibinfo{journal}{\emph{arXiv preprint arXiv:1810.03594}}
  (\bibinfo{year}{2018}).
\newblock


\end{thebibliography}

\appendix
\clearpage
\section{Proofs for Section 4.1}
\subsection{Preliminary Lemmas}
\begin{lemma}
\label{pre_lemma 1}
For any $t\ge 1$, we have
\begin{equation}
    \sum _{ i=1 }^{ t }{ \frac { 1 }{ \sqrt { i }  }  } \le 2\sqrt { t } -1
\end{equation}
\end{lemma}
\begin{lemma}
\label{pre_lemma 2}
(Proposition A.5 in \cite{fang2020online}) Let $R>0$ and any real numbers ${x}_{1},{x}_{2},...,{x}_{T}\in [0,R]$, then we have
\begin{equation}
    \sum _{ t=1 }^{ T }{ \frac { { x }_{ t } }{ \sqrt { R+\sum _{ i<t }{ { x }_{ i } }  }  }  } \le 2\sqrt { \sum _{ t=1 }^{ T }{ { x }_{ t } }  } 
\end{equation}
\end{lemma}
%
%
% The proof of Theorem 1 is exactly the same as of Theorem 2, and hence we omit the details.
%the proof mainly follows the proof of Theorem
%
\subsection{Proof of Lemma 1}
\emph{Proof.}
\begin{enumerate}
    \item We prove the inequality (\emph{1}) by induction. Assume ${\bm{\lambda}}(\tau)\ge \bm{0}$ holds for all $\tau\in \{0,1,...,t\}$, then for $\forall k$ we consider two cases.\\
    \textbf{Case 1:} If ${g}_{k,t}({x}_{t+1})\ge 0$, then we have 
    $${ \lambda  }_{k}(t+1)=max\{ { \lambda }_{k}(t)+{ \gamma}_{t}{ g }_{ k,t }({ x }_{ t+1 }),-{ \gamma }_{t}{g }_{k, t }({ x }_{ t+1 })\}\ge { \lambda  }_{k}(t)+{ \gamma}_{t}{ g }_{k, t }({ x }_{ t+1 })\ge 0$$
    \textbf{Case 2:} If ${g}_{k,t}({x}_{t+1})< 0$, then we have
    $${ \lambda  }_{k}(t+1)=max\{ { \lambda  }_{k}(t)+{ \gamma}_{t}{ g }_{k, t }({ x }_{ t+1 }),-{ \gamma }_{t}{g }_{k, t }({ x }_{ t+1 })\}\ge -{ \gamma }_{t}{g }_{ k,t }({ x }_{ t+1 })\ge 0$$
    Thus ${\bm{\lambda}}(t)\ge \bm{0}$ holds for $\forall t$. \\
    \item Since ${ \bm{\lambda}  }(t)=max\{ { \bm{\lambda}  }(t-1)+{ \gamma}_{t-1}{ \bm{g} }_{ t-1 }({ x }_{ t }),-{ \gamma }_{t-1}{\bm{g} }_{ t-1 }({ x }_{ t })\}\ge -{ \gamma }_{t-1}{\bm{g} }_{ t-1 }({ x }_{ t })$, then we can derive that ${ \bm{\lambda}  }(t)+{ {\gamma}_{t-1} \bm{g} }_{ t-1 }({ x }_{ t })\ge 0$,\;$\forall t$.
    \item It is obvious that (\emph{3}) holds if $t=1$, then for $t\ge 2$ and $\forall k$ we consider two cases.\\
    \textbf{Case 1:} If ${g}_{k,t}({x}_{t+1})\ge 0$, then we have 
    $${ \lambda  }_{k}(t)=max\{ { \lambda  }_{k}(t-1)+{ \gamma}_{t-1}{ g}_{ k,t-1 }({ x }_{ t }),-{ \gamma }_{t-1}{g }_{ ,t-1 }({ x }_{ t })\}$$
    $$\ge { \lambda }_{k}(t-1)+{ \gamma}_{t-1}{g }_{ k,t-1 }({ x }_{ t })\ge  {\gamma}_{t-1}|{ g }_{ k,t-1 }({ x }_{ t })|$$
    \textbf{Case 2:}  If ${g}_{k,t}({x}_{t+1})< 0$, then we have 
     $${ \lambda  }_{k}(t)=max\{ { \lambda  }_{k}(t-1)+{ \gamma}_{t-1}{ g}_{ k,t-1 }({ x }_{ t }),-{ \gamma }_{t-1}{g }_{ ,t-1 }({ x }_{ t })\}\ge -{ \gamma }_{t-1}{g }_{ ,t-1 }({ x }_{ t })= {\gamma}_{t-1}|{ g }_{ k,t-1 }({ x }_{ t })|$$
     Thus we have ${ \lambda  }_{k}(t)\ge {\gamma}_{t-1}|{ g }_{ k,t-1 }({ x }_{ t })|,\;\forall t$. Squaring both sides and summing over $k$, we obtain ${ ||\bm{\lambda}  }(t){ || }^{2}\ge { {\gamma}_{t-1}^{2} ||\bm{g} }_{ t-1 }({ x }_{ t }){ || }^{2}$, which is equivalent to the inequality (\emph{3}).
     \item Since ${ \bm{\lambda}  }(t)=max\{ { \bm{\lambda}  }(t-1)+{ \gamma}_{t-1}{ \bm{g} }_{ t-1 }({ x }_{ t }),-{ \gamma }_{t-1}{\bm{g} }_{ t-1 }({ x }_{ t })\}\ge { \bm{\lambda}  }(t-1)+{ \gamma}_{t-1}{ \bm{g} }_{ t-1 }({ x }_{ t })$, then we have ${ {\gamma}_{t-1} \bm{g} }_{ t-1}({ x }_{ t })\le { \bm{\lambda}  }(t)-{ \bm{\lambda}  }(t-1)$. Furthermore, 
     $${ \lambda  }_{ k }(t)=max\{ { \lambda  }_{ k }(t-1)+{ \gamma  }_{ t-1 }{ g }_{ k,t-1 }({ x }_{ t }),-{ \gamma  }_{ t-1 }{ g }_{ ,t-1 }({ x }_{ t })\}$$
     $$\le |{ \lambda  }_{ k }(t-1)|+|{ \gamma  }_{ t-1 }{ g }_{ k,t-1 }({ x }_{ t })|={ \lambda  }_{ k }(t-1)+{ \gamma  }_{ t-1 }|{ g }_{ k,t-1 }({ x }_{ t })|$$
     Squaring both sides and summing over $k$, we obtain
    $$||{ \bm{\lambda}  }(t){||}^{2}\le ||{ \bm{\lambda}  }(t-1)+{ \gamma  }_{ t-1 }{ \bm{g} }_{ t-1 }({ x }_{ t }){||}^{2}\Leftrightarrow  ||{ \bm{\lambda}  }(t)||\le ||{ \bm{\lambda}  }(t-1)+{ \gamma  }_{ t-1 }{ \bm{g} }_{ t-1 }({ x }_{ t })||$$
     By the triangle inequality we have  ${ ||\bm{\lambda}  }(t){ || }{ \le {\gamma}_{t-1} ||\bm{g} }_{ t-1 }({ x }_{ t }){ || }+{ ||\bm{\lambda}  }(t-1){ || }$.
     \item According to the above inequality ${ ||\bm{\lambda}  }(t){ || }{ \le {\gamma}_{t-1} ||\bm{g} }_{ t-1 }({ x }_{ t }){ || }+{ ||\bm{\lambda}  }(t-1){ || }$ , we have
     $${ ||\bm{\lambda}  }(t+1){ || }{ \le {\gamma}_{t} ||\bm{g} }_{ t }({ x }_{ t+1 }){ || }+{ ||\bm{\lambda}  }(t){ || }\quad \quad \quad \quad \quad \quad \quad \quad \quad \quad $$
     $$\Rightarrow ||{ \bm{\lambda}  }(t+1){||}^{2}\le ||{ \bm{\lambda}  }(t){||}^{2}+2{ \gamma  }_{ t }[{ \bm{\lambda}  }(t){ ] }^{ T }{ \bm{g} }_{ t }({ x }_{ t+1 })+{ \gamma  }_{ t }^{ 2 }||{ \bm{g} }_{ t }({ x }_{ t+1 }){ || }^{ 2 }$$
     Rearranging terms yields the inequality (\emph{5}).
\end{enumerate} 
\subsection{Proof of Lemma 2}
\emph{Proof.} Since ${ \nabla { f }_{ t }({ x }_{ t }) }^{ T }(x-{ x }_{ t })+[{ \bm{\lambda} }(t)+{\gamma}_{t-1} { { \bm{g} } }_{ t-1 }({ x }_{ t }){ ] }^{ T }({\gamma}_{t} { { \bm{g} } }_{ t }(x))+{ \alpha  }_{t}||x-{ x }_{ t }{ || }^{ 2 }$ is a $2{\alpha}_{t}$-strong convex function with respect to \emph{x} and ${x}_{t+1}$ minimizes this expression over $\chi$, we have
\begin{equation}
    \begin{split}
        \label{convexity analysis 1}
        &{ \nabla { f }_{ t }({ x }_{ t }) }^{ T }({ x }_{ t+1 }-{ x }_{ t })+[{ \bm{\lambda}  }(t)+{\gamma}_{t-1} { { \bm{ g } } }_{ t-1 }({ x }_{ t }){ ] }^{ T }({\gamma}_{t} { { \bm{ g } } }_{ t }({ x }_{ t+1 }))+{ { \alpha  }_{ t } }||{ x }_{ t+1 }-{ x }_{ t }{ || }^{ 2 }\\ 
        &\le { \nabla { f }_{ t }({ x }_{ t }) }^{ T }({ x }_{ t }^{ * }-{ x }_{ t })+[{ \bm{\lambda}  }(t)+{\gamma}_{t-1} { { \bm{ g } } }_{ t-1 }({ x }_{ t }){ ] }^{ T }({\gamma}_{t} { { \bm{ g } } }_{ t }({ x }_{ t }^{ * }))+{ { \alpha  }_{ t } }||{ x }_{ t }^{ * }-{ x }_{ t }{ || }^{ 2 }-{ { \alpha  }_{ t } }||{ x }_{ t+1 }-{ x }_{ t }^{ * }{ || }^{ 2 }\\ 
        &\overset { (a) }{ \le  } { \nabla { f }_{ t }({ x }_{ t }) }^{ T }({ x }_{ t }^{ * }-{ x }_{ t })+{ { \alpha  }_{ t } }||{ x }_{ t }^{ * }-{ x }_{ t }{ || }^{ 2 }-{ { \alpha  }_{ t } }||{ x }_{ t+1 }-{ x }_{ t }^{ * }{ || }^{ 2 }
    \end{split}
\end{equation}
Where (a) follows from the fact that ${ \bm{ g } }_{ t }({ x }_{ t }^{ * })\le 0$ and Lemma \ref{pre_lemma 2}. Adding ${f}_{t}({x}_{t})$ on both sides of (\ref{convexity analysis 1}) and using the convexity of ${f}_{t}$, we have
\begin{equation}
    \begin{split}
       \label{1_proof_of_lemma_1}
        &{ f }_{ t }({ x }_{ t })+{ \nabla { f }_{ t }({ x }_{ t }) }^{ T }({ x }_{ t+1 }-{ x }_{ t })+[{ \bm{\lambda } }(t)+{\gamma}_{t-1} { { \bm{ g } } }_{ t-1 }({ x }_{ t }){ ] }^{ T }({\gamma}_{t} { { \bm{ g } } }_{ t }({ x }_{ t+1 }))+{ { \alpha  }_{ t } }||{ x }_{ t+1 }-{ x }_{ t }{ || }^{ 2 }\\ 
        &\le { { f }_{ t }({ x }_{ t })+\nabla { f }_{ t }({ x }_{ t }) }^{ T }({ x }_{ t }^{ * }-{ x }_{ t })+{ { \alpha  }_{ t } }||{ x }_{ t }^{ * }-{ x }_{ t }{ || }^{ 2 }-{ { \alpha  }_{ t } }||{ x }_{ t+1 }-{ x }_{ t }^{ * }{ || }^{ 2 }\\ 
        &\le { f }_{ t }({ x }_{ t }^{ * })+{ { \alpha  }_{ t } }||{ x }_{ t }^{ * }-{ x }_{ t }{ || }^{ 2 }-{ { \alpha  }_{ t } }||{ x }_{ t+1 }-{ x }_{ t }^{ * }{ || }^{ 2 }
    \end{split}
\end{equation}
Rearranging terms in (\ref{1_proof_of_lemma_1}), we have
\begin{equation}
    \begin{split}
        \label{2_proof_of_lemma_1}
        &{ f }_{ t }({ x }_{ t })+[{ \bm{\lambda } }(t)]^{ T }({ \gamma  }_{ t }{\bm { g } }_{ t }({ x }_{ t+1 }))-{ f }_{ t }({ x }_{ t }^{ * })\\ 
        &\le { { \alpha  }_{ t } }||{ x }_{ t }^{ * }-{ x }_{ t }{ || }^{ 2 }-{ { \alpha  }_{ t } }||{ x }_{ t+1 }-{ x }_{ t }^{ * }{ || }^{ 2 }-{ { \alpha  }_{ t } }||{ x }_{ t+1 }-{ x }_{ t }{ || }^{ 2 }-{ \gamma  }_{ t-1 }{ \gamma  }_{ t }[{ \bm{ g } }_{ t-1 }({ x }_{ t }){ ] }^{ T }{ \bm{ g } }_{ t }({ x }_{ t+1 })-{ \nabla { f }_{ t }({ x }_{ t }) }^{ T }({ x }_{ t+1 }-{ x }_{ t })\\ 
        &\overset { (a) }{ \le  } { { \alpha  }_{ t } }||{ x }_{ t }^{ * }-{ x }_{ t }{ || }^{ 2 }-{ { \alpha  }_{ t } }||{ x }_{ t+1 }-{ x }_{ t }^{ * }{ || }^{ 2 }-{ { \alpha  }_{ t } }||{ x }_{ t+1 }-{ x }_{ t }{ || }^{ 2 }-{ \gamma  }_{ t-1 }{ \gamma  }_{ t }[{\bm { g } }_{ t-1 }({ x }_{ t }){ ] }^{ T }{ \bm{ g } }_{ t }({ x }_{ t+1 })+||{ \nabla { f }_{ t }({ x }_{ t }) }{ || }||{ x }_{ t+1 }-{ x }_{ t }{ || }\\ 
        &\overset { (b) }{ \le  } { { \alpha  }_{ t } }||{ x }^{ * }-{ x }_{ t }{ || }^{ 2 }-{ { \alpha  }_{ t } }||{ x }_{ t+1 }-{ x }_{ t }^{ * }{ || }^{ 2 }-{ { \alpha  }_{ t } }||{ x }_{ t+1 }-{ x }_{ t }{ || }^{ 2 }-{ \gamma  }_{ t-1 }{ \gamma  }_{ t }[{\bm { g } }_{ t-1 }({ x }_{ t }){ ] }^{ T }{\bm { g } }_{ t }({ x }_{ t+1 })
        +\frac { 1 }{ 2\delta  } ||\nabla { f }_{ t }({ x }_{ t }){ || }^{ 2 }+\frac { \delta  }{ 2 } ||{ x }_{ t+1 }-{ x }_{ t }{ || }^{ 2 }\\ &\overset { (c) }{ \le  } { { \alpha  }_{ t } }||{ x }_{ t }^{ * }-{ x }_{ t }{ || }^{ 2 }-{ { \alpha  }_{ t } }||{ x }_{ t+1 }-{ x }_{ t }^{ * }{ || }^{ 2 }-{ \alpha  }_{ t }||{ x }_{ t+1 }-{ x }_{ t }{ || }^{ 2 }-{ \gamma  }_{ t-1 }{ \gamma  }_{ t }[{ \bm{ g } }_{ t-1 }({ x }_{ t }){ ] }^{ T }{ \bm{ g } }_{ t }({ x }_{ t+1 })
        +\frac { 1 }{ 2\delta  } { G }^{ 2 }+\frac { \delta  }{ 2 } ||{ x }_{ t+1 }-{ x }_{ t }{ || }^{ 2 }
    \end{split}
\end{equation}
Where (a) holds by the Cauchy-Schwarz inequality; (b) comes from the AM–GM inequality; (c) holds due to the Assumption 1.
Based on Assumption 1, we note that
\begin{equation}
    \begin{split}
       \label{3_proof_of_lemma_1}
        &||{ x }_{ t }^{ * }-{ x }_{ t }{ || }^{ 2 }-||{ x }_{ t+1 }{ -{ x }_{ t }^{ * }|| }^{ 2 }
        =||{ x }_{ t }^{ * }-{ x }_{ t }{ || }^{ 2 }-||{ x }_{ t+1 }{ -{ x }_{ t+1 }^{ * }+{ x }_{ t+1 }^{ * }-{ x }_{ t }^{ * }|| }^{ 2 }\\ 
        &\le ||{ x }_{ t }^{ * }-{ x }_{ t }{ || }^{ 2 }-||{ x }_{ t+1 }{ -{ x }_{ t+1 }^{ * }|| }^{ 2 }-||{ { x }_{ t+1 }^{ * }-{ x }_{ t }^{ * }|| }^{ 2 }-2||{ x }_{ t+1 }{ -{ x }_{ t+1 } ^{ * }{ || }||{ x }_{ t+1 }^{ * }-{ x }_{ t }^{ * }|| }\\ &\le ||{ x }_{ t }^{ * }-{ x }_{ t }{ || }^{ 2 }-||{ x }_{ t+1 }{ -{ x }_{ t+1 }^{ * }|| }^{ 2 }+4R{ ||{ x }_{ t+1 }^{ * }-{ x }_{ t }^{ * }|| }
    \end{split}
\end{equation}
And
\begin{equation}
    \begin{split}
        \label{4_proof_of_lemma_1}
        &-[{ { \bm{g} } }_{ t-1 }({ x }_{ t }){ ] }^{ T }{ { \bm{g} } }_{ t }({ x }_{ t+1 })=-\frac { 1 }{ 2 } ||{ { \bm{g} } }_{ t-1 }({ x }_{ t }){ || }^{ 2 }-\frac { 1 }{ 2 } ||{ { \bm{g} } }_{ t }({ x }_{ t+1 }){ || }^{ 2 }+\frac { 1 }{ 2 } ||{ { \bm{g} } }_{ t-1 }({ x }_{ t })-{ { \bm{g} } }_{ t }({ x }_{ t+1 }){ || }^{ 2 }\\ 
        &=-\frac { 1 }{ 2 } ||{ { \bm{g} } }_{ t-1 }({ x }_{ t }){ || }^{ 2 }-\frac { 1 }{ 2 } ||{ { \bm{g} } }_{ t }({ x }_{ t+1 }){ || }^{ 2 }+\frac { 1 }{ 2 } ||{ { \bm{g} } }_{ t-1 }({ x }_{ t })-{ { \bm{g} } }_{ t }({ x }_{ t })+{ { \bm{g} } }_{ t }({ x }_{ t })-{ { \bm{g} } }_{ t }({ x }_{ t+1 }){ || }^{ 2 }\\
        &\overset { (a) }{ \le  } -\frac { 1 }{ 2 } ||{ { \bm{g} } }_{ t-1 }({ x }_{ t }){ || }^{ 2 }-\frac { 1 }{ 2 } ||{ { \bm{g} } }_{ t }({ x }_{ t+1 }){ || }^{ 2 }+\frac { 1 }{ 2 } [2||{ { \bm{g} } }_{ t-1 }({ x }_{ t })-{ { \bm{g} } }_{ t }({ x }_{ t }){ || }^{ 2 }+2||{ { \bm{g} } }_{ t }({ x }_{ t })-{ { \bm{g} } }_{ t }({ x }_{ t+1 }){ || }^{ 2 }]\\ 
        &\overset { (b) }{ \le  } -\frac { 1 }{ 2 } ||{ { \bm{g} } }_{ t-1 }({ x }_{ t }){ || }^{ 2 }-\frac { 1 }{ 2 } ||{ { \bm{g} } }_{ t }({ x }_{ t+1 }){ || }^{ 2 }+||{ { \bm{g} } }_{ t-1 }({ x }_{ t })-{ { \bm{g} } }_{ t }({ x }_{ t }){ || }^{ 2 }+{ \beta  }^{ 2 }||{ x }_{ t+1 }-{ x }_{ t }{ || }^{ 2 }
    \end{split}
\end{equation}
Where (a) follows from the AM-GM inequality; (b) holds by the Lipschitz continuity of ${\bm{g}}_{t}$ (Assumption 1). Substituting (\ref{3_proof_of_lemma_1}) and (\ref{4_proof_of_lemma_1}) into (\ref{2_proof_of_lemma_1}) we obtain
\begin{equation}
    \begin{split}
     \label{5_proof_of_lemma_1}
        &{ f }_{ t }({ x }_{ t })+[{ \bm{\lambda}  }(t)]^{ T }(\gamma { {\bm { g } } }_{ t }({ x }_{ t+1 }))-{ f }_{ t }({ x }_{ t }^{ * })\\ 
        &\le { \alpha  }_{ t }[||{ x }_{ t }^{ * }-{ x }_{ t }{ || }^{ 2 }-||{ x }_{ t+1 }{ -{ x }_{ t+1 }^{ * }|| }^{ 2 }+4R{ ||{ x }_{ t+1 }^{ * }-{ x }_{ t }^{ * }|| }+({ \beta  }^{ 2 }{ \gamma  }_{ t }{ \gamma  }_{ t-1 }+\frac { \delta  }{ 2 } -{ \alpha  }_{ t })||{ x }_{ t+1 }-{ x }_{ t }{ || }^{ 2 }\\
        &+\frac { 1 }{ 2\delta  } { G }^{ 2 }-\frac { 1 }{ 2 } { \gamma  }_{ t-1 }{ \gamma  }_{ t }||{ {\bm { g } } }_{ t-1 }({ x }_{ t }){ || }^{ 2 }-\frac { 1 }{ 2 } { \gamma  }_{ t-1 }{ \gamma  }_{ t }||{ { \bm{ g } } }_{ t }({ x }_{ t+1 }){ || }^{ 2 }+{ \gamma  }_{ t-1 }{ \gamma  }_{ t }||{ {\bm { g } } }_{ t-1 }({ x }_{ t })-{ {\bm { g } } }_{ t }({ x }_{ t }){ || }^{ 2 }\\ 
        &\overset { (a) }{ \le  } { \alpha  }_{ t }||{ x }_{ t }^{ * }-{ x }_{ t }{ || }^{ 2 }-{ \alpha  }_{ t+1 }||{ x }_{ t+1 }{ -{ x }_{ t+1 }^{ * }|| }^{ 2 }+4R{ \alpha  }_{ t }{ ||{ x }_{ t+1 }^{ * }-{ x }_{ t }^{ * }|| }+({ \beta  }^{ 2 }{ \gamma  }_{ t-1 }+\frac { \delta  }{ 2 } -{ \alpha  }_{ t })||{ x }_{ t+1 }-{ x }_{ t }{ || }^{ 2 }\\ 
        &+\frac { 1 }{ 2\delta  } { G }^{ 2 }-\frac { 1 }{ 2 } { \gamma  }_{ t-1 }{ \gamma  }_{ t }||{ { \bm{ g } } }_{ t-1 }({ x }_{ t }){ || }^{ 2 }-\frac { 1 }{ 2 } { \gamma  }_{ t-1 }{ \gamma  }_{ t }||{ { \bm{ g } } }_{ t }({ x }_{ t+1 }){ || }^{ 2 }+{ \gamma  }_{ t-1 }^{ 2 }||{ { \bm{ g } } }_{ t-1 }({ x }_{ t })-{ { \bm{ g } } }_{ t }({ x }_{ t }){ || }^{ 2 }
    \end{split}
\end{equation}
Where (a) comes from the fact that both $\{{\alpha}_{t}\}$ and $\{{\gamma}_{t}\}$ are non-increasing sequence. According to Lemma \ref{lemma 0} and adding Lyapunov drift term on both sides of (\ref{5_proof_of_lemma_1}) yields:
\begin{equation}
    \begin{split}
       &{ f }_{ t }({ x }_{ t })+\Delta (t)-{ f }_{ t }({ x }_{ t }^{ * })\\ 
       &\le { \alpha  }_{ t }||{ x }_{ t }^{ * }-{ x }_{ t }{ || }^{ 2 }-{ \alpha  }_{ t+1 }||{ x }_{ t+1 }{ -{ x }_{ t+1 }^{ * }|| }^{ 2 }+4R{ \alpha  }_{ t }{ ||{ x }_{ t+1 }^{ * }-{ x }_{ t }^{ * }|| }+({ \beta  }^{ 2 }{ \gamma  }_{ t-1 }+\frac { \delta  }{ 2 } -{ \alpha  }_{ t })||{ x }_{ t+1 }-{ x }_{ t }{ || }^{ 2 }\\ 
       &+\frac { 1 }{ 2\delta  } { G }^{ 2 }-\frac { 1 }{ 2 } { \gamma  }_{ t-1 }{ \gamma  }_{ t }||{ { \bm{ g } } }_{ t-1 }({ x }_{ t }){ || }^{ 2 }+({ \gamma  }_{ t }^{ 2 }-\frac { 1 }{ 2 } { \gamma  }_{ t-1 }{ \gamma  }_{ t })||{ {\bm { g } } }_{ t }({ x }_{ t+1 }){ || }^{ 2 }+{ \gamma  }_{ t-1 }^{ 2 }||{ {\bm { g } } }_{ t-1 }({ x }_{ t })-{ { \bm{ g } } }_{ t }({ x }_{ t }){ || }^{ 2 }\\
       &\overset { (a) }{ \le  } { \alpha  }_{ t }||{ x }_{ t }^{ * }-{ x }_{ t }{ || }^{ 2 }-{ \alpha  }_{ t+1 }||{ x }_{ t+1 }{ -{ x }_{ t+1 }^{ * }|| }^{ 2 }+4R{ \alpha  }_{ t }{ ||{ x }_{ t+1 }^{ * }-{ x }_{ t }^{ * }|| }+({ \beta  }^{ 2 }{ \gamma  }_{ t-1 }+\frac { \delta  }{ 2 } -{ \alpha  }_{ t })||{ x }_{ t+1 }-{ x }_{ t }{ || }^{ 2 }\\ 
       &+\frac { 1 }{ 2\delta  } { G }^{ 2 }-\frac { 1 }{ 2 } { \gamma  }_{ t-1 }{ \gamma  }_{ t }||{ { { g } } }_{ t-1 }({ x }_{ t }){ || }^{ 2 }+\frac { 1 }{ 2 } { \gamma  }_{ t }{ \gamma  }_{ t+1 }||{ {\bm { g } } }_{ t }({ x }_{ t+1 }){ || }^{ 2 }+{ \gamma  }_{ t-1 }^{ 2 }||{ {\bm { g } } }_{ t-1 }({ x }_{ t })-{ { \bm{ g } } }_{ t }({ x }_{ t }){ || }^{ 2 }
    \end{split}
\end{equation}
Where (a) is due to the fact that $2{ \gamma  }_{ t }\le { \gamma  }_{ t-1 }+{ \gamma  }_{ t+1 }$. This completes the proof.
\subsection{Proof of Lemma 3}
\emph{Proof.} According to Lemma 1, taking a telescoping sum over $t= 1, ..., T$, we obtain
\begin{equation}
    \begin{split}
    \label{1_proof_of_lemma_2}
        &\sum _{ t=1 }^{ T }{ { f }_{ t }({ x }_{ t }) } +\sum _{ t=1 }^{ T }{ \Delta (t) } \le \sum _{ t=1 }^{ T }{ { f }_{ t }({ x }_{ t }^{ * }) } +4R\sum _{ t=1 }^{ T }{ { \alpha  }_{ t }{ ||{ x }_{ t+1 }^{ * }-{ x }_{ t }^{ * }|| } } +{ \alpha  }_{ 1 }||{ x }_{ 1 }^{ * }-{ x }_{ 1 }{ || }^{ 2 }+\frac { T{ G }^{ 2 } }{ 2\delta  } \\
        &+\frac { 1 }{ 2 } { \gamma  }_{ T }{ \gamma  }_{ T+1 }||{ { \bm{ g } } }_{ T }({ x }_{ T+1 }){ || }^{ 2 }+\sum _{ t=1 }^{ T }{ { \gamma  }_{ t-1 }^{ 2 }||{ { \bm{ g } } }_{ t-1 }({ x }_{ t })-{ { \bm{ g } } }_{ t }({ x }_{ t }){ || }^{ 2 } } +\sum _{ t=1 }^{ T }{ ({ \beta  }^{ 2 }{ \gamma  }_{ t-1 }^{ 2 }+\frac { \delta  }{ 2 } -{ \alpha  }_{ t })||{ x }_{ t+1 }-{ x }_{ t }{ || }^{ 2 } } \\ 
        &\overset { (a) }{ \le  } \sum _{ t=1 }^{ T }{ { f }_{ t }({ x }_{ t }^{ * }) } +4R\sum _{ t=1 }^{ T }{ { \alpha  }_{ t }{ ||{ x }_{ t+1 }^{ * }-{ x }_{ t }^{ * }|| } } +{ \alpha  }_{ 1 }{ R }^{ 2 }+\frac { T{ G }^{ 2 } }{ 2\delta  } \\ 
        &+\frac { 1 }{ 2 } { \gamma  }_{ T }{ \gamma  }_{ T+1 }||{ { \bm{ g } } }_{ T }({ x }_{ T+1 }){ || }^{ 2 }+2F\sum _{ t=1 }^{ T }{ { \gamma  }_{ t-1 }^{ 2 }||{ { \bm{ g } } }_{ t-1 }({ x }_{ t })-{ { \bm{ g } } }_{ t }({ x }_{ t }){ || } } 
    \end{split}
\end{equation}
Where (a) holds by the Assumption 1 and the fact that ${\alpha}_{t}\ge { \beta  }^{ 2 }{ \gamma  }^{ 2 }+\frac { \delta  }{ 2 }$. Rearranging terms yields:
\begin{equation}
    \begin{split}
        &\sum _{ t=1 }^{ T }{ { f }_{ t }({ x }_{ t }) } \le \sum _{ t=1 }^{ T }{ { f }_{ t }({ x }_{ t }^{ * }) } +4R\sum _{ t=1 }^{ T }{ { \alpha  }_{ t }{ ||{ x }_{ t+1 }^{ * }-{ x }_{ t }^{ * }|| } } +{ \alpha  }_{ 1 }{ R }^{ 2 }+\frac { T{ G }^{ 2 } }{ 2\delta  } \\ 
        &+\frac { 1 }{ 2 } { \gamma  }_{ T }{ \gamma  }_{ T+1 }||{ { { g } } }_{ T }({ x }_{ T+1 }){ || }^{ 2 }+2F\sum _{ t=1 }^{ T }{ { \gamma  }_{ t-1 }^{ 2 }||{ { { g } } }_{ t-1 }({ x }_{ t })-{ { { g } } }_{ t }({ x }_{ t }){ || } } +\frac { 1 }{ 2 } ||\lambda (1){ || }^{ 2 }
    \end{split}
\end{equation}
It completes the proof.
\subsection{Proof of Lemma 4}
\emph{Proof.} According to Lemma 1, taking a telescoping sum over $t= 1, ..., T-1$ and using the fact that ${\alpha}_{t}\ge { \beta  }^{ 2 }{ \gamma  }_{t-1}^{ 2 }+\frac { \delta  }{ 2 }$, we obtain
\begin{equation}
    \begin{split}
        &\sum _{ t=1 }^{ T-1 }{ { f }_{ t }({ x }_{ t }) } +\sum _{ t=1 }^{ T-1 }{ \Delta (t) } \le \sum _{ t=1 }^{ T-1 }{ { f }_{ t }({ x }_{ t }^{ * }) } +4R\sum _{ t=1 }^{ T-1 }{ { \alpha  }_{ t }{ ||{ x }_{ t+1 }^{ * }-{ x }_{ t }^{ * }|| } } +{ \alpha  }_{ 1 }{ R }^{ 2 }\\ 
        &+\frac { (T-1){ G }^{ 2 } }{ 2\delta  } +\frac { 1 }{ 2 } { \gamma  }_{ T-1 }{ \gamma  }_{ T }||{ {\bm { g } } }_{ T-1 }({ x }_{ T }){ || }^{ 2 }+2F\sum _{ t=1 }^{ T-1 }{ { \gamma  }_{ t-1 }^{ 2 }||{ { \bm{ g } } }_{ t-1 }({ x }_{ t })-{ { \bm{ g } } }_{ t }({ x }_{ t }){ || } } 
    \end{split}
\end{equation}
Rearranging terms and multiplying both sides by 2 yields:
\begin{equation}
    \begin{split}
        &||\bm{\lambda} (t){ || }^{ 2 }\le 2(\sum _{ t=1 }^{ T-1 }{ { f }_{ t }({ x }_{ t }^{ * }) } -\sum _{ t=1 }^{ T-1 }{ { f }_{ t }({ x }_{ t }) } )+8R\sum _{ t=1 }^{ T-1 }{ { \alpha  }_{ t }{ ||{ x }_{ t+1 }^{ * }-{ x }_{ t }^{ * }|| } } +2{ \alpha  }_{ 1 }{ R }^{ 2 }\\ 
        &+\frac { (T-1){ G }^{ 2 } }{ \delta  } +{ \gamma  }_{ T-1 }{ \gamma  }_{ T }||{ { \bm{ g } } }_{ T-1 }({ x }_{ T }){ || }^{ 2 }+4F\sum _{ t=1 }^{ T-1 }{ { \gamma  }_{ t-1 }^{ 2 }||{ { \bm{ g } } }_{ t-1 }({ x }_{ t })-{ { \bm{ g } } }_{ t }({ x }_{ t }){ || } } \\ 
        &\overset { (a) }{ \le  } 2F(T-1)+8R\sum _{ t=1 }^{ T-1 }{ { \alpha  }_{ t }{ ||{ x }_{ t+1 }^{ * }-{ x }_{ t }^{ * }|| } } +2{ \alpha  }_{ 1 }{ R }^{ 2 }+\frac { (T-1){ G }^{ 2 } }{ \delta  } \\ 
        &+{ \gamma  }_{ T-1 }^{ 2 }||{ { \bm{ g } } }_{ T-1 }({ x }_{ T }){ || }^{ 2 }+4F\sum _{ t=1 }^{ T-1 }{ { \gamma  }_{ t-1 }^{ 2 }||{ { \bm{ g } } }_{ t-1 }({ x }_{ t })-{ { \bm{ g } } }_{ t }({ x }_{ t }){ || } } 
    \end{split}
\end{equation}
Where (a) holds by the Assumption 1 and the fact that ${\gamma}_{T}\le{\gamma}_{T-1}$. Taking the square root of both sides and using the fact that $\sqrt { \sum _{ i }{ { a }_{ i } }  } \le \sum _{ i }{ \sqrt { { a }_{ i } }  } ,\; \forall { a }_{ i }\ge 0$, we obtain
\begin{equation}
    \begin{split}
        &||\bm{\lambda} (t){ || }\le \sqrt { 2F(T-1) } +2\sqrt { 2R\sum _{ t=1 }^{ T-1 }{ { \alpha  }_{ t }{ ||{ x }_{ t+1 }^{ * }-{ x }_{ t }^{ * }|| } }  } +\sqrt { 2{ \alpha  }_{ 1 }{ R }^{ 2 } } \\ 
        &+\sqrt { \frac { (T-1){ G }^{ 2 } }{ \delta  }  } +{ \gamma  }_{ T-1 }||{ {\bm { g } } }_{ T-1 }({ x }_{ T }){ || }+2\sqrt { F\sum _{ t=1 }^{ T-1 }{ { \gamma  }_{ t-1 }^{ 2 }||{ { \bm{ g } } }_{ t-1 }({ x }_{ t })-{ { \bm{ g } } }_{ t }({ x }_{ t }){ || } }  } 
    \end{split}
\end{equation}
It completes the proof.
\subsection{Proof of Lemma 5}
\emph{Proof.} According to Lemma \ref{lemma 0}, we have ${ { \gamma  }_{ t-1 }g }_{ t-1,k }({ x }_{ t })\le { \lambda  }_{ k }(t)-{ \lambda  }_{ k }(t-1)\Leftrightarrow { g }_{ t-1,k }({ x }_{ t })\le \frac { { \lambda  }_{ k }(t) }{ { \gamma  }_{ t-1 } } -\frac { { \lambda  }_{ k }(t-1) }{ { \gamma  }_{ t-1 } } $. Adding ${g }_{ t,k }({ x }_{ t })$ on both sides of it and telescoping it over $t$ yields:
\begin{equation}
    \begin{split}
        &\sum _{ t=1 }^{ T }{ { g }_{ t-1,k }({ x }_{ t })+{ g }_{ t,k }({ x }_{ t }) } \le \sum _{ t=1 }^{ T }{ { g }_{ t,k }({ x }_{ t }) } +\sum _{ t=1 }^{ T }{ \frac { { \lambda  }_{ k }(t) }{ { \gamma  }_{ t-1 } } -\frac { { \lambda  }_{ k }(t-1) }{ { \gamma  }_{ t-1 } }  } \\ 
        &\overset { (a) }{ \le  } \sum _{ t=1 }^{ T }{ { g }_{ t,k }({ x }_{ t }) } +\sum _{ t=1 }^{ T }{ \frac { { \lambda  }_{ k }(t) }{ { \gamma  }_{ t } } -\frac { { \lambda  }_{ k }(t-1) }{ { \gamma  }_{ t-1 } }  }
        \le \sum _{ t=1 }^{ T }{ { g }_{ t,k }({ x }_{ t }) } +\frac { { \lambda  }_{ k }(T) }{ { \gamma  }_{ T } } \\ 
        &\Rightarrow \sum _{ t=1 }^{ T }{ { g }_{ t,k }({ x }_{ t }) } \le \sum _{ t=1 }^{ T }{ { g }_{ t,k }({ x }_{ t })-{ g }_{ t-1,k }({ x }_{ t }) } +\frac { { \lambda  }_{ k }(T) }{ { \gamma  }_{ T } } \\ 
        &\overset { (b) }{ \le  } \sum _{ t=1 }^{ T }{ { max\{ |g }_{ t,k }({ x })-{ g }_{ t-1,k }({ x })|\}  } +\frac { { \lambda  }_{ k }(T) }{ { \gamma  }_{ T } } 
        \overset { (c) }{ \le  } \frac { ||{ \bm{\lambda}  }(T)|| }{ { \gamma  }_{ T } } +{ V }_{ \bm{g} },\;\forall k=1,2,...,K
    \end{split}
\end{equation}
Where (a) comes from the fact that ${\gamma}_{t}$ is non-increasing with respect to $t$; (b) follows from the fact that
$ { |g }_{ t,k }({ x }_{ t })-{ g }_{ t-1,k }({ x }_{ t })|\le { max\{ |g }_{ t,k }({ x })-{ g }_{ t-1,k }({ x })|\} $; (c) is due to the fact that ${ |g }_{ t,k }({ x })-{ g }_{ t-1,k }({ x })|\le ||{ \bm{g} }_{ t }(x)-{ \bm{g} }_{ t-1 }(x){ || }_{ 2 }$ and ${\lambda}_{k}(T)\le||{\bm{\lambda}}(T)||$.
It completes the proof.
\subsection{Proof of Lemma 6}
\emph{Proof.} Combining Lemma 4 with Lemma 5, we have
\begin{equation}
    \begin{split}
        &\sum _{ t=1 }^{ T }{ { g }_{ t,k }({ x }_{ t }) } \le \frac { ||{ { \bm{\lambda}  } }(T)|| }{ { \gamma  }_{ T } } +{ V }_{\bm { g } } 
        \le \frac { 2 }{ { \gamma  }_{ T } } \sqrt { F(T-1) } +\frac { 1 }{ { \gamma  }_{ T } } \sqrt { 2{ \alpha  }_{ 1 }{ R }^{ 2 } } +\frac { { \gamma  }_{ T-1 } }{ { \gamma  }_{ T } } ||{ { {\bm { g } } } }_{ T-1 }({ x }_{ T }){ || }+\frac { G }{ { \gamma  }_{ T } } \sqrt { \frac { T-1 }{ \delta  }  } \\ 
        &+\frac { 2 }{ { \gamma  }_{ T } } \sqrt { 2R\sum _{ t=1 }^{ T-1 }{ { \alpha  }_{ t }||{ x }_{ t+1 }^{ * }-{ x }_{ t }^{ * }{ || } }  } +\frac { 2 }{ { \gamma  }_{ T } } \sqrt { F\sum _{ t=1 }^{ T-1 }{ { \gamma  }_{ t-1 }^{ 2 }||{ { { \bm{ g } } } }_{ t-1 }({ x }_{ t })-{ { { \bm{ g } } } }_{ t }({ x }_{ t }){ || } }  } +{ V }_{ \bm{ g } }
    \end{split}
\end{equation}
It completes the proof.
\section{Proofs for Section 4.2}
Note that both $\{{\alpha}_{t}\}$ and $\{{\gamma}_{t}\}$ are constant sequences in Algorithm 2, thus here we omit the subscript $t$.
We give the complete proofs of all lemmas in section 4.2 in the following.
\subsection{Proof of Lemma 8}
\emph{Proof.}
We conduct a similar derivation process as the proof of Lemma 1, note that ${ \nabla { f }_{ t }({ x }_{ t }) }^{ T }(x-{ x }_{ t })+[{ \bm{\lambda} }(t)+\gamma { { \bm{g} } }_{ t }({ x }_{ t }){ ] }^{ T }(\gamma { { \bm{g} } }_{ t }(x))+{ \alpha  }||x-{ x }_{ t }{ || }^{ 2 }$ is a $2{\alpha}$-strong convex function with respect to \emph{x} and ${x}_{t+1}$ minimizes this expression over $\chi$, we have
\begin{equation}
    \begin{split}
        \label{1_proof_of_lemma_B_2_2}
        & { \nabla { f }_{ t }({ x }_{ t }) }^{ T }({ x }_{ t+1 }-{ x }_{ t })+[{ \lambda  }(t)+\gamma { { \bm{g} } }_{ t}({ x }_{ t }){ ] }^{ T }(\gamma { { \bm{g} } }_{ t }({ x }_{ t+1 }))+{ \alpha  }||{ x }_{ t+1 }-{ x }_{ t }{ || }^{ 2 }\\ 
         &\le { \nabla { f }_{ t }({ x }_{ t }) }^{ T }({ x }_{ t }^{ * }-{ x }_{ t })+[{ \lambda  }(t)+\gamma { { \bm{g} } }_{ t }({ x }_{ t }){ ] }^{ T }(\gamma { { \bm{g} } }_{ t }({ x }_{ t }^{ * }))+{ \alpha  }||{ x }_{ t }^{ * }-{ x }_{ t }{ || }^{ 2 }-{ \alpha  }||{ x }_{ t+1 }-{ x }_{ t }^{ * }{ || }^{ 2 }\\ 
         &\overset { (a) }{ \le  } { \nabla { f }_{ t }({ x }_{ t }) }^{ T }({ x }_{ t }^{ * }-{ x }_{ t })+{ \alpha  }||{ x }_{ t }^{ * }-{ x }_{ t }{ || }^{ 2 }-{ \alpha  }||{ x }_{ t+1 }-{ x }_{ t }^{ * }{ || }^{ 2 }
    \end{split}
\end{equation}
Where (a) follows from the fact that ${ { g } }_{ t }({ x }_{ t }^{ * })\le 0$ and Lemma \ref{lemma 7}. Next we add ${f}_{t}({x}_{t})$ on both sides of (\ref{1_proof_of_lemma_B_2_2}) and use the convexity of ${f}_{t}$, then we obtain
\begin{equation}
    \begin{split}
       \label{2_proof_of_lemma_B_2_2}
        &{ f }_{ t }({ x }_{ t })+{ \nabla { f }_{ t }({ x }_{ t }) }^{ T }({ x }_{ t+1 }-{ x }_{ t })+[{ \lambda  }(t)+\gamma { { \bm{g} } }_{ t }({ x }_{ t }){ ] }^{ T }(\gamma { { \bm{g} } }_{ t }({ x }_{ t+1 }))+{ \alpha  }||{ x }_{ t+1 }-{ x }_{ t }{ || }^{ 2 }\\ 
        &\le { { f }_{ t }({ x }_{ t })+\nabla { f }_{ t }({ x }_{ t }) }^{ T }({ x }_{ t }^{ * }-{ x }_{ t })+{ \alpha  }||{ x }_{ t }^{ * }-{ x }_{ t }{ || }^{ 2 }-{ \alpha  }||{ x }_{ t+1 }-{ x }_{ t }^{ * }{ || }^{ 2 }\\ 
        &\le { f }_{ t }({ x }_{ t }^{ * })+{ \alpha  }||{ x }_{ t }^{ * }-{ x }_{ t }{ || }^{ 2 }-{ \alpha  }||{ x }_{ t+1 }-{ x }_{ t }^{ * }{ || }^{ 2 }
    \end{split}
\end{equation}
Rearranging terms in (\ref{2_proof_of_lemma_B_2_2}), we have
\begin{equation}
    \begin{split}
        \label{3_proof_of_lemma_B_2_2}
        &{ f }_{ t }({ x }_{ t })+[{ \bm{\lambda}  }(t)]^{ T }(\gamma { { \bm{g} } }_{ t }({ x }_{ t+1 }))\\ 
        &\le { f }_{ t }({ x }_{ t }^{ * })+{ \alpha  }||{ x }_{ t }^{ * }-{ x }_{ t }{ || }^{ 2 }-{ \alpha  }||{ x }_{ t+1 }-{ x }_{ t }^{ * }{ || }^{ 2 }-{ \alpha  }||{ x }_{ t+1 }-{ x }_{ t }{ || }^{ 2 }-{ \gamma  }^{ 2 }[{ { \bm{g} } }_{ t }({ x }_{ t }){ ] }^{ T }{ { \bm{g} } }_{ t }({ x }_{ t+1 })-{ \nabla { f }_{ t }({ x }_{ t }) }^{ T }({ x }_{ t+1 }-{ x }_{ t })\\ 
        &\overset { (a) }{ \le  } { f }_{ t }({ x }_{ t }^{ * })+{ \alpha  }||{ x }_{ t }^{ * }-{ x }_{ t }{ || }^{ 2 }-{ \alpha  }||{ x }_{ t+1 }-{ x }_{ t }^{ * }{ || }^{ 2 }-{ \alpha  }||{ x }_{ t+1 }-{ x }_{ t }{ || }^{ 2 }-{ \gamma  }^{ 2 }[{ { \bm{g} } }_{ t }({ x }_{ t }){ ] }^{ T }{ { \bm{g} } }_{ t }({ x }_{ t+1 })+||{ \nabla { f }_{ t }({ x }_{ t }) }{ || }||{ x }_{ t+1 }-{ x }_{ t }{ || }\\ 
        &\overset { (b) }{ \le  } { f }_{ t }({ x }_{ t }^{ * })+{ \alpha  }||{ x }_{ t }^{ * }-{ x }_{ t }{ || }^{ 2 }-{ \alpha  }||{ x }_{ t+1 }-{ x }_{ t }^{ * }{ || }^{ 2 }-{ \alpha  }||{ x }_{ t+1 }-{ x }_{ t }{ || }^{ 2 }-{ \gamma  }^{ 2 }[{ { \bm{g} } }_{ t }({ x }_{ t }){ ] }^{ T }{ { \bm{g} } }_{ t }({ x }_{ t+1 })\\
        &+\frac { 1 }{ 2\delta  } ||\nabla { f }_{ t }({ x }_{ t }){ || }^{ 2 }+\frac { \delta  }{ 2 } ||{ x }_{ t+1 }-{ x }_{ t }{ || }^{ 2 }\\ 
        &\overset { (c) }{ \le  } { f }_{ t }({ x }_{ t }^{ * })+{ \alpha  }||{ x }_{ t }^{ * }-{ x }_{ t }{ || }^{ 2 }-{ \alpha  }||{ x }_{ t+1 }-{ x }_{ t }^{ * }{ || }^{ 2 }-{ \alpha  }||{ x }_{ t+1 }-{ x }_{ t }{ || }^{ 2 }-{ \gamma  }^{ 2 }[{ { \bm{g} } }_{ t }({ x }_{ t }){ ] }^{ T }{ { \bm{g} } }_{ t }({ x }_{ t+1 })\\
        &+\frac { 1 }{ 2\delta  } { G }^{ 2 }+\frac { \delta  }{ 2 } ||{ x }_{ t+1 }-{ x }_{ t }{ || }^{ 2 }
    \end{split}
\end{equation}
Where (a) holds by the Cauchy-Schwarz inequality; (b) comes from the AM–GM inequality; (c) holds due to the Assumption 1.
Recall that we have the following inequality stated before,
\begin{equation}
    \begin{split}
       \label{4_proof_of_lemma_B_2_2}
        ||{ x }_{ t }^{ * }-{ x }_{ t }{ || }^{ 2 }-||{ x }_{ t+1 }{ -{ x }_{ t }^{ * }|| }^{ 2 }\le ||{ x }_{ t }^{ * }-{ x }_{ t }{ || }^{ 2 }-||{ x }_{ t+1 }{ -{ x }_{ t+1 }^{ * }|| }^{ 2 }+4R{ ||{ x }_{ t+1 }^{ * }-{ x }_{ t }^{ * }|| }
    \end{split}
\end{equation}
Furthermore, we have
\begin{equation}
    \begin{split}
        \label{5_proof_of_lemma_B_2_2}
        -[{ { \bm{g} } }_{ t }({ x }_{ t }){ ] }^{ T }{ { \bm{g} } }_{ t }({ x }_{ t+1 })&=-\frac { 1 }{ 2 } ||{ { \bm{g} } }_{ t }({ x }_{ t }){ || }^{ 2 }-\frac { 1 }{ 2 } ||{ { \bm{g} } }_{ t }({ x }_{ t+1 }){ || }^{ 2 }+\frac { 1 }{ 2 } ||{ { \bm{g} } }_{ t }({ x }_{ t })-{ { \bm{g} } }_{ t }({ x }_{ t+1 }){ || }^{ 2 }\\ 
        &\overset { (a) }{ \le  } -\frac { 1 }{ 2 } ||{ { \bm{g} } }_{ t }({ x }_{ t }){ || }^{ 2 }-\frac { 1 }{ 2 } ||{ { \bm{g} } }_{ t }({ x }_{ t+1 }){ || }^{ 2 }+\frac { 1 }{ 2 }{ \beta  }^{ 2 }||{ x }_{ t+1 }-{ x }_{ t }{ || }^{ 2 }
    \end{split}
\end{equation}
Where (a) holds by the Lipschitz continuity of ${\bm{g}}_{t}$ (Assumption 1). Substituting above two inequalities into (\ref{3_proof_of_lemma_B_2_2}) we obtain
\begin{equation}
    \begin{split}
     \label{6_proof_of_lemma_B_2_2}
        &{ f }_{ t }({ x }_{ t })+[{ { \bm{\lambda}  } }(t)]^{ T }(\gamma { { \bm{ g } } }_{ t }({ x }_{ t+1 }))\\ 
        &\le { f }_{ t }({ x }_{ t }^{ * })+{ \alpha  }_{ t }[||{ x }_{ t }^{ * }-{ x }_{ t }{ || }^{ 2 }-||{ x }_{ t+1 }{ -{ x }_{ t+1 }^{ * }|| }^{ 2 }+4R{ ||{ x }_{ t+1 }^{ * }-{ x }_{ t }^{ * }|| }]+(\frac { 1 }{ 2 } { \beta  }^{ 2 }{ \gamma  }^{ 2 }+\frac { \delta  }{ 2 } -{ \alpha })||{ x }_{ t+1 }-{ x }_{ t }{ || }^{ 2 }\\ 
        &+\frac { 1 }{ 2\delta  } { G }^{ 2 }-\frac { 1 }{ 2 } { \gamma  }^{ 2 }||{ { \bm{ g } } }_{ t }({ x }_{ t }){ || }^{ 2 }-\frac { 1 }{ 2 } { \gamma  }^{ 2 }||{ { \bm{ g } } }_{ t }({ x }_{ t+1 }){ || }^{ 2 }\\ &\overset { (a) }{ \le  } { f }_{ t }({ x }_{ t }^{ * })+{ \alpha  }||{ x }_{ t }^{ * }-{ x }_{ t }{ || }^{ 2 }-{ \alpha  }||{ x }_{ t+1 }{ -{ x }_{ t+1 }^{ * }|| }^{ 2 }+4R{ \alpha  }{ ||{ x }_{ t+1 }^{ * }-{ x }_{ t }^{ * }|| }+\frac { 1 }{ 2\delta  } { G }^{ 2 }-\frac { 1 }{ 2 } { \gamma  }^{ 2 }||{ { \bm{ g } } }_{ t }({ x }_{ t }){ || }^{ 2 }-\frac { 1 }{ 2 } { \gamma  }^{ 2 }||{ { \bm{ g } } }_{ t }({ x }_{ t+1 }){ || }^{ 2 }
    \end{split}
\end{equation}
Where (a) holds since $\alpha\ge \frac{1}{2}{\beta}^{2}{\gamma}^{2}+\frac{1}{2}\delta$. Based on Lemma \ref{lemma 6},
 adding Lyapunov drift term on both sides of (\ref{6_proof_of_lemma_B_2_2}) and rearranging terms yields:
\begin{equation}
    \begin{split}
       &{ f }_{ t }({ x }_{ t })+\Delta (t)\\ 
       &\le { f }_{ t }({ x }_{ t }^{ * })+{ \alpha  }||{ x }_{ t }^{ * }-{ x }_{ t }{ || }^{ 2 }-{ \alpha  }||{ x }_{ t+1 }{ -{ x }_{ t+1 }^{ * }|| }^{ 2 }+4R{ \alpha  }{ ||{ x }_{ t+1 }^{ * }-{ x }_{ t }^{ * }|| }+\frac { 1 }{ 2\delta  } { G }^{ 2 }+\frac { 1 }{ 2 } { \gamma  }^{ 2 }||{ { { \bm{ g } } } }_{ t+1 }({ x }_{ t+1 }){ || }^{ 2 }\\ 
       &-\frac { 1 }{ 2 } { \gamma  }^{ 2 }||{ { { \bm{ g } } } }_{ t }({ x }_{ t+1 }){ || }^{ 2 }+\frac { 1 }{ 2 } { \gamma  }^{ 2 }||{ { { \bm{ g } } } }_{ t+1 }({ x }_{ t+1 }){ || }^{ 2 }-\frac { 1 }{ 2 } { \gamma  }^{ 2 }||{ { { \bm{ g } } } }_{ t }({ x }_{ t }){ || }^{ 2 }+\gamma [{ { \bm{\lambda}  } }(t)]^{ T }({ { { \bm{ g } } } }_{ t+1 }({ x }_{ t+1 })-{ { { \bm{ g } } } }_{ t }({ x }_{ t+1 }))\\ 
       &\overset { (a) }{ \le  } { f }_{ t }({ x }_{ t }^{ * })+{ \alpha  }||{ x }_{ t }^{ * }-{ x }_{ t }{ || }^{ 2 }-{ \alpha  }||{ x }_{ t+1 }{ -{ x }_{ t+1 }^{ * }|| }^{ 2 }+4R{ \alpha  }{ ||{ x }_{ t+1 }^{ * }-{ x }_{ t }^{ * }|| }+\frac { 1 }{ 2\delta  } { G }^{ 2 }\\ 
       &+{ \gamma  }^{ 2 } F||{ { { \bm{ g } } } }_{ t+1 }({ x }_{ t+1 })-{ { { {\bm { g } } } }_{ t }({ x }_{ t+1 })|| }+\frac { 1 }{ 2 } { \gamma  }^{ 2 }||{ { { \bm{ g } } } }_{ t+1 }({ x }_{ t+1 }){ || }^{ 2 }-\frac { 1 }{ 2 } { \gamma  }^{ 2 }||{ { { \bm{ g } } } }_{ t }({ x }_{ t }){ || }^{ 2 }+\gamma [{ { \bm{\lambda}  } }(t)]^{ T }({ { { \bm{ g } } } }_{ t+1 }({ x }_{ t+1 })-{ { { \bm{ g } } } }_{ t }({ x }_{ t+1 }))\\ 
       &\overset { (b) }{ \le  } { f }_{ t }({ x }_{ t }^{ * })+{ \alpha  }||{ x }_{ t }^{ * }-{ x }_{ t }{ || }^{ 2 }-{ \alpha  }||{ x }_{ t+1 }{ -{ x }_{ t+1 }^{ * }|| }^{ 2 }+4R{ \alpha  }{ ||{ x }_{ t+1 }^{ * }-{ x }_{ t }^{ * }|| }^{ 2 }+\frac { 1 }{ 2\delta  } { G }^{ 2 }\\ 
       &+{ \gamma  }^{ 2 } F||{ { { \bm{ g } } } }_{ t+1 }({ x }_{ t+1 })-{ { { { \bm{ g } } } }_{ t }({ x }_{ t+1 })|| }+\frac { 1 }{ 2 } { \gamma  }^{ 2 }||{ { {\bm { g } } } }_{ t+1 }({ x }_{ t+1 }){ || }^{ 2 }-\frac { 1 }{ 2 } { \gamma  }^{ 2 }||{ { {\bm { g } } } }_{ t }({ x }_{ t }){ || }^{ 2 }+\gamma ||{ { \bm{\lambda}  } }(t){ || }||{ { { \bm{ g } } } }_{ t+1 }({ x }_{ t+1 })-{ { { \bm{ g } } } }_{ t }({ x }_{ t+1 }){ || }
    \end{split}
\end{equation}
Where (a) is due to $||{ { { { \bm{g} } } } }_{ t-1 }({ x }_{ t })-{ { { { \bm{g} } } } }_{ t }({ { x }_{ t } }){ || }\le 2F $; (b) holds by the Cauchy-Schwarz inequality. It completes the proof.
\subsection{proof of Lemma 9}
\emph{Proof.} According to Lemma \ref{lemma 7}, taking a telescoping sum over $t= 1, ..., T$, we obtain
\begin{equation}
    \begin{split}
        &\sum _{ t=1 }^{ T }{ { f }_{ t }({ x }_{ t }) } +\sum _{ t=1 }^{ T }{ \Delta (t) } \le \sum _{ t=1 }^{ T }{ { f }_{ t }({ x }_{ t }^{ * }) } +{ \alpha  }||{ x }_{ 1 }^{ * }-{ x }_{ 1 }{ || }^{ 2 }+4R\sum _{ t=1 }^{ T }{ { \alpha  }||{ x }_{ t+1 }^{ * }-{ x }_{ t }^{ * }{ || } } \\ 
        &+\frac { T{ G }^{ 2 } }{ 2\delta  } +{ \gamma  }^{ 2 } F{ V }_{\bm { g } }+\frac { 1 }{ 2 } { \gamma  }^{ 2 }||{ { { \bm{ g } } } }_{ T+1 }({ x }_{ T+1 }){ || }^{ 2 }-\frac { 1 }{ 2 } { \gamma  }^{ 2 }||{ { { \bm{ g } } } }_{ 1 }({ x }_{ 1 }){ || }^{ 2 }+\gamma \sum _{ t=1 }^{ T }{ ||{ { \bm{\lambda}  } }(t){ || }||{ { { \bm{ g } } } }_{ t+1 }({ x }_{ t+1 })-{ { {\bm{ g } } } }_{ t }({ x }_{ t+1 }){ || } } \\ 
        &\le \sum _{ t=1 }^{ T }{ { f }_{ t }({ x }_{ t }^{ * }) } +{ \alpha  }||{ x }_{ 1 }^{ * }-{ x }_{ 1 }{ || }^{ 2 }+4R\alpha {V}_{x} +\frac { T{ G }^{ 2 } }{ 2\delta  } +{ \gamma  }^{ 2 } F{ V }_{ \bm{ g } }+\gamma \max _{ t }{ ||{ { \bm{\lambda } } }(t){ || } } { V }_{ \bm{ g } }
    \end{split}
\end{equation}
Here we define ${\bm{g}}_{T+1}={\bm{g}}_{T}$, rearranging terms yields:
\begin{equation}
    \begin{split}
        &\sum _{ t=1 }^{ T }{ { f }_{ t }({ x }_{ t }) } \le \sum _{ t=1 }^{ T }{ { f }_{ t }({ x }_{ t }^{ * }) } +{ \alpha  }||{ x }_{ 1 }^{ * }-{ x }_{ 1 }{ || }^{ 2 }+4R\alpha {V}_{x}  +\frac { T{ G }^{ 2 } }{ 2\delta  } +{ \gamma  }^{ 2 } F{ V }_{ \bm{g} }+L(1)-L(T+1)+\gamma \max _{ t }{ ||{ { \bm{\lambda } } }(t){ || } } { V }_{ \bm{ g } }\\ 
        &\le \sum _{ t=1 }^{ T }{ { f }_{ t }({ x }_{ t }^{ * }) } +{ \alpha  }||{ x }_{ 1 }^{ * }-{ x }_{ 1 }{ || }^{ 2 }+4R\alpha {V}_{x} +\frac { T{ G }^{ 2 } }{ 2\delta  } +{ \gamma  }^{ 2 } F{ V }_{ \bm{g} }+\frac{1}{2}||\bm{\lambda}(1){||}^{2}+\gamma \max _{ t }{ ||{ {\bm{ \lambda}  } }(t){ || } } { V }_{ \bm{ g } }
    \end{split}
\end{equation}
This completes the proof.
\subsection{proof of Lemma 10}
\emph{Proof.} Recall that Lemma \ref{lemma 6} implies that ${ \gamma g }_{ t,k }({ x }_{ t })\le { \lambda  }_{ k }(t)-{ \lambda  }_{ k }(t-1),\;\forall k$. Telescoping it over $t$ yields:
\begin{equation}
    \begin{split}
        &\sum _{ t=1 }^{ T }{ { \gamma g }_{ t,k }({ x }_{ t }) } \le { \lambda  }_{ k }(T)-{ \lambda  }_{ k }(0),\; \forall k\in \{ 1,2,...,K\}. \\ 
        &\Rightarrow \sum _{ t=1 }^{ T }{ { g }_{ t,k }({ x }_{ t }) } \le \frac { { \lambda  }_{ k }(T) }{ \gamma  } \le \frac{||\bm{\lambda}(T)||}{\gamma},\; \forall k\in \{ 1,2,...,K\}. 
    \end{split}
\end{equation}
It completes the proof. 
\subsection{Proof of Lemma 11}
\textbf{Proof:} According to the strong convexity of ${ \nabla { f }_{ t }({ x }_{ t }) }^{ T }(x-{ x }_{ t })+[{ \bm{\lambda} }(t)+\gamma { { \bm{g} } }_{ t }({ x }_{ t }){ ] }^{ T }(\gamma { { \bm{g} } }_{ t }(x))+{ \alpha  }||x-{ x }_{ t }{ || }^{ 2 }$ with respect to \emph{x} and recalling that ${x}_{t+1}$ minimizes this expression over $\chi$, we have
\begin{equation}
    \begin{split}
       &{ \nabla { f }_{ t }({ x }_{ t }) }^{ T }({ x }_{ t+1 }-{ x }_{ t })+[{ \bm{\lambda}  }(t)+\gamma { { { \bm{g} } } }_{ t }({ x }_{ t }){ ] }^{ T }(\gamma { { { \bm{g} } } }_{ t }({ x }_{ t+1 }))+{ \alpha  }||{ x }_{ t+1 }-{ x }_{ t }{ || }^{ 2 }\\ 
      & \le { \nabla { f }_{ t }({ x }_{ t }) }^{ T }(\hat { x } -{ x }_{ t })+[{ \bm{\lambda}  }(t)+\gamma { { { \bm{g} } } }_{ t }({ x }_{ t }){ ] }^{ T }(\gamma { { { \bm{g} } } }_{ t }(\hat { x } ))+{ \alpha  }||\hat { x } -{ x }_{ t }{ || }^{ 2 }-{ \alpha  }||{ x }_{ t+1 }-\hat { x } { || }^{ 2 }\\ 
       &\overset { (a) }{ \le  } { \nabla { f }_{ t }({ x }_{ t }) }^{ T }(\hat { x } -{ x }_{ t })-\gamma \epsilon \sum _{ k=1 }^{ K }{ { [\lambda  }_{ k }(t)+\gamma { g }_{ t,k }({ x }_{ t })] } +{ \alpha  }||\hat { x } -{ x }_{ t }{ || }^{ 2 }-{ \alpha  }||{ x }_{ t+1 }-\hat { x } { || }^{ 2 }\\ 
       &\overset { (b) }{ = } { \nabla { f }_{ t }({ x }_{ t }) }^{ T }(\hat { x } -{ x }_{ t })-\gamma \epsilon ||{ \bm{\lambda} }(t)+\gamma { { { \bm{g} } } }_{ t }({ x }_{ t }){ || }_{ 1 }+{ \alpha  }||\hat { x } -{ x }_{ t }{ || }^{ 2 }-{ \alpha  }||{ x }_{ t+1 }-\hat { x } { || }^{ 2 }\\ 
       &\overset { (c) }{ \le  } { \nabla { f }_{ t }({ x }_{ t }) }^{ T }(\hat { x } -{ x }_{ t })-\gamma \epsilon ||{ \bm{\lambda} }(t)+\gamma { { { \bm{g} } } }_{ t }({ x }_{ t }){ || }+{ \alpha  }||\hat { x } -{ x }_{ t }{ || }^{ 2 }-{ \alpha  }||{ x }_{ t+1 }-\hat { x } { || }^{ 2 }\\ 
       &\overset { (d) }{ \le  } { \nabla { f }_{ t }({ x }_{ t }) }^{ T }(\hat { x } -{ x }_{ t })-\gamma \epsilon [||{ \bm{\lambda}  }(t){ || }-||\gamma { { { \bm{g} } } }_{ t }({ x }_{ t }){ || }]+{ \alpha  }||\hat { x } -{ x }_{ t }{ || }^{ 2 }-{ \alpha  }||{ x }_{ t+1 }-\hat { x } { || }^{ 2 }
    \end{split}
\end{equation}
Where (a) is due to the Slater condition (Assumption 2); (b) holds since ${\lambda}_{k}(t)+\gamma{g}_{t,k}({x}_{t})\ge 0,\;\forall k$; (c) holds due to the fact that $||x{ || }_{ 1 }\ge ||x{ || }$ for any vector $x\in \chi$; (d) holds by the triangle inequality $||u-v{ || }\ge ||u{ || }-||v{ || },\; \forall u,v\in \chi$. Base on Lemma \ref{lemma 6}, we add Lyapunov drift term $\Delta(t)$ on both sides and rearranging some terms yields:
\begin{equation}
    \begin{split}
       \label{1_of_lemma_B6}
       &\Delta (t)\le { \nabla { f }_{ t }({ x }_{ t }) }^{ T }(\hat { x } -{ x }_{ t })-{ \nabla { f }_{ t }({ x }_{ t }) }^{ T }({ x }_{ t+1 }-{ x }_{ t })-{ \gamma  }^{ 2 }{ { { { { \bm{g} } } } } }_{ t }({ x }_{ t })^{ T }{ { { { { \bm{g} } } } } }_{ t }({ x }_{ t+1 })-\gamma \epsilon [||{ { \bm{\lambda}  } }(t){ || }-||\gamma { { { { { \bm{g} } } } } }_{ t }({ x }_{ t }){ || }]\\ 
      & +{ \alpha  }||\hat { x } -{ x }_{ t }{ || }^{ 2 }-{ \alpha  }||{ x }_{ t+1 }-\hat { x } { || }^{ 2 }-{ \alpha  }||{ x }_{ t+1 }-{ x }_{ t }{ || }^{ 2 }+{ \gamma  }^{ 2 }||{ { \bm{g} } }_{ t+1 }({ x }_{ t+1 }){ || }^{ 2 }+\gamma [\bm{\lambda} (t){ ] }^{ T }({ { \bm{g} } }_{ t+1 }({ x }_{ t+1 })-{ { \bm{g} } }_{ t }({ x }_{ t+1 }))\\ 
       &\le { \nabla { f }_{ t }({ x }_{ t }) }^{ T }(\hat { x } -{ x }_{ t+1 })-\gamma \epsilon ||{ { \bm{\lambda}  } }(t){ || }+{ \gamma  }^{ 2 }\epsilon ||{ { { { { \bm{g} } } } } }_{ t }({ x }_{ t }){ || }+{ \gamma  }^{ 2 }||{ { \bm{g} } }_{ t+1 }({ x }_{ t+1 }){ || }^{ 2 }-{ \gamma  }^{ 2 }{ { { { { \bm{g} } } } } }_{ t }({ x }_{ t })^{ T }{ { { { { \bm{g} } } } } }_{ t }({ x }_{ t+1 })\\
      & +{ \alpha  }||\hat { x } -{ x }_{ t }{ || }^{ 2 }-{ \alpha  }||{ x }_{ t+1 }-\hat { x } { || }^{ 2 }-{ \alpha  }||{ x }_{ t+1 }-{ x }_{ t }{ || }^{ 2 }+\gamma [\bm{\lambda} (t){ ] }^{ T }({ { \bm{g} } }_{ t+1 }({ x }_{ t+1 })-{ { \bm{g} } }_{ t }({ x }_{ t+1 }))\\ 
      & \le { \nabla { f }_{ t }({ x }_{ t }) }^{ T }(\hat { x } -{ x }_{ t+1 })-\gamma \epsilon ||{ { \bm{\lambda} } }(t){ || }+{ \gamma  }^{ 2 }\epsilon ||{ { { { { \bm{g} } } } } }_{ t }({ x }_{ t }){ || }+{ \gamma  }^{ 2 }||{ { \bm{g} } }_{ t+1 }({ x }_{ t+1 }){ || }^{ 2 }-{ \gamma  }^{ 2 }{ { { { { \bm{g} } } } } }_{ t }({ x }_{ t })^{ T }{ { { { { \bm{g} } } } } }_{ t }({ x }_{ t+1 })\\ 
       &+{ \alpha  }||\hat { x } -{ x }_{ t }{ || }^{ 2 }+\gamma [\bm{\lambda} (t){ ] }^{ T }({ { \bm{g} } }_{ t+1 }({ x }_{ t+1 })-{ { \bm{g} } }_{ t }({ x }_{ t+1 }))\\ 
       &\overset { (a) }{ \le  } { \nabla { f }_{ t }({ x }_{ t }) }^{ T }(\hat { x } -{ x }_{ t+1 })-\gamma \epsilon ||{ { \bm{\lambda}  } }(t){ || }+{ \gamma  }^{ 2 }\epsilon  F+{ \gamma  }^{ 2 }{ F }^{ 2 }-{ \gamma  }^{ 2 }{ { { { { \bm{g} } } } } }_{ t }({ x }_{ t })^{ T }{ { { { { \bm{g} } } } } }_{ t }({ x }_{ t+1 })\\ 
       &+{ \alpha  }||\hat { x } -{ x }_{ t }{ || }^{ 2 }+\gamma [\bm{\lambda} (t){ ] }^{ T }({ { \bm{g} } }_{ t+1 }({ x }_{ t+1 })-{ { \bm{g} } }_{ t }({ x }_{ t+1 }))\\
       &\overset { (b) }{ \le  } ||\nabla { f }_{ t }({ x }_{ t }){ || }||\hat { x } -{ x }_{ t+1 }{ || }-\gamma \epsilon ||{ { \bm{\lambda}  } }(t){ || }+{ \gamma  }^{ 2 }\epsilon  F+{ \gamma  }^{ 2 }{ F }^{ 2 }+{ \gamma  }^{ 2 }||{ { { { { \bm{g} } } } } }_{ t-1 }({ x }_{ t }){ || }||{ { { { { \bm{g} } } } } }_{ t }({ x }_{ t+1 }){ || }\\ 
       &+{ \alpha  }||\hat { x } -{ x }_{ t }{ || }^{ 2 }+\gamma ||\bm{\lambda} (t){ || }||{ { g } }_{ t+1 }({ x }_{ t+1 })-{ { g } }_{ t }({ x }_{ t+1 }){ || }
    \end{split}
\end{equation}
\begin{equation}
    \begin{split}
       &\overset { (c) }{ \le  } -\gamma \epsilon ||{ { \bm{\lambda}  } }(t){ || }+GR+{ \gamma  }^{ 2 }\epsilon F+{ \gamma  }^{ 2 }{ F }^{ 2 }+{ \gamma  }^{ 2 }{ F }^{ 2 }+\alpha { R }^{ 2 }+\gamma ||\bm{\lambda} (t){ || }||{ { \bm{g} } }_{ t+1 }({ x }_{ t+1 })-{ { \bm{g} } }_{ t }({ x }_{ t+1 }){ || }\\ 
      & =-\gamma (\epsilon -||{ { \bm{g} } }_{ t+1 }({ x }_{ t+1 })-{ { \bm{g} } }_{ t }({ x }_{ t+1 }){ || })||{ { \bm{\lambda}  } }(t){ || }+GR+{ \gamma  }^{ 2 }\epsilon  F+2{ \gamma  }^{ 2 }{ F }^{ 2 }+{\alpha} { R }^{ 2 }\\ 
      &\le -\gamma (\epsilon -\max _{ t }{ \max _{ x\in \chi  }{ \{ ||{ { \bm{g} } }_{ t+1 }({ x })-{ { \bm{g} } }_{ t }({ x }){ || }\}  }  } )||{ { \bm{\lambda}  } }(t){ || }+GR+{ \gamma  }^{ 2 }\epsilon  F+2{ \gamma  }^{ 2 }{ F }^{ 2 }+{\alpha} { R }^{ 2 }
    \end{split}
\end{equation}
Where (a) holds by Assumption 1; (b) is due to the Cauchy–Schwarz inequality; (c) follows from the Assumption 1. We have $\hat{\epsilon}=\epsilon -\max _{ t }{ \max _{ x\in \chi  }{ \{ ||{ { \bm{g} } }_{ t+1 }({ x })-{ { \bm{g} } }_{ t }({ x }){ || }\}  }  }>0$ since $ \epsilon >{ \bar { V }  }_{ g }=\max _{ t }{ \max _{ x\in \chi  }{ \{ ||{ \bm{g} }_{ t+1 }(x)-{ \bm{g} }_{ t }(x){ || }\}  }  } $,
Next we perform a reduction to absurdity process to prove this Lemma. Recall that $ { \lambda  }_{ k }(1)=0$ for all $k \in \{1,2,...,K\}$ implies that $||{ \bm{\lambda}  }(1){ || }=0\le \gamma  F+\frac { GR+2{ \gamma  }^{ 2 }{ F }^{ 2 }+{\alpha} { R }^{ 2 }+{ \gamma  }^{ 2 }\epsilon  F }{ \gamma \hat{\epsilon}  }  $. Assume that there exists a $\tau \ge 2$ such that $||{ \bm{\lambda}  }(\tau ){ || }>\gamma  F+\frac { GR+2{ \gamma  }^{ 2 }{ F }^{ 2 }+{\alpha} { R }^{ 2 }+{ \gamma  }^{ 2 }\epsilon  F }{ \gamma \hat{\epsilon}  } $ and $||{ \bm{\lambda}  }(t){ || }\le \gamma F+\frac { GR+2{ \gamma  }^{ 2 }{ F }^{ 2 }+{\alpha} { R }^{ 2 }+{ \gamma  }^{ 2 }\epsilon  F }{ \gamma \hat{\epsilon}  } $ for all $t < \tau$, then we consider two cases about the value of $||{ \bm{\lambda}  }(\tau -1){ || }$. ~\\
\textbf{Case 1:} If $||{ \bm{\lambda}  }(\tau -1){ || }>\frac { GR+2{ \gamma  }^{ 2 }{ F }^{ 2 }+{\alpha} { R }^{ 2 }+{ \gamma  }^{ 2 }\epsilon  F  }{ \gamma \hat{\epsilon}  } $, then we can derive $\Delta (\tau-1)<0$. According to (\ref{1_of_lemma_B6}), we have
\begin{equation}
     \begin{split}
         &||{ { \lambda  } }(\tau ){ || }<||{ { \lambda  } }(\tau -1){ || }\le \frac { GR+2{ \gamma  }^{ 2 }{ F }^{ 2 }+{ \alpha  }{ R }^{ 2 }+{ \gamma  }^{ 2 }\epsilon  F }{ \gamma \hat{\epsilon}  } \\ 
         &<\gamma F+\frac { GR+2{ \gamma  }^{ 2 }{ F }^{ 2 }+{ \alpha  }{ R }^{ 2 }+{ \gamma  }^{ 2 }\epsilon  F }{ \gamma \hat{\epsilon}  } 
     \end{split}
\end{equation}
Which contradicts the definition of $\tau$. ~\\

\textbf{Case 2:} If $||{ \bm{\lambda}  }(\tau -1){ || }\le \frac { GR+2{ \gamma  }^{ 2 }{ F }^{ 2 }+{\alpha}{ R }^{ 2 }+{ \gamma  }^{ 2 }\epsilon  F }{ \gamma \hat{\epsilon}  }$, then according to Lemma \ref{lemma 6} we have
\begin{equation}
    \begin{split}
        &||{ { \lambda  } }(\tau ){ || }\le ||{ { \lambda  } }(\tau -1){ || }+\gamma ||{ { g } }_{ t-1 }({ x }_{ t }){ || }\le \gamma  F+\frac { GR+2{ \gamma  }^{ 2 }{ F }^{ 2 }+{ \alpha  }{ R }^{ 2 }+{ \gamma  }^{ 2 }\epsilon  F }{ \gamma \hat{\epsilon}  }
    \end{split}
\end{equation}
Which also contradicts the definition of $\tau$.
Hence $||{ \bm{\lambda}  }(t){ || }\le \gamma  F+\frac { GR+2{ \gamma  }^{ 2 }{ F }^{ 2 }+{\alpha} { R }^{ 2 }+{ \gamma  }^{ 2 }\epsilon  F }{ \gamma \hat{\epsilon}  } $ holds for all $t>1$. It completes the Proof.
\end{document}